\documentclass[10pt]{article}
\usepackage{amsfonts}
\usepackage{amssymb,amsmath,amsthm, amsfonts}

\def\sqr#1#2{{\vcenter{\vbox{\hrule height.#2pt
              \hbox{\vrule width.#2pt height#1pt \kern#1pt \vrule width.#2pt}
          \hrule height.#2pt}}}}
\def\signed #1{{\unskip\nobreak\hfil\penalty50
          \hskip2em\hbox{}\nobreak\hfil#1
          \parfillskip=0pt \finalhyphendemerits=0 \par}}
\def\endpf{\signed {$\sqr69$}}
\def\sqr#1#2{{\vcenter{\vbox{\hrule height.#2pt
              \hbox{\vrule width.#2pt height#1pt \kern#1pt \vrule width.#2pt}
              \hrule height.#2pt}}}}
\def\signed #1{{\unskip\nobreak\hfil\penalty50
              \hskip2em\hbox{}\nobreak\hfil#1
              \parfillskip=0pt \finalhyphendemerits=0 \par}}
\def\endpf{\signed {$\sqr69$}}
\def\3n{\negthinspace \negthinspace \negthinspace }
\def\2n{\negthinspace \negthinspace }
\def\1n{\negthinspace }

\def\={\buildrel \triangle \over =}

\def\esssup{\mathop{\rm esssup}}

\def\exp{\mathop{\rm exp}}
\def\sup{\mathop{\rm sup}}

\def\inf{\hbox{\rm inf$\,$}}
\def\esssup{\hbox{\rm ess$\,$\rm sup$\,$}}

\def\tr{\hbox{\rm tr$\,$}}

\def\|{\Big |}
\def\({\Big (}
\def\){\Big )}
\def\[{\Big[}
\def\]{\Big]}
\def\be{\begin{equation}}
\def\bel{\begin{equation}\label}
\def\ee{\end{equation}}
\def\bt{\begin{theorem}}
\def\bcd{\begin{condition}}
\def\ecd{\end{condition}}
\def\et{\end{theorem}}
\def\bc{\begin{corollary}}
\def\ec{\end{corollary}}
\def\bde{\begin{definition}}
\def\ede{\end{definition}}
\def\bl{\begin{lemma}}
\def\el{\end{lemma}}
\def\bp{\begin{proposition}}
\def\ep{\end{proposition}}
\def\bex{\begin{example}}
\def\eex{\end{example}}
\def\br{\begin{remark}}
\def\er{\end{remark}}
\def\ba{\begin{array}}
\def\ea{\end{array}}
\def\ed{\end{document}}

\def\square#1{\vbox{\hrule\hbox{\vrule height#1%
     \kern#1\vrule}\hrule}}
\def\rectangle#1#2{\vbox{\hrule\hbox{\vrule height#1%
     \kern#2\vrule}\hrule}}

\font\tenbb=msbm10 \font\sevenbb=msbm7 \font\fivebb=msbm5

\newfam\bbfam
\scriptscriptfont\bbfam=\fivebb \textfont\bbfam=\tenbb
\scriptfont\bbfam=\sevenbb

\newtheorem{lemma}{Lemma}[section]
\newtheorem{remark}{Remark}[section]
\newtheorem{example}{Example}[section]
\newtheorem{theorem}{Theorem}[section]
\newtheorem{corollary}{Corollary}[section]

\newtheorem{definition}{Definition}[section]
\newtheorem{proposition}{Proposition}[section]
\newtheorem{condition}{Condition}[section]

\makeatletter
   
   \@addtoreset{equation}{section}
\makeatother

\begin{document}

\title{ Optimal control problems of fully coupled FBSDEs and viscosity solutions of Hamilton-Jacobi-Bellman equations}
\author{Juan Li\\
{\small School of Mathematics and Statistics, Shandong University, Weihai, Weihai 264209, P. R. China.}\\
{\small{\it E-mail:  juanli@sdu.edu.cn}}\\
Qingmeng Wei{\footnote{Corresponding author.}}\\
{\small School of Mathematics, Shandong University, Jinan 250100, P. R. China.}\\
{\small {\it E-mail: weiqingmeng0207@163.com}}}
\date{November 11, 2012}

\maketitle \noindent{\bf Abstract.}

\hskip4mm  In this paper we study stochastic optimal control
problems of fully coupled forward-backward stochastic differential
equations (FBSDEs). The recursive cost functionals are defined by
controlled fully coupled FBSDEs. We study two cases of diffusion
coefficients $\sigma$ of FSDEs. We use a new method to prove that
the value functions are deterministic, satisfy the dynamic
programming principle (DPP), and are viscosity solutions to the
associated generalized Hamilton-Jacobi-Bellman (HJB) equations. The
associated generalized HJB equations are related with algebraic
equations when $\sigma$ depends on the second component of the
solution $(Y, Z)$ of the BSDE and doesn't depend on the control. For
this we adopt Peng's BSDE method, and so in particular, the notion
of stochastic backward semigroup in~\cite{Pe4}. We emphasize that
the fact that $\sigma$ also depends on $Z$ makes the stochastic
control much more complicate and has as consequence that the
associated HJB equation is combined with an algebraic equation,
which is inspired by Wu and Yu~\cite{WY}. We use the continuation method combined with the fixed point theorem to
prove that the algebraic equation has a unique
solution, and moreover, we also give the representation for this
solution. On the other hand, we prove
 some new basic estimates for fully coupled FBSDEs under the monotonic assumptions. In particular,
 we prove under the Lipschitz and linear growth conditions that fully coupled FBSDEs have a unique solution on the small time interval,
 if the Lipschitz constant of $\sigma$\ with respect to $z$ is sufficiently small.
  We also establish a generalized comparison theorem for such fully coupled FBSDEs.

\noindent{{\bf Keywords.}\small \ \ Fully coupled FBSDEs; value functions; stochastic backward semigroup; dynamic programming principle; viscosity solution} \\

\newpage

\section{\large{Introduction}}

\hskip1cm Nonlinear backward stochastic differential equations
(BSDEs) driven by a Brownian motion were first introduced by Pardoux
and Peng~\cite{PaPe1} in 1990. They got the uniqueness and the
existence theorem for nonlinear BSDEs under Lipschitz condition. The
theory of BSDEs has been studied since then by many authors and has
found various applications, namely in stochastic control (see
Peng~\cite{Pe1}), finance (see El Karoui, Peng and
Quenez~\cite{ELPeQu}), and partial differential equations (PDE)
theory (see Peng~\cite{Pe2}, etc).

Related with the BSDE theory, the theory of fully coupled
forward-backward stochastic differential equation (FBSDE) has been
developing very dynamically. We will usually meet fully coupled
FBSDEs which are used to describe state processes and the related
cost functionals when some optimization problems are studied (see
Cvitani\'{c} and Ma~\cite{CM}, Ma and Yong~\cite{MY}). There are
many results on the existence and the uniqueness of solutions of
fully coupled FBSDEs. Antonelli~\cite{An} first studied fully
coupled FBSDEs driven by Brownian motion on a ``small" time interval
with the fixed point theorem. As we know, there are mainly three
methods to study fully coupled FBSDEs on an arbitrarily given time
interval. The first one is a kind of ``four-step scheme" approach
(see Ma, Protter and Yong~\cite{MPY}) which combines PDE methods and
methods of probability. With these methods, the authors
of~\cite{MPY} proved the existence and the uniqueness for fully
coupled FBSDEs on an arbitrarily given time interval, but they
required  the equations to be non-degenerate, i.e., the diffusion
coefficients are non-degenerate. However, the PDE approach can not
be used to deal with the case, when the coefficients are random. The
second method is that of continuation which is  purely
probabilistic, see Hu and Peng~\cite{HP}, Pardoux and
Tang~\cite{PaT}, Peng and Wu~\cite{PW}, Yong~\cite{Y1}. They used
the ``monotonicity" condition on the coefficients which relaxes the
above assumptions. The third method is motivated by the numerical
approaches for some linear FBSDEs (see Delarue~\cite{D} and
Zhang~\cite{Z}). Using a probabilistic method, Wu~\cite{W} proved a
comparison theorem for FBSDEs. It is a useful tool to study fully
coupled FBSDEs. Recently, Ma, Wu, Zhang and Zhang~\cite{MWZZ} used a unified method to study fully coupled FBSDEs.
For more details on fully coupled FBSDEs, the reader
is referred to the book of Ma and Yong~\cite{MY}; more recent works
on FBSDEs refer to Yong~\cite{Y2}, or Ma, Wu, Zhang and Zhang~\cite{MWZZ}, and the references therein.

BSDE methods (for instance, a generalized dynamic programming
principle (DPP) or the maximum principle for control problem),
originally developed by Peng~\cite{Pe1},~\cite{Pe3},~\cite{Pe4}, for
the stochastic control theory. Pardoux and Tang~\cite{PaT}
associated fully coupled FBSDEs (without controls) with
quasilinear parabolic PDEs, and proved the existence of
viscosity solutions. The diffusion coefficient $\sigma$ of
the forward equation in the FBSDE they considered is supposed not to
depend on $Z$. In Wu and Yu~\cite{WY},~\cite{WY2}, they proved the
existence of a viscosity solution for quasilinear PDEs with the help
of fully coupled FBSDE when $\sigma$ depends on $z$, and the
stochastic systems without controls. Inspired by above works we want to study the optimal control problems of fully coupled FBSDEs. If one considers controlled
fully coupled FBSDEs the problem becomes more difficult, as well as
when one tries to give a probabilistic representation for systems of
generalized fully nonlinear Hamilton-Jacobi-Bellman (HJB) equations.
Using a method of Buckdahn
and Li~\cite{BL}, without assuming the coefficients to be
H\"{o}lder-continuous with respect to the control variable, we prove that the value function is deterministic
(Proposition 3.1) and it satisfies DPP (Theorem 3.1). Furthermore,
we obtain the existence of viscosity solutions of the associated HJB
equations (Theorem 4.1 and Theorem 4.2). For this we adopt Peng's
BSDE method (see \cite{Pe4}, or \cite{BL}). However, unlike the
existing literature, in the present work the stochastic backward
semigroup is defined through a fully coupled FBSDE. This makes
techniques which are implied by the work with the stochastic
backward semigroup much more subtle. This arises, in particular, in
the proofs of the Theorems 4.1 and 4.2 showing that the value
function is a viscosity solution: a ``classical" approach would lead
to BSDEs with quadratic growth in $(y,z)$. To avoid this, a deeper
study of properties of fully coupled FBSDEs on short time intervals
has to be done.  Hence we prove the existence and uniqueness for
such FBSDEs in the case of a sufficiently small Lipschitz constant
of $\sigma$ with respect to $z$ (see Proposition \ref{App-Pro4}), we
give $L^p$-estimates for the solution (see Proposition
\ref{App-Pro5}) and we also establish a new general comparison
result for such fully coupled FBSDEs (Theorem \ref{App-Th2}). We
also emphasize that when $\sigma$ depends on $z$ it makes the
stochastic control much more complicate and the associated HJB
equation is combined with an algebraic equation, which is inspired
by Wu and Yu~\cite{WY}. We use the continuation method combined with
the fixed point theorem to prove that the algebraic
equation has a unique solution whose representation is given
(Proposition 4.1).

Let us be more precise. We study a stochastic control problem of
fully coupled FBSDE. The cost functional is introduced by the
following fully coupled FBSDE: \be\label{ee1.1}
  \left \{
  \begin{array}{ll}
  dX_s^{t,x;u}  =  b(s,X_s^{t,x;u},Y_s^{t,x;u},Z_s^{t,x;u},u_s)ds +
          \sigma(s,X_s^{t,x;u},Y_s^{t,x;u},Z_s^{t,x;u},u_s) dB_s, \\
  dY_s^{t,x;u}  =  -f(s,X_s^{t,x;u},Y_s^{t,x;u},Z_s^{t,x;u},u_s)ds + Z_s^{t,x;u}dB_s, \ \ \ \ \ s\in[t,T], \\
   X_t^{t,x;u} =  x,\ Y_T^{t,x;u}  =  \Phi(X_T^{t,x;u}),
   \end{array}
   \right.
  \ee
\\where $T>0$ is an arbitrarily fixed finite time horizon,
$B=(B_s)_{s\in[0,T]}$ is a $d$-dimensional standard Brownian motion,
and $u=(u_s)_{s\in[t,T]}$ is an admissible control. Precise
assumptions on the coefficients $b,\  \sigma, \ f,  \  \Phi$ are
given in the next section. Under
 our assumptions, (\ref{ee1.1}) has a unique solution
  $(X_s^{t,x;u},Y_s^{t,x;u},Z_s^{t,x;u})_{s\in[t,T]}$ and the cost
  functional is defined by \be\label{ee1.2} J(t,x;u)=Y_t^{t,x;u}. \ee

\noindent We define the value function of our stochastic control
problems as follows: \be\label{ee1.3} W(t,x):=\esssup_{u\in\mathcal
{U}_{t,T}}J(t,x;u). \ee The objective of our paper is to investigate
this value function. The main results of the paper state that $W$
is, deterministic (Proposition 3.1), continuous viscosity solution
of the associated HJB equations (Theorem 4.1 and Theorem 4.2). The
associated HJB equations are very complicated, we consider
two cases of $\sigma$ for the existence of a viscosity solution.

\textbf{Case 1.}  $\sigma$ does not depend on $z$, but depends on
$u$.

The associated HJB equation is then the following:\be\label{ee1.4}
 \left \{\begin{array}{ll}
 &\!\!\!\!\! \frac{\partial }{\partial t} W(t,x) +  H(t, x, W(t,x), DW(t,x), D^2W(t,x))=0,
 \hskip 0.5cm   (t,x)\in [0,T)\times {\mathbb{R}}^n ,  \\
 &\!\!\!\!\!  W(T,x) =\Phi (x), \hskip0.5cm   x \in {\mathbb{R}}^n,
 \end{array}\right.
\ee with
$$\begin{array}{llll}
H(t, x, y, p, X)=\sup\limits_{u \in U}\{p.b(t, x, y,  p.\sigma,
u)+\frac{1}{2}\tr(\sigma\sigma^{T}(t, x, y, u)X)+f(t, x, y,
p.\sigma, u)\},
 \end{array}$$
where $t\in [0, T],\ x\in{\mathbb{R}}^n,\ y\in \mathbb{R},\ p\in
{\mathbb{R}}^n$, and $X\in\mathbb{S}^n$ ($\mathbb{S}^n$
denotes the set of $n \times n$ symmetric matrices).

\textbf{Case 2.}  $\sigma$ does not depend on $u$, but depends on
$z$.

This case is more complicate than the former one. The associated HJB
equation is combined with an algebraic equation as follows:
\be\label{ee1.5}
 \left \{\begin{array}{ll}
 &\!\!\!\!\! \frac{\partial }{\partial t} W(t,x) +  H(t, x, W(t,x),
 V(t,x))=0,
   \\
&\!\!\!\!\!V(t,x)=DW(t,x).\sigma(t,x,W(t,x),V(t,x)),\hskip 0.5cm   (t,x)\in [0,T)\times {\mathbb{R}}^n ,\\
 &\!\!\!\!\!  W(T,x) =\Phi (x),\hskip 0.5cm   x\in {\mathbb{R}}^n.
 \end{array}\right.
\ee
 In this case
$$\begin{array}{lll}
&H(t, x, W(t,x), V(t,x))= \sup\limits_{u \in U}\{DW(t,x).b(t, x,
W(t,x), V(t,x), u)\\
 &\qquad\qquad\qquad\qquad\qquad\quad\ +\frac{1}{2}\tr(\sigma\sigma^{T}(t, x, W(t,x),
V(t,x))D^2W(t,x))
  +f(t, x, W(t,x), V(t,x), u)\},
 \end{array}$$
where $t\in [0, T], x\in{\mathbb{R}}^n.$

The Case 2 is more complicate, the associated HJB equation is
combined with an algebraic equation, which is inspired by Wu and
Yu~\cite{WY},~\cite{WY2}, we use a new method-the
continuation method combined with the fixed point theorem in order
to prove for the first time that the algebraic equation
has a unique solution, and give the representation for the solution
(see Proposition 4.1) which makes that other proofs are available.
But both cases require new estimates and new generalized comparison
theorem for small time interval FBSDEs, which are discussed in the
Appendix.

Our paper is organized as follows: Section 2 recalls some elements
of the theory of fully coupled FBSDEs which will be used
in what follows. Section 3 introduces the setting of the stochastic
control problems. We prove that the value function $W$ is a
deterministic function (Proposition 3.1) which is Lipschitz in $x$
(Lemma 3.2), monotonic (Lemma 3.3) and continuous in $t$ (Theorem
3.2). Moreover, it satisfies the DPP (Theorem 3.1). In Section 4, by
using the DPP we prove that $W$ is a viscosity solution of the
associated HJB equation in the two
 cases (Theorem 4.1 and Theorem 4.2) described above. In Section 5, we give two examples.
 In Appendix we prove some basic important estimates for fully coupled FBSDEs under the monotonic assumptions,
  and new estimates and new generalized comparison theorem for small time interval FBSDEs.

\section{ {\large Preliminaries}}
  \hskip1cm
Let $(\Omega, {\cal{F}},P)$ be the  Wiener space, where $\Omega$ is
the set of continuous functions from $[0, T]$ to ${\mathbb{R}}^d$
starting from $0$ ($\Omega= C_0([0, T];{\mathbb{R}}^d)$), $
{\cal{F}} $ the completed Borel $\sigma$-algebra over $\Omega$, and
$P$ the Wiener measure. Let $B$ be the canonical process:
$B_s(\omega)=\omega_s,\ s\in [0, T],\ \omega\in \Omega$. We denote
by ${\mathbb{F}}=\{{\mathcal{F}}_s,\ 0\leq s \leq T\}$\ the natural
filtration generated by $\{B_t\}_{t\geq0}$ and augmented by all
$P\mbox{-}$null sets, i.e., $$ {\cal{F}}_s=\sigma\{B_r,r\leq s\}\vee
\mathcal {N}_P,\ \ \ \ s\in [0,T],$$ where $\mathcal {N}_P$ is the
set of all $P\mbox{-}$null subsets and $T$ is a fixed real time
horizon. We introduce the following two spaces of processes which will be used frequently: for $t_0\in [0,T],$

\noindent ${\cal{S}}^2(t_0, T; {\mathbb{R}^n})$\ is the set of $\mathbb{R}^n\mbox{-valued}\ \mathbb{F}\mbox{-adapted
    continuous process}\ (\psi_t)_{t_0\leq t\leq T}\ \mbox{with}\ E[\sup\limits_{t_0\leq t\leq T}| \psi_{t} |^2]< +\infty$;\\

\noindent ${\cal{H}}^{2}(t_0,T;{\mathbb{R}}^{n})$\ is the set of ${\mathbb{R}}^{n}\mbox{-valued}\ \mathbb{F}
    \mbox{-progressively meas. process}\ (\psi_t)_{t_0\leq t\leq T}\  \mbox{with}\
     E[\int^T_{t_0}| \psi_t| ^2dt]<+\infty$.

 Now we introduce the following  fully coupled FBSDE associated with $(b, \sigma,f, \zeta, \Phi)$
  \be
  \left \{
  \begin{array}{ll}
   dX_s  =  b(s,X_s,Y_s,Z_s) ds +
          \sigma(s,X_s,Y_s,Z_s) dB_s, \\
   dY_s  = -f(s,X_s,Y_s,Z_s) ds +Z_s dB_s,\ \ \   s\in[t,T],\\
    X_t  =  \zeta,\ Y_T  =  \Phi(X_T),
   \end{array}
   \right.
  \ee
  where $(X,Y,Z) \in \mathbb{R}^n \times  \mathbb{R}^m \times  \mathbb{R}^{m\times d},\
  T>0,$\\
$b: \Omega \times [0,T] \times \mathbb{R}^n \times  \mathbb{R}^m
\times  \mathbb{R}^{m\times
  d}
  \rightarrow \mathbb{R}^n, \ \ \ \ \ \ \ \ \ \  \sigma: \Omega \times [0,T] \times \mathbb{R}^n \times  \mathbb{R}^m \times  \mathbb{R}^{m\times
  d}
  \rightarrow \mathbb{R}^{n\times d}, $\\
$f: \Omega \times [0,T] \times \mathbb{R}^n \times \mathbb{R}^m
\times  \mathbb{R}^{m\times
  d}
  \rightarrow \mathbb{R}^m,  \ \ \ \ \ \  \ \ \   \Phi: \Omega \times \mathbb{R}^n
  \rightarrow \mathbb{R}^m, $\\
$(b(t,x,y,z))_{t\in[0,T]},\ (\sigma(t,x,y,z))_{t\in[0,T]},\ (f(t,x,y,z))_{t\in[0,T]}$ are $\mathbb {F}$-progressively measurable for each $(x,y,z)\in \mathbb{R}^n \times  \mathbb{R}^m \times  \mathbb{R}^{m\times d}$, and  $\Phi(x)$\ is ${\cal F}_T$-measurable for each $x\in \mathbb{R}^n$. In this paper
we use the usual inner product and the Euclidean norm in
$\mathbb{R}^n, \ \mathbb{R}^m$ and $\mathbb{R}^{m\times d},$
respectively. Given an $m \times n$ full-rank matrix $G$, we define:
$$\lambda= \
            \left(\begin{array}{c}
            x\\
            y\\
            z
            \end{array}\right)
            \ , \ \ \ \ \ \ \ \ \ \
A(t,\lambda)= \
            \left(\begin{array}{c}
            -G^Tf\\
            Gb\\
            G\sigma
            \end{array}\right)(t,\lambda),$$where $G^T$ is the transposed matrix of $G$.

We assume that

\noindent (B1) (i) $A(t,\lambda)$ is uniformly Lipschitz  with respect to
$\lambda$, and for any $\lambda$, $A(\cdot,\lambda)\in {\cal{H}}^{2}(0,T;{\mathbb{R}}^{n}\times{\mathbb{R}}^{m}\times{\mathbb{R}}^{m\times
d});$

\noindent\ \ \ \ \ \ \ (ii) $\Phi(x)$ is uniformly Lipschitz
with respect to $x\in \mathbb{R}^n$, and for any $x\in
\mathbb{R}^n,\ \Phi(x)\in L^2(\Omega,\mathcal
{F}_T,P;\mathbb{R}^m).$

\vskip0.2cm

 The following monotonicity conditions are also necessary:

\noindent (B2) (i) $\langle A(t,\lambda)-A(t,\overline{\lambda}),\lambda-\overline{\lambda} \rangle \leq -\beta_1|G\widehat{x}|^2-\beta_2(|G^T \widehat{y}|^2+|G^T
 \widehat{z}|^2),$

\noindent\ \ \ \ \ \ \  (ii) $ \langle
\Phi(x)-\Phi(\overline{x}),G(x-\overline{x}) \rangle \geq
\mu_1|G\widehat{x}|^2,\ \ \widehat{x}=x-\overline{x},\
\widehat{y}=y-\overline{y},\ \widehat{z}=z-\overline{z}$,
\\where $\beta_1,\ \beta_2,\ \mu_1$ are nonnegative constants with
$\beta_1 + \beta_2>0,\ \beta_2 + \mu_1>0$. Moreover, we have
$\beta_1>0,\ \mu_1>0 \ (\mbox{resp., }\beta_2>0)$, when $m>n$ (resp.,
$m<n$).

\br \rm{(B2)'-(ii)} $ \langle
\Phi(x)-\Phi(\overline{x}),G(x-\overline{x}) \rangle \geq0.$ \er

 When $\Phi$ does not depend on $x$, i.e., $\Phi(x)=\xi\in
L^2(\Omega,\mathcal {F}_T,P;\mathbb{R}^m)$, the
monotonicity condition (B2)-(i) can be weakened as follows

\noindent(B3) (i) $\langle
A(t,\lambda)-A(t,\overline{\lambda}),\lambda-\overline{\lambda}
\rangle \leq -\beta_1|G\widehat{x}|^2-\beta_2|G^T \widehat{y}|^2,
\ \ \widehat{x}=x-\overline{x},\ \widehat{y}=y-\overline{y},$\\
where $\beta_1,\ \beta_2$ are nonnegative constants with $\beta_1 +
\beta_2>0$. Moreover, we have $\beta_1>0$ (resp., $\beta_2>0$), when
$m>n$ (resp., $m<n$).

\bl Under the assumptions (B1) and (B2), for any initial state
$\zeta \in L^2(\Omega,\mathcal {F}_t,P;\mathbb{R}^n)$, FBSDE (2.2)
associated with $(b, \sigma, f, \zeta, \Phi)$ has a unique adapted
solution $(X_s,Y_s,Z_s)_{s \in [t,T]}\in {\cal{S}}^2(t, T;
{\mathbb{R}^n})\times {\cal{S}}^2(t, T;
{\mathbb{R}^m})\times{\cal{H}}^2(t, T; {\mathbb{R}^{m\times d}}).$
\el

When $\Phi$ does not depend on $x$, i.e.,
$\Phi(x)=\xi\in
L^2(\Omega,\mathcal {F}_T,P;\mathbb{R}^m)$, there is a corresponding result for FBSDE (2.2).
 \bl Under
the assumptions (B1) and (B3), for any initial state $\zeta \in
L^2(\Omega,\mathcal {F}_t,P;\mathbb{R}^n)$\ and the terminal
condition $\Phi(x)=\xi\in L^2(\Omega,\mathcal
{F}_T,P;\mathbb{R}^m)$, FBSDE (2.2) associated with $(b, \sigma, f,
\zeta, \xi)$  has a unique adapted solution $(X_s,Y_s,Z_s)_{s \in
[t,T]}\in {\cal{S}}^2(t, T; {\mathbb{R}^n})\times {\cal{S}}^2(t, T;
{\mathbb{R}^m})\times{\cal{H}}^2(t, T; {\mathbb{R}^{m\times d}}).$
\el

The reader can find the proofs of the Lemmas 2.1 and 2.2 in Peng and
Wu~\cite{PW}.

Now we give the comparison theorem for FBSDEs which will be used in
the later section.
 \bl
(Comparison Theorem) Let $m=1$ and assume that
$(b,\sigma,f,a,\Phi^i)$, for $i=1,2,$ satisfy (B1) and (B2), where
$a\in \mathbb{R}^n$ is the initial state for SDE. Let $(X_s^i,Y_s^i,Z_s^i)_{t\leq
 s\leq T}$ be the solution of FBSDE (2.2) associated with
$(b,\sigma,f,a,\Phi^i)$, respectively. If $\Phi^1(x)\geq \Phi^2(x),\ \mbox{P-a.s.}$
for all $x\in \mathbb{R}^n$. Then, $Y_t^1\geq Y_t^2,\ \mbox{P-a.s.}$
 \el

As a special case, when $\Phi^1(x)\equiv\xi^1,\
\Phi^2(x)\equiv\xi^2$ and $\xi^1 \geq \xi^2,$ we have
 \bl Let $m=1$ and assume
that $(b, \sigma, f, a, \xi^i)$, for $i=1,2,$ satisfy (B1) and (B3),
where $a\in \mathbb{R}^n$ is the initial state for SDE, and $\xi^1,
\xi^2\in L^2(\Omega,\mathcal {F}_T,P;\mathbb{R})$ are the terminal
conditions for related BSDEs, respectively. Let
$(X_s^i,Y_s^i,Z_s^i)_{t\leq s \leq T}$ be the solution of FBSDE
(2.2) associated with $(b,\sigma,f,a,\xi^i).$  If $\xi^1 \geq
\xi^2,\ \mbox{P-a.s.}$ Then, $Y_t^1\geq Y_t^2,\ \mbox{P-a.s.}$
 \el
The Lemmas 2.3 and 2.4 can be found in Wu~\cite{W}.

\section{\large{A DPP for stochastic
optimal control problems of FBSDEs}}
 \hskip1cm
 In this section, we prove the DPP for fully coupled FBSDEs. First we introduce the
background of stochastic optimal control problems. We suppose that
the control state space $U$ is a compact metric space.
${\mathcal{U}}$ is the set of all $U$-valued $ {{\mathbb {F}}
}$-progressively measurable processes. If $u\in \mathcal {U}$, we
call $u$ an admissible control.

For a given admissible control $u(\cdot)\in {\mathcal{U}}$, we
regard $t$ as the initial time and $\zeta \in L^2 (\Omega
,{\mathcal{F}}_t, P;{\mathbb{R}}^n)$ as the initial state. We
consider the following fully coupled forward-backward stochastic
control system
\be
  \left \{
  \begin{array}{ll}
  dX_s^{t,\zeta,u}  =  b(s,X_s^{t,\zeta;u},Y_s^{t,\zeta;u},Z_s^{t,\zeta;u},u_s)ds +
          \sigma(s,X_s^{t,\zeta;u},Y_s^{t,\zeta;u},Z_s^{t,\zeta;u},u_s) dB_s, \\
  dY_s^{t,\zeta;u}  =  -f(s,X_s^{t,\zeta;u},Y_s^{t,\zeta;u},Z_s^{t,\zeta;u},u_s)ds + Z_s^{t,\zeta;u}dB_s, \ \ \ \ \ \ \ \ \ \  s\in[t,T],\\
   X_t^{t,\zeta;u} =  \zeta,\  Y_T^{t,\zeta;u}  =  \Phi(X_T^{t,\zeta;u}),
   \end{array}
   \right.
  \ee
  where the deterministic mappings

  $b: [0,T] \times \mathbb{R}^n \times  \mathbb{R} \times  \mathbb{R}^d  \times U
  \rightarrow \mathbb{R}^n, \ \ \ \ \ \ \ \ \ \ \ \  \  \sigma:  [0,T] \times \mathbb{R}^n \times  \mathbb{R} \times  \mathbb{R}^d \times U
  \rightarrow \mathbb{R}^{n\times d}, $

  $f:  [0,T] \times \mathbb{R}^n \times  \mathbb{R} \times  \mathbb{R}^d \times U
  \rightarrow \mathbb{R},\ \  \ \ \ \ \ \ \ \ \ \ \ \ \Phi: \mathbb{R}^n
  \rightarrow \mathbb{R}$\\
  are continuous to $(t, u)$, and satisfy the assumptions (B1) and (B2), for each $u\in U$, and
  also\\
(B4) there exists a constant $K\geq 0$ such that, for all $t\in [0,
T],\ u\in U,\ x_1,x_2\in\mathbb{R}^n, \ y_1, y_2\in\mathbb{R},\ z_1,
z_2\in\mathbb{R}^d,$
$$|l(t,x_1,y_1,z_1,u)-l(t,x_2,y_2,z_2,u)|\leq K(|x_1-x_2|+|y_1-y_2|+|z_1-z_2|),$$
\ \ \ \ \ \ $l=b,\ \sigma, \ f$, respectively, and $|\Phi(x_1)-\Phi(x_2)|\leq K|x_1-x_2|$.
\br Under our assumptions, it is obvious that there exists a constant $C\geq 0$ such that,
$$|b(t,x,y,z,u)|+|\sigma(t,x,y,z,u)|+|f(t,x,y,z,u)|+|\Phi(x)|\leq C(1+|x|+|y|+|z|),$$
for all $(t,x,y,z,u)\in [0,T]\times
\mathbb{R}^n\times\mathbb{R}\times\mathbb{R}^d\times U.$ Also notice that now (B4) implies (B1).\er

  Hence, for any $u(\cdot) \in \mathcal
  {U},$ from Lemma 2.1, FBSDE (3.1) has a unique solution.

   From Proposition \ref{App-Pro1} in Appendix, there exists $C \in \mathbb{R}^+$ such that, for any $t \in
  [0,T]$, $\zeta, \zeta' \in L^2(\Omega,\mathcal
  {F}_t,P;\mathbb{R}^n),$ $u(\cdot) \in \mathcal
  {U},$ we have, \mbox{P-a.s.}: \be
  \begin{array}{llll}
 &\mbox{(i)}E[\sup\limits_{t\leq s\leq T}|X_s^{t,\zeta;u}-X_s^{t,\zeta';u}|^2+\sup\limits_{t\leq s\leq T}|Y_s^{t,\zeta;u}-Y_s^{t,\zeta';u}|^2+\int_t^T|Z_s^{t,\zeta;u}-Z_s^{t,\zeta';u}|^2ds\mid\mathcal {F}_t]  \leq  C|\zeta - \zeta'|^2, \\
&\mbox{(ii)} E[\sup\limits_{t\leq s\leq T}|X_s^{t,\zeta;u}|^2+\sup\limits_{t\leq s\leq T}|Y_s^{t,\zeta;u}|^2+ \int_t^T|Z_s^{t,\zeta;u}|^2ds \mid\mathcal {F}_t] \leq   C(1 +|\zeta|^2). \\
  \end{array}
  \ee
Therefore, we get   \be\begin{array}{llll}&&{\rm(i)}\ \ \
  |Y_t^{t,\zeta;u}| \leq C(1+|\zeta|), \ \mbox{P-a.s.};\\
  &&{\rm(ii)}\ \ |Y_t^{t,\zeta;u}-Y_t^{t,\zeta';u}| \leq C|\zeta - \zeta'|,\ \mbox{P-a.s.}  \end{array}\ee

We now introduce the subspaces of admissible controls. An admissible
control process $u=(u_r)_{r\in [t,s]}$ on $[t,s]$ is an
$\mathbb{F}$-progressively measurable, $U$-valued process. The set
of all admissible controls on $[t,s]$ is denote $\mathcal
{U}_{t,s},\ t\leq s\leq T.$

For a given process $u(\cdot) \in \mathcal {U}_{t,T}$, we define the
associated cost functional as follows:
  \be J(t,x;u):=Y_s^{t,x;u}\mid_{s=t}, \ \ (t,x)\in[0,T] \times
  \mathbb{R}^n,\ee where the process $Y^{t,x;u}$ is defined by FBSDE
  (3.1).

From Theorem \ref{App-Th1} we have, for any $t \in
  [0,T]$ and $\zeta \in L^2(\Omega,\mathcal {F}_t,P;\mathbb{R}^n),$
 \be J(t,\zeta;u)=Y_t^{t,\zeta;u},  \mbox{ P-a.s.} \ee
For $ \zeta = x \in \mathbb{R}^n,$ we define the value
  function as \be W(t,x) := \esssup_{u\in\mathcal {U}_{t,T}}J(t,x;u). \ee
  \br Thanks to the assumptions (B1) and (B2), the value
  function $W(t,x)$ is well defined and it is a bounded ${\mathcal {F}}_t$-measurable random variable. But it turns out to be deterministic. \er
  Inspired by the method in Buckdahn and Li~\cite{BL}, we can prove
  that $W$ is deterministic.
  \bp We assume the assumptions (B1) and (B2) hold. Then, for any $(t,x) \in [0,T] \times \mathbb{R}^n,$ $W(t,x)$ is a deterministic function in the sense that $W(t,x) =
  E[W(t,x)],\
  \mbox{P-a.s.}$
  \ep
  \noindent \textbf{Proof}. Let $H$ denote the Cameron-Martin space of all absolutely
  continuous elements $h \in \Omega$ whose derivative $\dot{h}$
  belongs to $L^2([0,T];\mathbb{R}^d).$

  For any $h \in H$, we define the mapping $\tau_h \omega := \omega + h, \ \omega \in
  \Omega.$ It is easy to check that $\tau_h : \Omega \rightarrow \Omega$ is a
  bijection, and its law is given by $P \circ [\tau_h]^{-1} = \exp\{\int_0^T \dot{h}_sdB_s - \frac{1}{2}\int_0^T
  |\dot{h}_s|^2ds\}P.$  For any fixed $(t,x) \in [0,T] \times \mathbb{R}^n$, set $H_t = \{h \in H |h(\cdot) =h(\cdot \wedge
  t)\}$. The proof can be separated into the following three steps:

(1). For all $u\in\mathcal {U}_{t,T}, \ h\in H_t, \
J(t,x;u)(\tau_h)=J(t,x;u(\tau_h)),  \ \mbox{P-a.s.}$

  In fact, using the Girsanov transformation to FBSDE (3.1) (with $\zeta =
  x$) and comparing the obtained equation with the FBSDE obtained from
  (3.1) by replacing the transformed control process $u(\tau_h)$
  for $u$, due to the uniqueness of the solution of (3.1) we obtain
   $$\begin{array}{llll}X_s^{t,x;u}(\tau_h) &=& X_s^{t,x;u(\tau_h)}, \ \mbox{for any} \ s \in
   [t,T],\
  \mbox{P-a.s.},\\
 Y_s^{t,x;u}(\tau_h) &=& Y_s^{t,x;u(\tau_h)}, \ \mbox{for any} \ s  \in
  [t,T], \ \mbox{P-a.s.},\\
  Z_s^{t,x;u}(\tau_h) &=& Z_s^{t,x;u(\tau_h)}, \ \mbox{dsdP-a.e. on} \ [0,T] \times \Omega.
  \end{array}$$
  Hence, $J(t,x;u)(\tau_h) = J(t,x;u(\tau_h)), \ \mbox{P-a.s.}$

  (2). For any $h \in H_t,$ we have $$\{\esssup_{u \in \mathcal
  {U}_{t,T}}J(t,x;u)\}(\tau_h) = \esssup_{u \in \mathcal
  {U}_{t,T}}\{J(t,x;u)(\tau_h)\}, \ \mbox{P-a.s.}$$
  In fact, for convenience, setting $I(t,x) = \esssup_{u \in \mathcal
  {U}_{t,T}}J(t,x;u)$, we have $I(t,x) \geq J(t,x;u)$. Then, $I(t,x)(\tau_h) \geq J(t,x;u)(\tau_h),
 \ \mbox{P-a.s.},$ for all  $u \in \mathcal
  {U}_{t,T}.$
 Therefore,
 $\{\esssup_{u \in \mathcal
  {U}_{t,T}}J(t,x;u)\}(\tau_h) \geq \esssup_{u \in \mathcal
  {U}_{t,T}}\{J(t,x;u)(\tau_h)\}, $ $\ \mbox{P-a.s.}$\ On the other hand, for any random variable $\xi$ which satisfies $\xi \geq
  J(t,x;u)(\tau_h)$, we have $\xi(\tau_{-h}) \geq
  J(t,x;u),\ \mbox{P-a.s.}$, for all $u \in \mathcal
  {U}_{t,T}$. So $\xi(\tau_{-h}) \geq
  I(t,x), \ \mbox{P-a.s.}$, i.e. $\xi \geq
  I(t,x)(\tau_h), \ \mbox{P-a.s.}$ Thus,
   $J(t,x;u)(\tau_h) \geq \{\esssup_{u \in \mathcal
  {U}_{t,T}}J(t,x;u)\}(\tau_h), \ \mbox{P-a.s., for any }u\in\mathcal {U}_{t,T}.$\
  Therefore, $$\esssup_{u \in \mathcal {U}_{t,T}}\{J(t,x;u)(\tau_h)\} \geq \{\esssup_{u \in \mathcal {U}_{t,T}}J(t,x;u)\}(\tau_h),\  \mbox{P-a.s.}$$
  From above we get $\{\esssup_{u \in \mathcal {U}_{t,T}}J(t,x;u)\}(\tau_h) = \esssup_{u \in \mathcal {U}_{t,T}}\{J(t,x;u)(\tau_h)\},\  \mbox{P-a.s.}$

  (3). Under the Girsanov
  transformation $\tau_h$, $W(t,x)$ is invariant, i.e., $$W(t,x)(\tau_h) = W(t,x),
 \ \mbox{P-a.s.}, \
  \mbox{for any} \ h \in H.$$
  From the first step and the second one, for all $h \in
  H_t$, we have
 $$ \begin{array}{llll}
  W(t,x)(\tau_h) & = & \esssup_{u \in \mathcal {U}_{t,T}}J(t,x;u)(\tau_h)
                  =  \esssup_{u \in \mathcal {U}_{t,T}}\{J(t,x;u)(\tau_h)\} \\
                 & = & \esssup_{u \in \mathcal {U}_{t,T}}J(t,x;u(\tau_h))
                  =  W(t,x), \ \mbox{P-a.s.}
   \end{array}$$
  In the latter equality we have used $\{u(\tau_h)\mid u(\cdot) \in \mathcal {U}_{t,T}\} = \mathcal
  {U}_{t,T}$. Therefore, for any $h \in H_t,\  W(t,x)(\tau_h) =
  W(t,x),\
  \mbox{P-a.s.},$ and since $W(t,x)$ is $\mathcal {F}_t$-measurable, we
  have this relation for all $ h \in H$.

 Combined with the following auxiliary
  lemma we can complete the proof. \endpf

  \bl  Let $\zeta$ be a random variable defined over our
classical Wiener space $(\Omega,\mathcal {F}_T, P)$, such that
$\zeta(\tau_h) = \zeta, \ \mbox{P-a.s.}$, for any $h \in H$. Then
$\zeta = E\zeta,\ \mbox{P-a.s.}$ \el
\noindent Its proof can be found in Buckdahn and
Li~\cite{BL}.\\

 From (3.3) and (3.6)-the definition of the value function
$W(t,x)$, we get the following property:

\bl There exists a constant $C>0$ such that, for all $0 \leq t \leq
T,\ x,x' \in \mathbb{R}^n,$ \be \begin{array}{llll}&& {\rm(i)} \ \
|W(t,x)-W(t,x')| \leq C|x-x'|; \\
&&  {\rm(ii)}\ \ |W(t,x)| \leq C(1+|x|). \end{array}\ee \el

\bl Under the assumptions (B1) and (B2), the cost functional
$J(t,x;u),$\ for any $u\in \mathcal {U}_{t,T}$, and the value
function $W(t,x)$ are monotonic in the following sense: for each $x,
\bar{x} \in \mathbb{R}^n,$ $t\in [0,T],$
$$\begin{array}{llll}
&&{\rm(i)}\ \  \langle
J(t,x;u)-J(t,\bar{x};u),\ G(x-\bar{x})\rangle\geq 0,\ \mbox{P-a.s.};\\
&&{\rm(ii)}\ \ \langle
W(t,x)-W(t,\bar{x}),\ G(x-\bar{x})\rangle \geq 0.\end{array}$$  \el
 \noindent \textbf{Proof}. We define
$\hat{X}_s=X_s^{t,x;u}-X_s^{t,\bar{x};u}, \
\hat{Y}_s=Y_s^{t,x;u}-Y_s^{t,\bar{x};u},\
\hat{Z}_s=Z_s^{t,x;u}-Z_s^{t,\bar{x};u}, $\ and $\Delta
h(s)=h(s,X_s^{t,x;u},Y_s^{t,x;u},Z_s^{t,x;u},u_s)-h(s,X_s^{t,\bar{x};u},Y_s^{t,\bar{x};u},Z_s^{t,\bar{x};u},u_s),$\ for $h=b,\ \sigma,\ f,\ A,$ respectively.\\
 Applying It\^o's
formula to $\langle \hat{Y}_s,G\hat{X}_s\rangle$, we
get immediately from (B2)
 $$\begin{array}{llll} &&\langle
J(t,x;u)-J(t,\bar{x};u),G(x-\bar{x})\rangle =E[\langle
Y_t^{t,x;u}-Y_t^{t,\bar{x};u},G(x-\bar{x})\rangle\mid\mathcal
{F}_t]\geq 0,\ \ \  \mbox{for any }u\in\mathcal
{U}_{t,T}.\\
\end{array}$$
From the definition of  $W(t,x)$, we always have
$W(t,x)\geq J(t,x;u),\ \mbox{P-a.s.}$ for any $u\in\mathcal {U}_{t,T}$.\ On the other hand, similar to Remark 3.5-(ii), we can get, for any $\varepsilon
>0,$ the existence of $u^\varepsilon\in\mathcal {U}_{t,T}$, such that
$W(t,\bar{x}) \leq
J(t,\bar{x};u^\varepsilon)+\varepsilon.$\\
If $G(x-\bar{x})\geq 0,$ then $\begin{array}{llll} \langle W(t,x)-W(t,\bar{x}),G(x-\bar{x})\rangle\geq
(J(t,x;u^\varepsilon)-J(t,\bar{x};u^\varepsilon)-\varepsilon)G(x-\bar{x})\geq
-\varepsilon G(x-\bar{x}).
\end{array}$\\
If $G(x-\bar{x})\leq 0,$ then for $u^\varepsilon$ such that
$W(t,x)\leq J(t,x;u^\varepsilon)+\varepsilon,$\ $\langle W(t,x)-W(t,\bar{x}),G(x-\bar{x})\rangle
\geq -\varepsilon G(\bar{x}-{x}).$\\
Therefore, $\langle
W(t,x)-W(t,\bar{x}),G(x-\bar{x})\rangle \geq
 -\varepsilon |G(x-\overline{x})|,$ for any
$x, \bar{x}\in\mathbb{R}^n, \ t\in[0,T].$\
Consequently, letting $\varepsilon\downarrow 0$, $\langle
W(t,x)-W(t,\bar{x}),G(x-\bar{x})\rangle \geq 0,$ for any $x,
\bar{x}\in\mathbb{R}^n,$ $t\in[0,T].$\endpf

\br (1) From (B2)-(i) we see that if $\sigma$\ doesn't depend on $z$, then $\beta_2=0$. Furthermore, we assume that:\\ $\begin{array}{llll}
&&(B5)\  \mbox{the Lipschitz constant}\ L_\sigma\geq 0 \ \mbox{of}\ \sigma\   \mbox{with respect to}\ z\ \mbox{is sufficiently small, i.e., there exists some}\\
&&\ \ \ \ \ \ L_\sigma\geq 0\ \mbox{small enough such that}, \ \mbox{for all}\ t\in[0, T],\ u\in U,\ x_1, x_2\in\mathbb{R}^n,\ y_1, y_2\in\mathbb{R},\ z_1, z_2\in\mathbb{R}^d,\\
&&\ \ \ \ \ \ |\sigma(t,x_1,y_1,z_1,u)-\sigma(t,x_2,y_2,z_2,u)|\leq K(|x_1- x_2|+|y_1-y_2|)+L_\sigma|z_1-z_2|.
\end{array}$\\
(2) On the other hand, notice that when $\sigma$\ doesn't depend on
$z$ it's obvious that (B5) always holds true. \er

 The notation of stochastic backward
  semigroup was first introduced by Peng~\cite{Pe4} and was applied to
  prove the DPP for stochastic control problems. Now we discuss a generalized DPP for our stochastic optimal control problem (3.1), (3.6).
   For this we have to adopt Peng's notion of stochastic backward
  semigroup, and to define the family of (backward) semigroups associated with FBSDE (3.1).

    For given initial data $(t,x)$, a real number $\delta\in (0, T -t]$, an admissible control process $u(\cdot)\ \in \ \mathcal
    {U}_{t,t+\delta}$ and a real-valued random function $\Psi: \Omega\times \mathbb{R}^n\rightarrow \mathbb{R}$,
    ${\cal F}_{t+\delta}\otimes {\cal B}(\mathbb{R}^n)$-measurable such that (B2)-(ii) holds, we put
    $$G_{s,t+\delta}^{t,x;u}[\Psi(t+\delta, \widetilde{X}_{t+\delta}^{t,x;u})]:=\widetilde{Y}_s^{t,x;u}, \  s \in
    [t,t+\delta],$$
    where $(\widetilde{X}_s^{t,x;u},\widetilde{Y}_s^{t,x;u},\widetilde{Z}_s^{t,x;u})_{t \leq s \leq
    t+\delta}$ is the solution of the following FBSDE with the time
    horizon $t+\delta$:
\be
  \left \{
  \begin{array}{ll}
d\widetilde{X}_s^{t,x;u}  =b(s,\widetilde{X}_s^{t,x;u},\widetilde{Y}_s^{t,x;u},\widetilde{Z}_s^{t,x;u},u_s)ds + \sigma(s,\widetilde{X}_s^{t,x;u},\widetilde{Y}_s^{t,x;u},\widetilde{Z}_s^{t,x;u},u_s)dB_s, \\
 d\widetilde{Y}_s^{t,x;u}  =  -f(s,\widetilde{X}_s^{t,x;u},\widetilde{Y}_s^{t,x;u},\widetilde{Z}_s^{t,x;u},u_s)ds + \widetilde{Z}_s^{t,x;u}dB_s, \ \ \ \ \ s\in [t,t+\delta],\\
  \widetilde{X}_t^{t,x;u}  =  x,\ \widetilde{Y}_{t+\delta}^{t,x;u}  =  \Psi(t+\delta, \widetilde{X}_{t+\delta}^{t,x;u}).
   \end{array}
   \right.
  \ee
\br (1) From Lemma 2.1 and Lemma 2.2 if $\Psi$\ doesn't depend on $x$, we know FBSDE (3.8) has a unique solution
$(\widetilde{X}^{t,x;u},\widetilde{Y}^{t,x;u},\widetilde{Z}^{t,x;u}).$\\
 (2) We also point out that if $\Psi$\ is Lipschitz with respect to $x$, FBSDE (3.8) can be also solved under the assumptions (B4) and (B5)
  on the small interval $[t, t+\delta]$, for any $0\leq \delta\leq \delta_0,$ where small enough $\delta_0>0$
   is independent of $(t, x)$ and the control $u$, from Proposition \ref{App-Pro4}. \er

Since $\Phi$\ satisfies (B2)-(ii) the solution
$(X^{t,x;u},Y^{t,x;u},Z^{t,x;u})$ of FBSDE (3.1) exists and we get
$$ G_{t,T}^{t,x;u}[\Phi(X_T^{t,x;u})] =
G_{t,t+\delta}^{t,x;u}[Y_{t+\delta}^{t,x;u}]. $$ Moreover, we have
\be J(t,x;u)=Y_t^{t,x;u}=G_{t,T}^{t,x;u}[\Phi(X_T^{t,x;u})]
=G_{t,t+\delta}^{t,x;u}[Y_{t+\delta}^{t,x;u}]=G_{t,{t+\delta}}^{t,x;u}[J(t+\delta,X_{t+\delta}^{t,x;u};u)].\ee
 \bt Under the assumptions (B2), (B4) and (B5), the value function $W(t,x)$
 satisfies the following DPP: there exists sufficiently small $\delta_0>0$, such that for any $0\leq \delta \leq \delta_0,\ t\in [0, T-\delta],\ x \in \mathbb{R}^n,$
$$W(t,x)=\esssup_{u\in \mathcal
{U}_{t,t+\delta}}G_{t,{t+\delta}}^{t,x;u}[W(t+\delta,\widetilde{X}_{t+\delta}^{t,x;u})].$$
\et
\noindent \textbf{Proof}. With the help of Lemma 3.2, (3.9), Theorem 5.2, Corollary 5.1, and Proposition 5.4, adapting the method of the proof of Theorem 3.6 in~\cite{BL}, we can complete the proof.\endpf
From its proof we can get
\begin{remark}\label{r3.5}  {\rm(i)} For all $u\in {\cal{U}}_{t, t+\delta},$
$$\label{56} W(t,x)(=W_\delta(t, x))\geq G^{t,x; u}_{t,t+\delta}
      [W(t+\delta, \widetilde{X}^{t,x;u}_{t+\delta})],\quad \mbox{{\it P}-a.s.}
$$
 {\rm(ii)} For any $(t, x)\in
[0,T]\times {\mathbb{R}}^n,$\ $\delta \in [0, \delta_0]$\ and
$\varepsilon>0$, there exists some $u^{\varepsilon}(\cdot) \in {\cal{U}}_{t,
t+\delta}$\ such that
 $$\label{55} W(t,x)(=W_\delta(t, x))\leq G^{t,x;
u^{\varepsilon}}_{t,t+\delta}
      [W(t+\delta, \widetilde{X}^{t,x; u^{\varepsilon}}_{t+\delta})]+\varepsilon,\ \mbox{{\it P}-a.s.}
$$
\end{remark}
 Notice that from the definition of our stochastic backward
  semigroup we know here $$G_{s,t+\delta}^{t,x;u}[W(t+\delta,\widetilde{X}_{t+\delta}^{t,x;u})]=\widetilde{Y}_s^{t,x;u}, \  s \in
    [t,t+\delta],\ u(\cdot) \in  \mathcal
    {U}_{t,t+\delta},$$
    where $(\widetilde{X}_s^{t,x;u},\widetilde{Y}_s^{t,x;u},\widetilde{Z}_s^{t,x;u})_{t \leq s \leq
    t+\delta}$ is the solution of the following FBSDE with the time
    horizon $t+\delta$:
\be
  \left \{
  \begin{array}{ll}
d\widetilde{X}_s^{t,x;u}  =  b(s,\widetilde{X}_s^{t,x;u},\widetilde{Y}_s^{t,x;u},\widetilde{Z}_s^{t,x;u},u_s)ds + \sigma(s,\widetilde{X}_s^{t,x;u},\widetilde{Y}_s^{t,x;u},\widetilde{Z}_s^{t,x;u},u_s)dB_s, \\
 d\widetilde{Y}_s^{t,x;u}  =  -f(s,\widetilde{X}_s^{t,x;u},\widetilde{Y}_s^{t,x;u},\widetilde{Z}_s^{t,x;u},u_s)ds + \widetilde{Z}_s^{t,x;u}dB_s, \ \ \ \ \ s\in [t,t+\delta],\\
  \widetilde{X}_t^{t,x;u}  =  x,\
  \widetilde{Y}_{t+\delta}^{t,x;u}  =  W(t+\delta,\widetilde{X}_{t+\delta}^{t,x;u}).
   \end{array}
   \right.
  \ee
Due to Proposition \ref{App-Pro4} there exists sufficiently small
$\delta_0>0$, such that for any $0\leq \delta \leq \delta_0,$\ the
above equation (3.10) has a unique solution
$(\widetilde{X}^{t,x;u},\widetilde{Y}^{t,x;u},\widetilde{Z}^{t,x;u})$\
on the time interval $[t, t+\delta]$.

From Lemma 3.2, we get the value function $W(t,x)$ is Lipschitz
continuous in $x$, uniformly in $t$. Now we can get the continuity
property of $W(t,x)$ in $t$ with the help of Theorem 3.1.

 \bt  Under (B2), (B4) and (B5), the
 value function $W(t,x)$ is continuous in $t$.\et
 \noindent \textbf{Proof}. Let $(t,x)\in[0,T]\times\mathbb{R}^n$. In order to obtain $W$ is continuous in $t$, it is sufficient to prove the following inequality: there exists some constant $C$, such that
 $$-C(1+|x|)\delta^{\frac{1}{2}} \leq W(t,x)-W(t+\delta,x) \leq
 C(1+|x|)\delta^{\frac{1}{2}},\ \mbox{for all}\ 0\leq\delta \leq T-t\ \mbox{sufficiently small}.$$
We will only prove the second inequality, the proof of the first one is similar.

From Remark \ref{r3.5}, there exists $u^\varepsilon\in\mathcal
{U}$, such that
 $$G_{t,{t+\delta}}^{t,x;u^\varepsilon}[W(t+\delta,\widetilde{X}_{t+\delta}^{t,x;u^\varepsilon})]+\varepsilon
 \geq W(t,x) \geq
 G_{t,{t+\delta}}^{t,x;u^\varepsilon}[W(t+\delta,\widetilde{X}_{t+\delta}^{t,x;u^\varepsilon})].$$
Therefore, $W(t,x)-W(t+\delta,x) \leq
 G_{t,{t+\delta}}^{t,x;u^\varepsilon}[W(t+\delta,\widetilde{X}_{t+\delta}^{t,x;u^\varepsilon})]+\varepsilon-W(t+\delta,x)
 =I_\delta^1+I_\delta^2+\varepsilon,$\\
 where$$\begin{array}{llll}I_\delta^1& =&
 G_{t,{t+\delta}}^{t,x;u^\varepsilon}[W(t+\delta,\widetilde{X}_{t+\delta}^{t,x;u^\varepsilon})]-
 G_{t,{t+\delta}}^{t,x;u^\varepsilon}[W(t+\delta,x)],\\
 I_\delta^2 &=&
 G_{t,{t+\delta}}^{t,x;u^\varepsilon}[W(t+\delta,x)]-
 W(t+\delta,x).\end{array}$$
 Notice also that $G_{s,{t+\delta}}^{t,x;u}[W(t+\delta,x)]=\hat{Y}_s^{t,x;u}, \  s \in
    [t,t+\delta],\ u(\cdot)\in  \mathcal
    {U}_{t,t+\delta},$
    where $(\hat{X}_s^{t,x;u},\hat{Y}_s^{t,x;u},\hat{Z}_s^{t,x;u})_{t \leq s \leq
    t+\delta}$ is the solution of the following FBSDE with the time
    horizon $t+\delta$:
\be
  \left \{
  \begin{array}{ll}
d\hat{X}_s^{t,x;u}  =  b(s,\hat{X}_s^{t,x;u},\hat{Y}_s^{t,x;u},\hat{Z}_s^{t,x;u},u_s)ds + \sigma(s,\hat{X}_s^{t,x;u},\hat{Y}_s^{t,x;u},\hat{Z}_s^{t,x;u},u_s)dB_s, \\
 d\hat{Y}_s^{t,x;u}  =  -f(s,\hat{X}_s^{t,x;u},\hat{Y}_s^{t,x;u},\hat{Z}_s^{t,x;u},u_s)ds + \hat{Z}_s^{t,x;u}dB_s, \ \ \ \ \ s\in [t,t+\delta],\\
  \hat{X}_t^{t,x;u}  =  x,\
  \hat{Y}_{t+\delta}^{t,x;u}  =  W(t+\delta,x).
   \end{array}
   \right.
  \ee
 Applying It\^o's formula to $e^{\beta s}|\widetilde{Y}_s^{t,x;u^\varepsilon}-\hat{Y}_s^{t,x;u^\varepsilon}|^2$,
 by taking $\beta$\ large enough and using standard methods for BSDEs, we get with the help of (3.7) and Proposition \ref{App-Pro5}-(ii) for equations (3.10) and (3.11) that
 \be
  \begin{array}{llll}
& &|\widetilde{Y}_t^{t,x;u^\varepsilon}-\hat{Y}_t^{t,x;u^\varepsilon}|^2\\
 &\leq & CE[|W(t+\delta,\widetilde{X}_{t+\delta}^{t,x;u^\varepsilon})-W(t+\delta,x)|^2\mid\mathcal
{F}_t]+CE[\int_t^{t+\delta}|\widetilde{X}_r^{t,x;u^\varepsilon}-\hat{X}_r^{t,x;u^\varepsilon}|^2dr\mid\mathcal
{F}_t]\\
 &\leq & CE[|\widetilde{X}_{t+\delta}^{t,x;u^\varepsilon}-x|^2\mid\mathcal
{F}_t]+C\delta(E[\sup\limits_{t\leq r\leq
t+\delta}|\widetilde{X}_{r}^{t,x;u^\varepsilon}-x|^2\mid\mathcal
{F}_t]+ E[\sup\limits_{t\leq r\leq
t+\delta}|\hat{X}_{r}^{t,x;u^\varepsilon}-x|^2\mid\mathcal
{F}_t])\\
 &\leq & C\delta (1+|x|^2),\ \mbox{P-a.s.}
   \end{array}
  \ee
That is, there exists some
constant $C$ independent of the controls such that
$$|I_\delta^1| =|\widetilde{Y}_t^{t,x;u^\varepsilon}-\hat{Y}_t^{t,x;u^\varepsilon}| \leq C(1+|x|)\delta^\frac{1}{2},\ \mbox{P-a.s.}$$
From equation (3.11), Remark 3.1, and Proposition \ref{App-Pro5}-(i) (the estimates for FBSDE (3.11))
$$\begin{array}{llll}
|I_\delta^2|&=&|E[W(t+\delta,x)+\int_t^{t+\delta}f(s,\hat{X}_s^{t,x;u^\varepsilon},\hat{Y}_s^{t,x;u^\varepsilon},
\hat{Z}_s^{t,x;u^\varepsilon},u_s^\varepsilon)ds\mid\mathcal
{F}_t]-W(t+\delta,x)|\\
&\leq&
C\delta^\frac{1}{2}E[\int_t^{t+\delta}(1+|\hat{X}_s^{t,x;u^\varepsilon}|+|\hat{Y}_s^{t,x;u^\varepsilon}|
+|\hat{Z}_s^{t,x;u^\varepsilon}|)^2ds\mid\mathcal
{F}_t]^\frac{1}{2}\\
&\leq &C(1+|x|)\delta^\frac{1}{2}.\end{array}$$
Therefore, $W(t,x)-W(t+\delta,x) \leq
 C(1+|x|)\delta^{\frac{1}{2}}+\varepsilon.$ Letting $\varepsilon\downarrow0$, we complete the proof. \endpf

 \section{\large Viscosity solutions of HJB equations}
  \hskip1cm
In this section we show that the value function $W(t, x)$ defined in
(3.6) is a viscosity solution of the corresponding HJB equation. For
this we use Peng's BSDE
 approach~\cite{Pe4} developed from stochastic control problems of decoupled
 FBSDEs, but still more difficulties for two cases, especially for
 Case 2.

 \textbf{Case 1}. We suppose that $\sigma$ does not depend on $z$, but depends on $u$.

Then the equation (3.1) becomes the following equation (4.1):
  \be
  \left \{
  \begin{array}{ll}
  dX_s^{t,x;u}  =  b(s,X_s^{t,x;u},Y_s^{t,x;u},Z_s^{t,x;u},u_s)ds +
          \sigma(s,X_s^{t,x;u},Y_s^{t,x;u},u_s) dB_s, \\
  dY_s^{t,x;u}  =  -f(s,X_s^{t,x;u},Y_s^{t,x;u},Z_s^{t,x;u},u_s)ds + Z_s^{t,x;u}dB_s,\ \ \ s\in [t,T], \\
   X_t^{t,x;u} =  x,\
   Y_T^{t,x;u}  =  \Phi(X_T^{t,x;u}).
   \end{array}
   \right.
  \ee
 We consider the following HJB equation: \be
 \left \{\begin{array}{ll}
 &\!\!\!\!\! \frac{\partial }{\partial t} W(t,x) +  H(t, x, W(t,x), DW(t,x), D^2W(t,x))=0,
 \hskip 0.5cm   (t,x)\in [0,T)\times {\mathbb{R}},  \\
 &\!\!\!\!\!  W(T,x) =\Phi (x), \hskip0.5cm   x \in {\mathbb{R}}.
 \end{array}\right.
\ee
 In this case the Hamiltonian is given by
$$\begin{array}{lll}
& H(t, x, y, p, X)= \sup\limits_{u \in U}\{ p.b(t, x, y,  p.\sigma,
u)+\frac{1}{2}\tr(\sigma\sigma^{T}(t,
x, y,u)X)+f(t, x, y, p.\sigma, u)\},
 \end{array}$$
where $t\in [0, T],\ x\in{\mathbb{R}},\ y\in \mathbb{R},\ p\in
{\mathbb{R}}$, and $X\in{\mathbb{R}}$.

  Let us first recall the definition of a viscosity
solution of equation (4.2). More details on viscosity solutions can
be found in Crandall, Ishill and Lions~\cite{CIL}.

 \bde\mbox{ }A real-valued
continuous function $W\in C([0,T]\times {\mathbb{R}}^k )$ is called \\
   \noindent{\rm(i)} a viscosity subsolution of equation (4.2) if $W(T,x) \leq \Phi (x),\ \mbox{for all}\ x \in
  {\mathbb{R}}^k$, and if for all functions $\varphi \in C^3_{l, b}([0,T]\times
  {\mathbb{R}}^k)$ and for all $(t,x) \in [0,T) \times {\mathbb{R}}^k$ such that $W-\varphi $\ attains
  a local maximum at $(t, x)$,
   $$
     \frac{\partial \varphi}{\partial t} (t,x) + H(t,x,\varphi,D\varphi,D^2\varphi) \geq
     0;$$
\noindent{\rm(ii)} a viscosity supersolution of equation (4.2) if
$W(T,x) \geq \Phi (x), \mbox{for all}\ x \in
  {\mathbb{R}}^k$, and if for all functions $\varphi \in C^3_{l, b}([0,T]\times
  {\mathbb{R}}^k)$ and for all $(t,x) \in [0,T) \times {\mathbb{R}}^k$\ such that $W-\varphi $\ attains
  a local minimum at $(t, x)$,
$$\frac{\partial \varphi}{\partial t} (t,x) +
H(t,x,\varphi,D\varphi,D^2\varphi) \leq 0;$$
 \noindent{\rm(iii)} a viscosity solution of equation (4.2) if it is both a viscosity sub- and supersolution of equation
     (4.2). \ede

\br \mbox{  }$C^3_{l, b}([0,T]\times {\mathbb{R}}^k)$ denotes the
set of the real-valued functions that are continuously
differentiable up to the third order and whose derivatives of order
from 1 to 3 are bounded.\er

\bt  Under the assumptions (B2) and (B4), the value function
$W(t,x)$ defined in (3.6) is a viscosity solution of (4.2).\et
\noindent\textbf{Proof}.  Obviously, $W(T,x)=\Phi(x),\
x\in\mathbb{R}$. Let us show
  that $W$ is a viscosity subsolution, the proof for the viscosity supersolution is similar. We suppose that $\varphi \in
  C_{l,b}^3([0,T]\times\mathbb{R})$ and
 that $(t,x)\in[0,T)\times\mathbb{R}$ is such that $W-\varphi$
  attains its maximum at $(t,x)$. Since $W$ is continuous
  and of at most linear growth, we only need to consider
  the global maximum at $(t,x)$. Without loss of generality we may assume that $\varphi(t,x)=W(t,x).$ We consider the
  following equation:
 \be
  \left \{\begin{array}{ll}
  dX^u_s  =  b(s,X^u_s,Y^u_s,Z^u_s,u_s)ds + \sigma(s,X^u_s,Y^u_s,u_s)dB_s, \\
 dY^u_s  =  -f(s,X^u_s,Y^u_s,Z^u_s,u_s)ds + Z^u_s dB_s, \ \  \ \ \ \  s\in [t,t+\delta],\\
 X_t^u = x,\
  Y^u_{t+\delta} =  \varphi(t+\delta,X_{t+\delta}^{u}),\ \ 0\leq
  \delta \leq T-t.
   \end{array}
   \right.
  \ee
From Propositions \ref{App-Pro4} and \ref{App-Pro5} in Appendix, we
know there exists sufficiently small $0< {\delta}_1< T-t$, such that
for any $0\leq \delta\leq {\delta}_1$, FBSDE (4.3) has a unique
solution $(X_s^u,Y_s^u,Z_s^u)_{t\leq s\leq t+\delta}\in \mathcal
{S}^2(t,t+\delta;\mathbb{R})\times
\mathcal{S}^2(t,t+\delta;\mathbb{R})\times
\mathcal{H}^2(t,t+\delta;\mathbb{R}^d),$ and for $p\geq 2,$
\be\begin{array}{llll} {\rm(i)} &&E[\sup\limits_{t\leq s\leq
t+\delta}|X_s^u|^p+\sup\limits_{t\leq s\leq
t+\delta}|Y_s^u|^p+(\int_t^{t+\delta}|Z_s^u|^2ds)^{\frac{p}{2}}\mid
\mathcal
{F}_t]\leq C(1+|x|^p),\ \mbox{P-a.s.};\\
{\rm(ii)}&&E[\sup\limits_{t\leq s\leq t+\delta}|X_s^u-x|^p\mid
\mathcal {F}_t]\leq C\delta^{\frac{p}{2}}(1+|x|^p),\ \mbox{P-a.s.};\\
{\rm(iii)}&&E[(\int_t^{t+\delta}|Z_s^u|^2ds)^{\frac{p}{2}}\mid
\mathcal {F}_t]\leq C\delta^{\frac{p}{2}}(1+|x|^p),\ \mbox{P-a.s.}
\end{array}
  \ee
From the definition of the backward stochastic semigroup for fully
coupled FBSDE, we have \be
G_{s,t+\delta}^{t,x;u}[\varphi(t+\delta,X_{t+\delta}^{u})]=Y_s^u,\ \
s\in [t,t+\delta],\ 0\leq \delta \leq {\delta}_1.\ee From the DPP
(Theorem 3.1), we have
  $$\varphi(t,x)=W(t,x)=\esssup_{u\in \mathcal
{U}_{t,t+\delta}}G_{t,{t+\delta}}^{t,x;u}[W(t+\delta,\widetilde{X}_{t+\delta}^{t,x,u})],\
0\leq \delta \leq  {\delta}_1,$$
where $\widetilde{X}^{t,x,u}$\ is defined by FBSDE (3.10).\\
 From $W(s,y)\leq\varphi(s,y),\ (s,y)\in[0,T)\times\mathbb{R},$ and the monotonicity
property of $G_{t,{t+\delta}}^{t,x;u}[\cdot]$ (see Theorem
\ref{App-Th2} in Appendix) we have,  \be \esssup_{u\in \mathcal
{U}_{t,t+\delta}}G_{t,{t+\delta}}^{t,x;u}[\varphi(t+\delta,X_{t+\delta}^{u})]-\varphi(t,x)\geq0,\
\ 0\leq \delta \leq \delta_1.\ee Now we define \be
\begin{array}{llll}Y_s^{1,u}&=&Y_s^u-\varphi(s,X_s^u)\\
&=&\int_s^{t+\delta}f(r,X_r^u,Y_r^u,Z_r^u,u_r)dr-\int_s^{t+\delta}Z_r^udB_r+\varphi(t+\delta,X_{t+\delta}^u)-\varphi(s,X_s^u).
\end{array}\ee
Using It\^o's formula to $\varphi(s,X_s^u)$, and setting
$Z_s^{1,u}=Z_s^u-D\varphi(s,X_s^u).\sigma(s,X_s^u,Y_s^u,u_s)$, we
have \be\left
\{\begin{array}{llll}Y_s^{1,u}&=&\int_s^{t+\delta}[\frac{\partial}{\partial
r}\varphi(r,X_r^u)+D\varphi(r,X_r^u).b(r,X_r^u,Y_r^u,Z_r^u,u_r)
+\frac{1}{2}\tr(\sigma\sigma^T(r,X_r^u,Y_r^u,u_r)D^2\varphi(r,X_r^u))\\
&&+f(r,X_r^u,Y_r^u,Z_r^u,u_r)]dr-\int_s^{t+\delta}Z_r^{1,u}dB_r,\\
Z_s^{1,u}&=&Z_s^u-D\varphi(s,X_s^u).\sigma(s,X_s^u,Y_s^u,u_s),\ \
t\leq s\leq t+\delta.  \end{array}\right.\ee From (4.5), (4.6),
(4.7), we have \be \esssup_{u\in
\mathcal{U}_{t,t+\delta}}Y_t^{1,u}\geq0, \ \mbox{P-a.s.}\ee
For $ (s,x,y,z,u) \in
[0,T]\times\mathbb{R}\times\mathbb{R}\times\mathbb{R}^d\times
 U,$\ we define\\
$
\begin{array}{lll}  L(s,x,y,z,u)&=&\frac{\partial}{\partial
s}\varphi(s,x)+ D\varphi(s,x).b(s, x, y+\varphi(s,x), z, u) +
\frac{1}{2}\tr(\sigma\sigma^{T}(s, x,
y+\varphi(s,x),u)D^2\varphi(s,x))\\ &&+
 f(s, x, y+\varphi(s,x), z, u),\\
 F(s,x,y,z,u)&=&
 L(s,x,y,z+D\varphi(s,x).\sigma(s,x,y+\varphi(s,x),u),u).
 \end{array}$\\

Then equation (4.8) can be reformulated as: \be\left
\{\begin{array}{llll}
dY_s^{1,u}&=&-F(s,X_s^u,Y_s^{1,u},Z_s^{1,u},u_s)ds+Z_s^{1,u}dB_s,\ \
s\in [t,t+\delta],\\
Y_{t+\delta}^{1,u}&=&0,
\end{array}\right.\ee
where $Y_s^{1,u}=Y_s^u-\varphi(s,X_s^u)$, $Z_s^{1,u}=Z_s^u-D\varphi(s,X_s^u).\sigma(s,X_s^u,Y_s^u,u_s)$.\\
Obviously, equation (4.10) has a unique solution
$(Y_s^{1,u},Z_s^{1,u})_{s\in [t,t+\delta]}\in \mathcal
{S}^2(t,t+\delta;\mathbb{R})\times \mathcal
{H}^2(t,t+\delta;\mathbb{R}^d).$ Indeed, equation (4.10) has a
solution
$(Y_s^u-\varphi(s,X_s^u),Z_s^u-D\varphi(s,X_s^u).\sigma(s,X_s^u,Y_s^u,u_s))_{s\in
[t,t+\delta]}.$ If (4.10) has another solution
$(\tilde{Y}_s^{1,u},\tilde{Z}_s^{1,u})_{s\in [t,t+\delta]}\in
\mathcal {S}^2(t,t+\delta;\mathbb{R})\times \mathcal
{H}^2(t,t+\delta;\mathbb{R}^d),$ then
$(X_s^u,\tilde{Y}_s^{1,u}+\varphi(s,X_s^u),\tilde{Z}_s^{1,u}+D\varphi(s,X_s^u).\sigma(s,X_s^u,Y_s^u,u_s))_{s\in
[t,t+\delta]}$ is the solution of equation (4.3), from the
uniqueness of the solution of FBSDE (4.3), we have
$\tilde{Y}_s^{1,u}+\varphi(s,X_s^u)=Y_s^u,\
\tilde{Z}_s^{1,u}+D\varphi(s,X_s^u).\sigma(s,X_s^u,Y_s^u,u_s)=Z_s^u,$
i.e., $\tilde{Y}_s^{1,u}=Y_s^u-\varphi(s,X_s^u),\ \mbox{P-a.s.},\
\tilde{Z}_s^{1,u}=Z_s^u-D\varphi(s,X_s^u).\sigma(s,X_s^u,Y_s^u,u_s),\
\mbox{a.s., a.e.}, \ s\in [t,t+\delta].$

Now we need to study the following BSDE:
 \be\left
\{\begin{array}{llll}
dY_s^{2,u}&=&-L(s,x,0,\hat{Z}_s^{u},u_s)ds+Z_s^{2,u}dB_s,\ \
s\in [t,t+\delta],\\
Y_{t+\delta}^{2,u}&=&0,
\end{array}\right.\ee
where
$\hat{Z}_s^{u}=Z_s^{1,u}+D\varphi(s,x).\sigma(s,x,Y_s^{1,u}+\varphi(s,x),u_s),\
s\in [t,t+\delta].$ Notice that \be
\begin{array}{llll}
\mbox{(i)}&&|L(s,x,y,z,u)-L(s,x,y',z',u)|\leq
C(1+|x|+|y|+|y'|)(|y-y'|+|z-z'|);\\
\mbox{(ii)}&&|L(s,x,y,z,u)|\leq C(1+|x|^2+|y|^2+|z|),\ \ \forall
(s,x,y,z,u)\in [0,T]\times
\mathbb{R}^n\times\mathbb{R}\times\mathbb{R}^d\times U.
\end{array}\ee
Since $\hat{Z}^u\in \mathcal {H}^2(t,t+\delta;\mathbb{R}^d)$, and
$E[\int_t^{t+\delta}|L(s,x,0,\hat{Z}_s^u,u_s)|^2ds]<+\infty,$
from Lemma 2.1 in~\cite{BL}, BSDE (4.11) has a unique solution
$(Y_s^{2,u},Z_s^{2,u})_{s\in [t,t+\delta]}\in \mathcal
{S}^2(t,t+\delta;\mathbb{R})\times\mathcal
{H}^2(t,t+\delta;\mathbb{R}^d).$\\
Now we prove some lemmas for the proof of Theorem 4.1.
 \bl For all
$u\in\mathcal {U}_{t,t+\delta}$, we have
$$E[\int_t^{t+\delta}(|Y_s^{1,u}|+|Z_s^{1,u}|)ds\mid \mathcal
{F}_t]\leq C\delta^{\frac{5}{4}},\ \mbox{P-a.s.},\ 0\leq \delta\leq
\delta_1,$$ where the constant $C$ is independent of the control $u$
and of $\delta>0.$ \el \noindent \textbf{Proof}. From (4.7) and (4.4),
 \be
\begin{array}{llll}
|Y_s^{1,u}|&=&|E[\int_s^{t+\delta}f(r,X_r^u,Y_r^u,Z_r^u,u_r)dr\mid\mathcal
{F}_s]+E[\varphi(t+\delta,X_{t+\delta}^u)-\varphi(s,X_s^u)\mid\mathcal
{F}_s]|\\
&\leq &
CE[\int_s^{t+\delta}(1+|X_r^u|+|Y_r^u|+|Z_r^u|)dr\mid\mathcal
{F}_s]+C\delta+CE[|X_{t+\delta}^u-X_s^u|\mid\mathcal {F}_s]\\
&\leq &
C\delta^\frac{1}{2}(E[\int_s^{t+\delta}(1+|X_r^u|^2+|Y_r^u|^2+|Z_r^u|^2)dr\mid\mathcal
{F}_s])^\frac{1}{2}+C\delta+C\delta^\frac{1}{2}(1+|X_s^u|)\\
&\leq &C\delta^\frac{1}{2}(1+|X_s^u|),\ \ \mbox{P-a.s.},\ s\in
[t,t+\delta],
\end{array}\ee
where notice that
$(X_s^u,Y_s^u,Z_s^u)=(X_s^{s,X_s^u},Y_s^{s,X_s^u},Z_s^{s,X_s^u}),\
\mbox{P-a.s.},\ s\in [t,t+\delta],$ from the uniqueness of the
solution of FBSDE (4.3) on $[t,t+\delta]$.\ Furthermore, from (4.7), we get $|Y_s^{u}|\leq C(1+|X_s^u|),\
\mbox{P-a.s.},\ s\in [t,t+\delta].$\ On the other hand,
$Z_s^{1,u}=Z_s^u-D\varphi(s,X_s^u).\sigma(s,X_s^u,Y_s^u,u_s)$, we
have
\be|Z_s^{1,u}|\leq C(1+|X_s^u|+|Z_s^u|),\ \
\mbox{P-a.s.},\ s\in [t,t+\delta].\ee
From (4.10), (4.12), (4.13), (4.14) and (4.4)
\be\begin{array}{llll}
&&|Y_t^{1,u}|^2+E[\int_t^{t+\delta}|Z_r^{1,u}|^2dr\mid\mathcal
{F}_t]\\
&=&2E[\int_t^{t+\delta}Y_r^{1,u}F(r,X_r^{u},Y_r^{1,u},Z_r^{1,u},u_r)dr\mid\mathcal
{F}_t]\\
&\leq&C\delta^\frac{1}{2}E[\int_t^{t+\delta}(1+|X_r^{u}|^2+|X_r^{u}|^3)dr\mid\mathcal
{F}_t]+C\delta^\frac{1}{2}E[\int_t^{t+\delta}|Z_r^{u}|^2dr\mid\mathcal
{F}_t]\\
&\leq&C\delta^\frac{3}{2},\ \mbox{P-a.s.}
\end{array}\ee
Therefore, from (4.13), (4.15) and (4.4)$$\begin{array}{llll}
&&E[\int_t^{t+\delta}(|Y_s^{1,u}|+|Z_s^{1,u}|)ds\mid\mathcal
{F}_t]\\
&\leq&C\delta^\frac{1}{2}E[\int_t^{t+\delta}(1+|X_r^{u}|)dr\mid\mathcal
{F}_t]+C\delta^\frac{1}{2}(E[\int_t^{t+\delta}|Z_r^{1,u}|^2dr\mid\mathcal
{F}_t])^\frac{1}{2}\\
&\leq&C\delta^\frac{5}{4},\ \mbox{P-a.s.},\ 0\leq \delta \leq
\delta_1.
\end{array}$$ \endpf
\br From (4.4), and (4.14), \be E[(\int_t^{t+\delta}|Z_r^{1,u}|^2dr)^2\mid\mathcal
{F}_t]\leq C\delta^2.\ee\er
\bl  For all $u\in\mathcal {U}_{t,t+\delta}$, we have
$$|Y_t^{1,u}-Y_t^{2,u}|\leq C\delta^{\frac{5}{4}},\ \
\mbox{P-a.s.},\ 0\leq \delta\leq \delta_1,$$ where
 $C$ is independent of the control process $u$ and of $\delta>0$.
\el
  \noindent \textbf{Proof}. We define: $g(s)=L(s,X_s^{u},0,Z_s^{u},u_s)-L(s,x,0,Z_s^{u},u_s),\ \rho_0(r)=(1+|x|^2+|Z_r^{u}|)(r+r^2),\
  r\geq0,$ Then, $|g(s)|\leq C\rho_0(|X_s^{u}-x|),\ s\in
[t,t+\delta].$ From (4.4), (4.10), (4.11), the definitions of $F$ and
$L$, we have
 \be\begin{array}{llll}
 |Y_t^{1,u}-Y_t^{2,u}|&=&|E[(Y_t^{1,u}-Y_t^{2,u})|\mathcal
{F}_t]|\\
&=&|E[\int_t^{t+\delta}(L(s,X_s^{u},Y_s^{1,u},Z_s^{u},u_s)-L(s,x,0,\hat{Z}_s^{u},u_s))ds\mid\mathcal
{F}_t]|\\
&\leq&
E[\int_t^{t+\delta}[C(1+|X_s^{u}|+|Y_s^{1,u}|)|Y_s^{1,u}|+C(1+|x|)|Z_s^{u}-\hat{Z}_s^{u}|+g(s)]ds\mid\mathcal
{F}_t].\\
\end{array}\ee
Notice that from (4.13), $E[\int_t^{t+\delta}(1+|X_s^{u}|+|Y_s^{1,u}|)|Y_s^{1,u}|ds\mid\mathcal{F}_t]\leq C\delta^\frac{3}{2};$
from (4.8), (4.11) and (4.4)
 $$\begin{array}{llll}
&&
E[\int_t^{t+\delta}|Z_s^{u}-\hat{Z}_s^{u}|ds\mid\mathcal
{F}_t]\\
&=&E[\int_t^{t+\delta}|Z_s^{1,u}+D\varphi(s,X_s^u).\sigma(s,X_s^u,Y_s^{1,u}+\varphi(s,X_s^u),u_s)\\
 & & \ \ \ \  \ \ \ -(Z_s^{1,u}+D\varphi(s,x).\sigma(s,x,Y_s^{1,u}+\varphi(s,x),u_s))|ds\mid\mathcal
{F}_t]\\
&\leq&
CE[\int_t^{t+\delta}|X_s^{u}-x|(1+|X_s^{u}|+|Y_s^{1,u}|)ds\mid\mathcal
{F}_t]\\
&\leq& C\delta E[\sup\limits_{t\leq s\leq
t+\delta}|X_s^{u}-x|(1+|X_s^{u}|+|Y_s^{1,u}|)|\mathcal {F}_t]\leq
C\delta^\frac{3}{2},\ \mbox{P-a.s.};\end{array}$$ furthermore, from (4.4),\\
$$\begin{array}{llll}& & E[\int_t^{t+\delta}g(s)ds\mid\mathcal {F}_t]
 \leq \delta^\frac{1}{2}(E[\int_t^{t+\delta}g(s)^2ds\mid\mathcal
{F}_t])^\frac{1}{2}\\
&\leq& C
\delta^\frac{1}{2}(E[\int_t^{t+\delta}(1+|x|^4+|Z_r^{u}|^2)(|X_r^{u}-x|^2+|X_r^{u}-x|^4)dr\mid\mathcal
{F}_t])^\frac{1}{2}
 \leq  C\delta^\frac{5}{4},\ \ \mbox{P-a.s.}
\end{array}$$
From (4.17), the proof is complete. \endpf

 Now we consider the following equation:
\be
 \left \{\begin{array}{llll}
dY_s^{3,u}&=&-L(s,x,0,D\varphi(s,x).\sigma(s,x,\varphi(s,x),u_s),u_s)ds+Z_s^{3,u}dB_s,\ \  s\in[t,t+\delta],  \\
 Y_{t+\delta}^{3,u} &= &0,   \end{array}\right. \ee
where $u(\cdot)\in\mathcal {U}_{t,t+\delta}$. Notice that
$$L(s,x,0,D\varphi(s,x).\sigma(s,x,\varphi(s,x),u_s),u_s)=F(s,x,0,0,u_s).$$

\bl  For all $u\in\mathcal {U}_{t,t+\delta}$, we have
$$|Y_t^{2,u}-Y_t^{3,u}|\leq C\delta^{\frac{5}{4}},\ \
\mbox{P-a.s.},\ 0\leq \delta\leq \delta_1,$$ where
 $C$ is independent of the control process $u$ and of $\delta>0$.
\el
  \noindent \textbf{Proof}. From (4.11), (4.18), (4.12) and Lemma 4.1,
$$\begin{array}{llll}
 |Y_t^{2,u}-Y_t^{3,u}|
&=&|E[\int_t^{t+\delta}(L(s,x,0,\hat{Z}_s^{u},u_s)-L(s,x,0,D\varphi(s,x).\sigma(s,x,\varphi(s,x),u_s),u_s))ds\mid\mathcal
{F}_t]|\\
&\leq&
CE[\int_t^{t+\delta}(1+|x|)(|Y_s^{1,u}|+|Z_s^{1,u}|)ds\mid\mathcal
{F}_t]\\
&\leq&C\delta^\frac{5}{4},\ \ \mbox{P-a.s.},\  0\leq \delta\leq
\delta_1.
\end{array}$$ \endpf
\bl  Let  $Y_0(\cdot)$ be the solution of the following ordinary
differential equation: \be
 \left \{\begin{array}{llll}
 \dot{Y}_0(s) &=& -F_0(s,x,0,0),\  s\in[t,t+\delta],\\
Y_0(t+\delta) &=& 0,  \end{array}\right. \ee
 where \be F_0(s,x,0,0)=\sup\limits_{u\in U}F(s,x,0,0,u). \ee
 Then, P-a.s.,\be \esssup_{u\in\mathcal {U}_{t,t+\delta}}Y_t^{3,u}=Y_0(t).\ee \el
 \noindent \textbf{Proof}. Obviously, (4.19) has a unique solution. From the definition of $F_0(s,x,0,0)$, we know
$$F_0(s,x,0,0)\geq F(s,x,0,0,u_s), \ \ \  \mbox{for any} \ u \in \mathcal {U}_{t,t+\delta}.$$
Therefore, from the comparison theorem of BSDE (see, Lemma 2.2 in~\cite{BL}),
$$ \tilde{Y_0}(s) \geq Y_s^{3,u},\ \mbox{P-a.s.},
\ \ \mbox{for any}\ s \in [t,t+\delta],\ \mbox{for any} \ u \in
\mathcal {U}_{t,t+\delta},$$ where
$(\tilde{Y_0}(\cdot),\tilde{Z_0}(\cdot))$ is the solution of the
following BSDE:
$$ \left \{\begin{array}{llll}
  d\tilde{Y_0}(s) &=& -F_0(s,x,0,0)ds+\tilde{Z_0}(s)dB_s, \ s\in[t,t+\delta],\\
 \tilde{Y_0}(t+\delta) &=& 0.  \end{array}\right. $$
In fact, $(\tilde{Y_0}(s),\tilde{Z_0}(s))=(Y_0(s),0)$. Therefore,
$Y_0(t)\geq Y_t^{3,u}, \mbox{P-a.s.}, \mbox{ for any} \ u\in
\mathcal {U}_{t,t+\delta}$.

 On the other hand, since
$F_0(s,x,0,0)=\sup\limits_{u\in U}F(s,x,0,0,u)$, there exists some
measurable function $\tilde{u}(s,x): [t,t+\delta]\times\mathbb{R}^n
\rightarrow U$, such that $F(s,x,0,0,\tilde{u}(s,x))=F_0(s,x,0,0)$.
Define $\tilde{u}_s^0=\tilde{u}(s,x),\ s\in[t,t+\delta],$ then
$\tilde{u}^0\in \mathcal {U}_{t,t+\delta}$ and
$F_0(s,x,0,0)=F(s,x,0,0,\tilde{u}_s^0),\ s\in[t,t+\delta].$
Consequently, from the uniqueness of the solution of the BSDE it
follows that $Y_0(t)=Y_t^{3,\tilde{u}^0},\ \mbox{P-a.s.}$ \
Therefore, $\esssup_{u\in\mathcal {U}_{t,t+\delta}}Y_t^{3,u}=Y_0(t),\ \ \mbox{P-a.s.}$ \endpf

 Now we are able to complete the proof of Theorem 4.1 as follows:

Indeed, from (4.9) we know that $\esssup_{u\in\mathcal
{U}_{t,t+\delta}}Y_t^{1,u}\geq 0, \ \ \mbox{P-a.s.}$\ Therefore,
from the Lemmas 4.2 and 4.3 we get $\esssup_{u\in\mathcal
{U}_{t,t+\delta}}Y_t^{3,u}\geq -C\delta^\frac{5}{4}, \ \
\mbox{P-a.s.}$\ Thus, from Lemma 4.4, $Y_0(t)\geq
-C\delta^\frac{5}{4},$ where $Y_0$ is the solution of (4.19). Then,
$$\frac{1}{\delta}\int_t^{t+\delta}F_0(s,x,0,0)ds\geq -C\delta^\frac{1}{4},\ \ 0\leq \delta\leq \delta_1.$$
It follows by letting $\delta \rightarrow 0$ that
$$\sup_{u\in U}F(t,x,0,0,u)=F_0(t,x,0,0)\geq 0.$$
From the definition of $F$ we see that $W$ is a subsolution of
(4.2). Similarly, we can prove that $W$ is a viscosity supersolution
of (4.2). Therefore, $W$ is a viscosity solution of (4.2).\\
\endpf

 \textbf{Case 2}. We suppose that $\sigma$ depends on $z$, and does not depend on $u$.

 Now equation (3.1) becomes the following one \be
  \left \{
  \begin{array}{ll}
  dX_s^{t,x;u}  =  b(s,X_s^{t,x;u},Y_s^{t,x;u},Z_s^{t,x;u},u_s)ds +
          \sigma(s,X_s^{t,x;u},Y_s^{t,x;u},Z_s^{t,x;u}) dB_s, \\
  dY_s^{t,x;u}  =  -f(s,X_s^{t,x;u},Y_s^{t,x;u},Z_s^{t,x;u},u_s)ds + Z_s^{t,x;u}dB_s, \ \ \ \  s\in [t,T],\\
   X_t^{t,x;u} =  x,\
   Y_T^{t,x;u}  =  \Phi(X_T^{t,x;u}).
   \end{array}
   \right.
  \ee
The related HJB equation is the following PDE combined with the
algebraic equation: \be
 \left \{\begin{array}{ll}
 &\!\!\!\!\! \frac{\partial }{\partial t} W(t,x) +  H(t, x, W(t,x),
 V(t,x))=0,
   \\
&\!\!\!\!\!V(t,x)=DW(t,x).\sigma(t,x,W(t,x),V(t,x)),\hskip 0.5cm   (t,x)\in [0,T)\times {\mathbb{R}}^n ,\\
 &\!\!\!\!\!  W(T,x) =\Phi (x),\ \ \ \ x\in \mathbb{R}^n.
 \end{array}\right.
\ee
 In this case
$$\begin{array}{lll}
 H(t, x, W(t,x), V(t,x))= \sup\limits_{u \in U}\{ DW(t,x).b(t, x,
W(t,x), V(t,x), u)\\
 \hskip 4cm +\frac{1}{2}\tr(\sigma\sigma^{T}(t, x, W(t,x),
V(t,x))D^2W(t,x))+f(t, x, W(t,x), V(t,x), u)\},
 \end{array}$$
where $t\in [0, T], x\in{\mathbb{R}}^n$.

We also give the definition of viscosity solution for this kind of
PDE.

\bde\mbox{ }A real-valued
continuous function $W\in C([0,T]\times {\mathbb{R}}^n )$ is called \\
  {\rm(i)} a viscosity subsolution of equation (4.23) if $W(T,x) \leq \Phi (x),\ \mbox{for all}\ x \in
  {\mathbb{R}}^n$, and if for all functions $\varphi \in C^3_{l, b}([0,T]\times
  {\mathbb{R}}^n)$ satisfying the monotonicity condition (B2)'-(ii) and for all $(t,x) \in [0,T) \times {\mathbb{R}}^n$ such that $W-\varphi $\ attains
  a local maximum at $(t, x)$,
 $$\left \{\begin{array}{ll}
 &\!\!\!\!\! \frac{\partial \varphi}{\partial t} (t,x) + H(t,x,\varphi(t,x),\psi(t,x)) \geq
 0,\\
&\!\!\!\!\!\mbox{where }\psi \mbox{ is the unique solution of the
following algebraic
equation:}\\
&\!\!\!\!\!\psi(t,x)=D\varphi(t,x).\sigma(t,x,\varphi(t,x),\psi(t,x)).
\end{array}\right.
$$
{\rm(ii)} a viscosity supersolution of equation (4.23) if $W(T,x)
\geq \Phi (x), \mbox{for all}\ x \in
  {\mathbb{R}}^n$, and if for all functions $\varphi \in C^3_{l, b}([0,T]\times
  {\mathbb{R}}^n)$ satisfying the monotonicity condition (B2)'-(ii) and for all $(t,x) \in [0,T) \times {\mathbb{R}}^n$\ such that $W-\varphi $\ attains
  a local minimum at $(t, x)$,
$$\left \{\begin{array}{ll}
 &\!\!\!\!\! \frac{\partial \varphi}{\partial t} (t,x) + H(t,x,\varphi(t,x),\psi(t,x)) \leq
 0,\\
 &\!\!\!\!\!\mbox{where }\psi \mbox{ is the unique solution of the following algebraic
equation:}\\
&\!\!\!\!\!
\psi(t,x)=D\varphi(t,x).\sigma(t,x,\varphi(t,x),\psi(t,x)).
\end{array}\right.
$$
{\rm(iii)} a viscosity solution of equation (4.23) if it is both a
viscosity sub- and supersolution of equation (4.23).
\ede

\br In this case we need the following technical assumption:
\\
(B6) $\beta_2>0$;\\
(B7) $G\sigma(s, x, y, z)$ is continuous in $s$, uniformly with
respect to $(x, y, z)\in {\mathbb{R}}^n\times{\mathbb{R}}\times
{\mathbb{R}}^d$.\er
 \bt  Under
the assumptions (B2), (B4), (B5), (B6) and (B7), the value function $W$
is a viscosity solution of (4.23).\et

 \noindent \textbf{Proof}. Obviously, $W(T,x)=\Phi(x),\ x\in\mathbb{R}^n.$ We prove only that $W$ is a viscosity subsolution,
 that it is also a viscosity supersolution can be
proved similarly. We suppose that $\varphi \in
  C_{l,b}^3([0,T]\times\mathbb{R}^n)$ satisfying the monotonicity condition (B2)'-(ii) and that
  $(t,x)\in[0,T)\times\mathbb{R}^n$ is such that $W-\varphi$
  attains its maximum at $(t,x)$. Since $W$ is continuous
  and of at most linear growth, we can replace the
  condition of a local maximum by that of a global one in the
  definition of the viscosity subsolution. Without loss of generality we may assume that $\varphi(t,x)=W(t,x).$
  We consider the following equation: \be
  \left \{
  \begin{array}{ll}
  d\overline{X}_s^{u}  = b(s,\overline{X}_s^{u},\overline{Y}_s^{u},\overline{Z}_s^{u},u_s)ds +
          \sigma(s,\overline{X}_s^{u},\overline{Y}_s^{u},\overline{Z}_s^{u})dB_s, \\
  d\overline{Y}_s^{u}  =  -f(s,\overline{X}_s^{u},\overline{Y}_s^{u},\overline{Z}_s^{u},u_s)ds + \overline{Z}_s^{u}dB_s, \ \ s\in[t,t+\delta], \\
   \overline{X}_t^{u} =  x,\
   \overline{Y}_{t+\delta}^{u}  =
   \varphi(t+\delta,\overline{X}_{t+\delta}^{u}),\ \ 0\leq\delta\leq
   T-t.
   \end{array}
   \right.\ee
From Proposition \ref{App-Pro4} and Proposition \ref{App-Pro5} in
Appendix, we know there exists sufficiently small
$0<{\bar{\delta}}_1<T-t$ such that for any $0\leq \delta\leq
{\bar{\delta}}_1$, FBSDE (4.24) has a unique solution
$(\overline{X}_s^{u},\overline{Y}_s^{u},\overline{Z}_s^{u})_{s\in[t,t+{\delta}]}\in\mathcal
{S}^2(t,t+ {\delta};\mathbb{R}^n)\times \mathcal {S}^2(t,t+
{\delta};\mathbb{R})\times\mathcal {H}^2(t,t+
{\delta};\mathbb{R}^d)$, and for $p\geq 2,$
 \be\begin{array}{llll}
\mbox{(i)}&&E[\sup\limits_{t\leq s\leq
t+{\delta}}|\overline{X}_s^{u}|^p+\sup\limits_{t\leq s\leq
t+{\delta}}|\overline{Y}_s^{u}|^p+(\int_t^{t+{\delta}}|\overline{Z}_s^{u}|^2ds)^{\frac{p}{2}}\mid\mathcal
{F}_t]\leq C(1+|x|^p),\ \ \mbox{P-a.s.};\\
\mbox{(ii)}&&E[\sup\limits_{t\leq s\leq
t+\delta}|\overline{X}_s^{u}-x|^p\mid\mathcal {F}_t]\leq
C\delta^{\frac{p}{2}}(1+|x|^p),\ \ \mbox{P-a.s.};\\
\mbox{(iii)}&&E[(\int_t^{t+{\delta}}|\overline{Z}_s^{u}|^2ds)^{\frac{p}{2}}\mid\mathcal
{F}_t]\leq C\delta^{\frac{p}{2}}(1+|x|^p),\ \ \mbox{P-a.s.}
 \end{array}\ee
According to the definition of the backward stochastic semigroup for
fully coupled FBSDE, we have
   $$G_{s,t+\delta}^{t,x;u}[\varphi(t+\delta,\overline{X}_{t+\delta}^u)]=\overline{Y}_s^u,\ \ s\in [t,t+\delta].$$
And due to the DPP (Theorem 3.1), we have $$
\varphi(t,x)=W(t,x)=\esssup_{u\in\mathcal
{U}_{t,t+\delta}}G_{t,t+\delta}^{t,x;u}[W(t+\delta,\widetilde{X}_{t+\delta}^{t,x;u})],\
0\leq \delta \leq  {\bar{\delta}}_1,$$
where $\widetilde{X}^{t,x;u}$\ is defined by FBSDE (3.10).

From $\varphi(s,y) \geq W(s,y), \ (s,y) \in [0,T) \times
\mathbb{R}^n,$ and the monotonicity property of
$G_{t,t+\delta}^{t,x;u}[\cdot]$ (see Theorem \ref{App-Th2}) we
obtain \be \esssup_{u\in\mathcal
{U}_{t,t+\delta}}{G_{t,t+\delta}^{t,x;u}[\varphi(t+\delta,\overline{X}_{t+\delta}^u)]-\varphi(t,x)}\geq
0,\ 0\leq \delta \leq {\bar{\delta}}_1.\ee  Now we set \be
\begin{array}{llll}
\overline{Y}^{1,u}_s &=& \overline{Y}^{u}_s-\varphi(s,\overline{X}^{u}_s) \\
&=&\int_s^{t+\delta}f(r,\overline{X}_r^u,\overline{Y}_r^u,\overline{Z}_r^u,u_r)dr-\int_s^{t+\delta}\overline{Z}_r^udB_r+\varphi(t+\delta,\overline{X}_{t+\delta}^u)-\varphi(s,\overline{X}_s^u).
\end{array}
\ee
 Applying It\^o's formula to $\varphi(s,\overline{X}_s^u),$
and setting
$\overline{Z}_s^{1,u}=\overline{Z}_s^u-D\varphi(s,\overline{X}_s^u).\sigma(s,\overline{X}_s^u,\overline{Y}_s^u,\overline{Z}_s^u)$,
we obtain \be
  \left \{
  \begin{array}{llll}
  \overline{Y}_s^{1,u} & = &\int_s^{t+\delta}[\frac{\partial}{\partial r}\varphi(r, \overline{X}_r^u)
  +\frac{1}{2}\tr(\sigma
  \sigma^T(r,\overline{X}_r^u,\overline{Y}_r^u,\overline{Z}_r^u)D^2\varphi(r,
  \overline{X}_r^u))\\
  & &+D\varphi(r,\overline{X}_r^u).b(r,\overline{X}_r^u,\overline{Y}_r^u,\overline{Z}_r^u,u_r)+f(r,\overline{X}_r^u,\overline{Y}_r^u,\overline{Z}_r^u,u_r)]dr
          -\int_s^{t+\delta}\overline{Z}_r^{1,u}dB_r, \\
  \overline{Z}_s^{1,u}&=&\overline{Z}_s^u-D\varphi(s,\overline{X}_s^u).\sigma(s,\overline{X}_s^u,\overline{Y}_s^u,\overline{Z}_s^u).
   \end{array}
   \right.
  \ee
From (4.26) and (4.27), we have \be \esssup_{u\in\mathcal
{U}_{t,t+\delta}}\overline{Y}_t^{1,u}\geq 0.\ee Define\\
$\begin{array}{llll}L(s,x,y,z,u)&=&\frac{\partial}{\partial
s}\varphi(s, x)
  +D\varphi(s,x).b(s,x,y+\varphi(s,x),z,u)+\frac{1}{2}\tr(\sigma
  \sigma^T(s,x,y+\varphi(s,x),z)D^2\varphi(s,
  x))\\&&+f(s,x,y+\varphi(s,x),z,u),\ (s, x, y, z, u)\in [0, T]\times \mathbb{R}^n\times \mathbb{R}\times \mathbb{R}^d\times U.\end{array}$\\
Notice now that
$$
\begin{array}{llll}
\mbox{(i)}&&|L(s,x,y,z,u)-L(s,x,y',z',u)|\leq
C(1+|x|+|y|+|y'|+|z|+|z'|)(|y-y'|+|z-z'|);\\
\mbox{(ii)}&&|L(s,x,y,z,u)|\leq C(1+|x|^2+|y|^2+|z|^2),\ \ \forall
(s,x,y,z,u)\in [0,T]\times
\mathbb{R}^n\times\mathbb{R}\times\mathbb{R}^d\times U.
\end{array}$$
Therefore, equation (4.28) can be written into the following form:
 \be
  \left \{
  \begin{array}{llll}
  d\overline{Y}_s^{1,u} & =&-L(s,\overline{X}_s^u,\overline{Y}_s^{1,u},\overline{Z}_s^{u},u_s)ds+\overline{Z}_s^{1,u}dB_s, \ \ s\in[t,t+\delta], \\
 \overline{Z}_s^{u}&=&\overline{Z}_s^{1,u}+D\varphi(s,\overline{X}_s^u).\sigma(s,\overline{X}_s^u,\overline{Y}_s^{1,u}+\varphi(s,\overline{X}_s^u),\overline{Z}_s^u), \ \ s\in[t,t+\delta],\\
\overline{Y}_{t+\delta}^{1,u}&=&0,
   \end{array}
   \right.
  \ee
  Obviously, (4.30) has a solution
  $(\overline{Y}_s^{1,u},\overline{Z}_s^{1,u})\in\mathcal
{S}^2(t,t+\bar{\delta}_1;\mathbb{R})\times\mathcal
{H}^2(t,t+\bar{\delta}_1;\mathbb{R}^d)$, because (4.24) has
  a unique solution
  $(\overline{X}_s^u,\overline{Y}_s^u,\overline{Z}_s^u)_{s\in [ t, t+{\delta}]}$, and
  $\overline{Y}_s^{1,u}=\overline{Y}_s^u-\varphi(s,\overline{X}_s^u),$
   $\overline{Z}_s^{1,u}=
   \overline{Z}_s^u-D\varphi(s,\overline{X}_s^u).\sigma(s,\overline{X}_s^u,\overline{Y}_s^{1,u}
   +\varphi(s,\overline{X}_s^u),\overline{Z}_s^u)$\ solves (4.30).

In order to complete the proof of Theorem 4.2, we need the following
lemmas.

We need to consider the following BSDE combined with an algebraic
equation: \be
  \left \{
  \begin{array}{llll}
  d\overline{Y}_s^{2,u} & = & -L(s,x,0,\hat{Z}_s^{u},u_s)ds+
  \overline{Z}_s^{2,u}dB_s,
\ \ \ s\in[t,t+\delta],\\
 \hat{Z}_s^{u} & = &\overline{Z}_s^{1,u} + D\varphi(s,x).\sigma(s,x,\overline{Y}_s^{1,u}+\varphi(s,x),\hat{Z}_s^{u}),\ \ \ s\in[t,t+\delta],\\
   \overline{Y}_{t+\delta}^{2,u} & = & 0,
   \end{array}
   \right.\ee
where $u(\cdot) \in \mathcal {U}_{t,t+\delta}.$ For it, we have to
study first the algebraic equation.

\br For $m=1,$ the  matrix $G$ becomes a vector in $\mathbb{R}^n$,
and without loss of generality, we may assume
$G=(1,0,\cdots,0)\in\mathbb{R}^n$. Thus, we have the following
conditions from the monotonicity condition (B2):
 \be
\begin{array}{llll}{\rm(i)}& & \langle
\sigma_1(t,x,y,z)-\sigma_1(t,x,y,\bar{z}),z-\bar{z}\rangle \leq
-\beta_2|z-\bar{z}|^2;\\
{\rm(ii)}& & G^TD\varphi(s,x) \geq 0,\ i.e., D_{x_1}\varphi(s,x)\geq
0,\ D_{x_i}\varphi(s,x)= 0,\ 2\leq i \leq n.\end{array} \ee Indeed,
(i) follows from (B2)-(i), and (ii) follows from (B2)'-(ii)
satisfied by $\varphi: \langle
\varphi(s,x)-\varphi(s,\bar{x}),G(x-\bar{x})\rangle\geq 0,\
\varphi\in C_{l,b}^3([0,T]\times\mathbb{R}^n)$. \er

We have the following important Representation Theorem for the
solution of the algebraic equation.

\bp For any $s\in [0, T],\ \zeta\in\mathbb{R}^d,\ y\in \mathbb{R},\
\bar{x}\in \mathbb{R}^n$, there exists a unique $z$ such that
$z=\zeta+D\varphi(s,\bar{x}).\sigma(s,\bar{x},y+\varphi(s,\bar{x}),z).
$ That means, the solution $z$ can be written as $z=h(s, \bar{x},y,
\zeta)$, where the function $h$ is Lipschitz with respect to $y,
\zeta,$ and $|h(s, \bar{x},y, \zeta)|\leq
C(1+|\bar{x}|+|y|+|\zeta|).$ The constant $C$\ is independent of $s,
\bar{x}, y, \zeta$. And $z=h(s, \bar{x},y, \zeta)$ is continuous
with respect to $s$. \ep

 \noindent \textbf{Proof}. \emph{First Step.} From Remark 4.4, we can prove that the equation
$z=D\varphi(s,\bar{x}).\sigma(s,\bar{x},\varphi(s,\bar{x}),z)$ has a
unique
solution $z$.\\
Indeed, by fixing $(s,\bar{x})$ and setting
$a:=D_{x_1}\varphi(s,\bar{x}),\
\sigma_1(z)=\sigma_1(s,\bar{x},\varphi(s,\bar{x}),z)$, we only need
to consider the equation $z=a\sigma_1(z)$.

(1) If $a=0,$ then $z=0$ is the solution.

(2) Let $a>0,$ and $ \zeta\in\mathbb{R}^d.$ We choose a small
$\delta\in (0, 1)$\ to make the mapping $z\mapsto
a\delta\sigma_1(z)$ Lipschitz with Lipschitz constant $L_\delta<1,$
and we set $\delta_0=\frac{1}{1+[\frac{1}{\delta}]}$, where $[z]$
represents the integer part of the real nonnegative number $z$. Then
the mapping $z\mapsto a\delta_0\sigma_1(z)$ is still Lipschitz with
Lipschitz constant $L_{\delta_0}<1.$ For simplicity of the notations
we still use $\delta$ to denote $\delta_0$. Obviously, there exists
a unique fixed point $z_\delta$ such that $z_\delta=\zeta+\delta
a\sigma_1(z_\delta).$ We now consider $z_\delta^{n+1}=(\zeta+\delta
a\sigma_1(z^n_\delta))+\delta a\sigma_1(z^{n+1}_\delta),\ n\geq1,\ z^{1}_\delta=0,$\
and we put $ \bar{z}_\delta^n:=z_\delta^{n}-z_\delta^{n-1},\ n>1.$
Then
$$\begin{array}{llll} |\bar{z}_\delta^{n+1}|^2&=&\delta
a\langle\sigma_1(z^n_\delta)-\sigma_1(z^{n-1}_\delta),\bar{z}^{n+1}_\delta\rangle+\delta
a\langle\sigma_1(z^{n+1}_\delta)-\sigma_1(z^{n}_\delta),\bar{z}^{n+1}_\delta\rangle\\
&\leq& L_\delta |\bar{z}^{n}_\delta|\cdot
|\bar{z}^{n+1}_\delta|-\beta_2\delta
a|\bar{z}^{n+1}_\delta|^2.\end{array}$$ Thus, putting
$\bar{\beta}_2:=a\delta\beta_2>0,$ we have
$(1+\bar{\beta}_2)|\bar{z}^{n+1}_\delta|^2\leq
|\bar{z}^{n}_\delta|\cdot|\bar{z}^{n+1}_\delta|\leq
\frac{1}{2}|\bar{z}^{n}_\delta|^2+\frac{1}{2}|\bar{z}^{n+1}_\delta|^2,$
from where we get
$(\frac{1}{2}+\bar{\beta}_2)|\bar{z}^{n+1}_\delta|^2\leq
\frac{1}{2}|\bar{z}^{n}_\delta|^2,\ n\geq 1.$
Therefore, there exists a
unique $z_\delta^{'}$, such that
$z_\delta^{'}=\zeta+2\delta a \sigma_1(z_\delta^{'}).$\\
Let $N\geq 1$ such that $N\delta=1$ and $1\leq k\leq N$. Suppose
that $\zeta\in\mathbb{R}^d$, there exists a unique $z_\delta$ such
that $z_\delta=\zeta+k\delta a\sigma_1(z_\delta)$. Now we consider the
equation $z_\delta^{n+1}=(\zeta+\delta
a\sigma_1(z^n_\delta))+k\delta a\sigma_1(z^{n+1}_\delta),\ n\geq1,\ z^{1}_\delta=0.$
Using the above argument we see that also the equation
$z_\delta=\zeta+(k+1)\delta a \sigma_1(z_\delta)$ has a unique fixed
point $z_\delta,$ for all $\zeta\in\mathbb{R}^d$. This completes the
proof of the first step.

 \emph{Second step.} From the above, since for any $\zeta\in\mathbb{R}^d$, there
exists a unique $z$ of the equation
$z=\zeta+D\varphi(s,\bar{x}).\sigma(s,\bar{x},y+\varphi(s,\bar{x}),z)(=\zeta
+D_{x_1}\varphi(s,\bar{x})\sigma_1(s,\bar{x},y+\varphi(s,\bar{x}),z)).$
$z$ is uniquely determined by $(s, \bar{x},y, \zeta)$, and we can
put $z=h(s, \bar{x},y, \zeta)$. This function $h$\ is measurable and
it is Lipschitz with respect to $y, \ \zeta.$

 Indeed, for any $\bar{y}, \ \bar{\zeta},\ \hat{y}, \ \hat{\zeta},$ we consider: $$\bar{z}=\bar{\zeta}+D_{x_1}\varphi(s,\bar{x})\sigma_1(s,\bar{x},\bar{y}+\varphi(s,\bar{x}),\bar{z}),\ \hat{z}=\hat{\zeta}+D_{x_1}\varphi(s,\bar{x})\sigma_1(s,\bar{x},\hat{y}+\varphi(s,\bar{x}),\hat{z}).$$
Then, taking into account Remark 4.4-(i) and $\varphi\in
C^3_{l,b}([0,T]\times\mathbb{R}^n),$
\\$\begin{array}{llll}
&&\langle\bar{z}-\hat{z},\bar{z}-\hat{z}\rangle\\
&=&\langle \bar{\zeta}-\hat{\zeta}+
D_{x_1}\varphi(s,\bar{x})\sigma_1(s,\bar{x},\bar{y}+\varphi(s,\bar{x}),\bar{z})
-D_{x_1}\varphi(s,\bar{x})\sigma_1(s,\bar{x},\hat{y}+\varphi(s,\bar{x}),\hat{z}),
\bar{z}-\hat{z}\rangle\\
&=&\langle \bar{\zeta}-\hat{\zeta},
\bar{z}-\hat{z}\rangle + \langle D_{x_1}\varphi(s,\bar{x})\sigma_1(s,\bar{x},\bar{y}+\varphi(s,\bar{x}),\bar{z})-
D_{x_1}\varphi(s,\bar{x})\sigma_1(s,\bar{x},\bar{y}+\varphi(s,\bar{x}),\hat{z}), \bar{z}-\hat{z}\rangle\\
& &+\langle D_{x_1}\varphi(s,\bar{x})\sigma_1(s,\bar{x},\bar{y}+\varphi(s,\bar{x}),\hat{z})
-D_{x_1}\varphi(s,\bar{x})\sigma_1(s,\bar{x},\hat{y}+\varphi(s,\bar{x}),\hat{z}), \bar{z}-\hat{z}\rangle\\
&\leq& C |\bar{\zeta}-\hat{\zeta}|^2+\frac{1}{8} |\bar{z}-\hat{z}|^2+C(-\beta_2)|\bar{z}-\hat{z}|^2
 +C |\bar{y}-\hat{y}||\bar{z}-\hat{z}| \\
 &\leq& C |\bar{\zeta}-\hat{\zeta}|^2+\frac{1}{4} |\bar{z}-\hat{z}|^2+C|\bar{y}-\hat{y}|^2.\end{array}$\\
Therefore, we have $|\bar{z}-\hat{z}|\leq
C(|\bar{\zeta}-\hat{\zeta}|+|\bar{y}-\hat{y}|).$

Similarly we can
prove $|h(s, \bar{x},y, \zeta)|\leq C(1+|\bar{x}|+|y|+|\zeta|),$
where the constant $C$\ is independent of  $s,\ \bar{x},\ y,\
\zeta$. Indeed, we have
\\$\begin{array}{llll}
&&\langle z,z\rangle=\langle  \zeta+
D_{x_1}\varphi(s,\bar{x})\sigma_1(s,\bar{x},y+\varphi(s,\bar{x}),{z})
, z\rangle\\
&=&\langle  {\zeta} ,
z\rangle + \langle D_{x_1}\varphi(s,\bar{x})\sigma_1(s,\bar{x}, {y}+\varphi(s,\bar{x}),{z})-
D_{x_1}\varphi(s,\bar{x})\sigma_1(s,\bar{x}, {y}+\varphi(s,\bar{x}),0), z\rangle\\
& &+\langle D_{x_1}\varphi(s,\bar{x})\sigma_1(s,\bar{x}, {y}+\varphi(s,\bar{x}),0), z\rangle\\
&\leq& C |{\zeta}|^2+\frac{1}{4} |z|^2+C(-\beta_2)|z|^2
 +C(1+|\bar{x}|+| {y}|)^2,
 \end{array}$\\
which implies $|z|\leq C(1+|\bar{x}|+|y|+|\zeta|).$

The fact that $z=h(s, \bar{x}, y, \zeta)$ is continuous with respect
to $s$ can be proved similarly.\\ Indeed, let
$ z_1=\zeta+D_{x_1}\varphi(s_1,x)\sigma_1(s_1,x,y+\varphi(s_1,x),z_1),\ z_2=\zeta+D_{x_1}\varphi(s_2,x)\sigma_1(s_2,x,y+\varphi(s_2,x),z_2),$
then, \\$\begin{array}{llll} &&\langle z_1-z_2,
z_1-z_2\rangle\\&=&\langle
D_{x_1}\varphi(s_1,x)\sigma_1(s_1,x,y+\varphi(s_1,x),z_1)-D_{x_1}\varphi(s_1,x)\sigma_1(s_1,x,y+\varphi(s_1,x),z_2),z_1-z_2\rangle\\
&&+\langle
D_{x_1}\varphi(s_1,x)\sigma_1(s_1,x,y+\varphi(s_1,x),z_2)-D_{x_1}\varphi(s_2,x)\sigma_1(s_2,x,y+\varphi(s_2,x),z_2),z_1-z_2\rangle\\
&\leq&\langle
D_{x_1}\varphi(s_1,x)\sigma_1(s_1,x,y+\varphi(s_1,x),z_2)-D_{x_1}\varphi(s_2,x)\sigma_1(s_2,x,y+\varphi(s_2,x),z_2),z_1-z_2\rangle\\
&\leq&C|s_1-s_2|(1+|x|+|y|)|z_1-z_2|+C(|\sigma_1(s_2,x,y+\varphi(s_1,x),z_2)-\sigma_1(s_1,x,y+\varphi(s_1,x),z_2)|\\
&&+|s_1-s_2|)|z_1-z_2|\\
&\leq&C|s_1-s_2|^2(1+|x|^2+|y|^2)+C(|\sigma_1(s_2,x,y+\varphi(s_1,x),z_2)-\sigma_1(s_1,x,y+\varphi(s_1,x),z_2)|^2\\
&&+|s_1-s_2|^2)+\frac{1}{2}|z_1-z_2|^2,\\
 \end{array}$\\
therefore, we have \\$\begin{array}{llll} |z_1-z_2|^2&\leq&
C|s_1-s_2|^2(1+|x|^2+|y|^2)+C|\sigma_1(s_2,x,y+\varphi(s_1,x),z_2)-\sigma_1(s_1,x,y+\varphi(s_1,x),z_2)|^2\\
&\leq& C|s_1-s_2|^2(1+|x|^2+|y|^2)+C|\rho(\delta)|^2,\end{array}$\\
where $\rho(\delta):=\sup\limits_{x\in \mathbb{R}^n,\ y\in \mathbb{R},\ z\in \mathbb{R}^d}|\sigma_1(s_2,x,y,z)-\sigma_1(s_1,x,y,z)|,\ s_1,\ s_2\in [0, T]$, for $|s_1-s_2|\leq \delta,$\ then from (B7) we have $\rho(\delta)\rightarrow 0,$\ as $\delta\rightarrow 0$.
It follows that $z=h(s,\bar{x},y,\zeta)$ is continuous with respect
to $s$.  \endpf

\bl For every $u\in \mathcal {U}_{t,t+\delta},$
 \be E[\int_t^{t+\delta}(|\overline{Y}_s^{1,u}|+|\overline{Z}_s^{1,u}|)ds\mid\mathcal
  {F}_t] \leq C\delta^\frac{5}{4},\ \  \mbox{P-}a.s.,\ 0\leq \delta\leq\bar{\delta}_1,\ee
  where the constant $C$ is independent of the control $u$ and of $\delta>0$.
\el  \noindent \textbf{Proof}. From equation (4.30), we have that
$\overline{Z}_s^{u}$ can be written as
$\overline{Z}_s^{u}=h(s,\overline{X}_s^u,\overline{Y}_s^{1,u},\overline{Z}_s^{1,u})$,
where $h$ satisfies the properties given in Proposition 4.1. Let
$F(s,x,y,z,u)=L(s,x,y,h(s,x,y,z),u)$. Then, (4.30) can be rewritten
as follows \be
  \left \{
  \begin{array}{llll}
  d\overline{Y}_s^{1,u} & =&-F(s,\overline{X}_s^u,\overline{Y}_s^{1,u},\overline{Z}_s^{1,u},u_s)ds+\overline{Z}_s^{1,u}dB_s, \ \ s\in[t,t+\delta], \\
\overline{Y}_{t+\delta}^{1,u}&=&0.
   \end{array}
   \right.
  \ee
Then the proof is similar to the proof of Lemma 4.1, we can get \be
|\overline{Y}_s^{1,u}|\leq
C\delta^\frac{1}{2}(1+|\overline{X}_s^{u}|),\ \ \mbox{P-a.s.},\ \
s\in [t,t+\delta].\ee On the other hand,
 \be |\overline{Z}_s^{1,u}|=|\overline{Z}_s^{u}-D\varphi(s,\overline{X}_s^u).\sigma(s,\overline{X}_s^u,\overline{Y}_s^{u},\overline{Z}_s^u)|\\
 \leq C(1+|\overline{X}_s^u|+|\overline{Y}_s^{u}|+|\overline{Z}_s^u|),\ \ \mbox{P-a.s.}\ee
And, from (4.27) and (4.35) we know, \be |\overline{Y}_s^{u}|\leq
C(1+|\overline{X}_s^{u}|),\ \ \mbox{P-a.s.}, \ s\in[t,t+\delta]. \ee
From (4.34), (4.35), (4.36), (4.37) and
(4.25),\be\begin{array}{llll}
&&|\overline{Y}_t^{1,u}|^2+E[\int_t^{t+\delta}|\overline{Z}_r^{1,u}|^2dr\mid\mathcal
{F}_t]=2E[\int_t^{t+\delta}\overline{Y}_r^{1,u}F(r,\overline{X}_r^{u},\overline{Y}_r^{1,u},\overline{Z}_r^{1,u},u_r)dr\mid\mathcal
{F}_t]\\
&\leq&CE[\int_t^{t+\delta}|\overline{Y}_r^{1,u}|(1+|\overline{X}_r^{u}|^2+|\overline{Y}_r^{1,u}|^2+|\overline{Z}_r^{1,u}|^2)dr\mid\mathcal
{F}_t]\\
&\leq&C\delta^\frac{1}{2}E[\int_t^{t+\delta}(1+|\overline{X}_r^{u}|^2+|\overline{X}_r^{u}|^3)dr\mid\mathcal
{F}_t]+C\delta^\frac{1}{2}E[\int_t^{t+\delta}(1+|\overline{X}_r^{u}|)|\overline{Z}_r^{u}|^2dr\mid\mathcal
{F}_t]\leq C\delta^\frac{3}{2},\ \ \mbox{P-a.s.}\\
\end{array}\ee
Therefore, \\$\begin{array}{llll}
&&E[\int_t^{t+\delta}(|\overline{Y}_s^{1,u}|+|\overline{Z}_s^{1,u}|)ds\mid\mathcal
{F}_t] \leq C\delta^\frac{1}{2}E[\int_t^{t+\delta}(1+|\overline{X}_r^{u}|)dr\mid\mathcal
{F}_t]+C\delta^\frac{1}{2}(E[\int_t^{t+\delta}|\overline{Z}_r^{1,u}|^2dr\mid\mathcal
{F}_t])^\frac{1}{2}\\
&\leq&C\delta^\frac{5}{4},\ \ \mbox{P-a.s.},\ 0\leq \delta
<\bar{\delta}_1.
\end{array}$\\   \endpf

\br From (4.36), (4.37)
and (4.25) we get \be
E[(\int_t^{t+\delta}|\overline{Z}_s^{1,u}|^2ds)^2\mid\mathcal
{F}_t]\leq C\delta^2,\ \ \mbox{P-a.s.} \ee \er
\br From Proposition 4.1 and the fact that the solution
$(\overline{Y}^{1,u},\overline{Z}^{1,u})$ belongs to $\mathcal
{S}^2(t, t + \delta;\mathbb{R})\times\mathcal {H}^2(t, t +
\delta;\mathbb{R}^d)$, we see that the unique solution $\hat{Z}_s^u$
of the equation
$$\hat{Z}_s^u=\overline{Z}_s^{1,u}+D\varphi(s,x).\sigma(s,x,\overline{Y}_s^{1,u}+\varphi(s,x),\hat{Z}_s^u),\
s\in[t,t+\delta],$$ belongs to  $\mathcal
{H}^2(t,t+\delta;\mathbb{R}^d)$ and
$\hat{Z}_s^u=h(s,x,\overline{Y}_s^{1,u},\overline{Z}_s^{1,u})$.
Similar to Remark 4.5, we know
\\$\begin{array}{ll}
{\rm{(i)}}\ \  E[\int_t^{t+\delta}|\hat{Z}_s^{u}|^2ds\mid\mathcal
{F}_t]\leq
C\delta,\ \ \mbox{P-a.s.};\ \
{\rm{(ii)}}\ \
E[(\int_t^{t+\delta}|\hat{Z}_s^{u}|^2ds)^2\mid\mathcal {F}_t]\leq
C\delta^2,\ \ \mbox{P-a.s.}
\end{array}$\\
Then, from Lemma 2.1 in~\cite{Pe4} BSDE (4.31) has a  unique solution
$(\overline{Y}^{2,u},\overline{Z}^{2,u}).$ \er

  \bl For every $u \in \mathcal {U}_{t,t+\delta}$, we have $$|\overline{Y}_t^{1,u}-\overline{Y}_t^{2,u}|\leq C\delta^{\frac{5}{4}},\ \ \mbox{P-a.s.},\ \ 0\leq\delta\leq\bar{\delta}_1,$$ where $C$ is
 independent of the control process $u$ and of $\delta>0$.  \el
  \noindent \textbf{Proof}. Similar to the proof of Lemma 4.2
we set $g(s)=L(s, \overline{X}_s^u,0,\overline{Z}_s^u,u_s)-L(s,
x,0,\overline{Z}_s^u,u_s),\
\rho_0(r)=(1+|x|^2+|\overline{Z}_s^u|^2)(r+r^2),\ r\geq 0.$
   Obviously, $|g(s)|\leq
   C\rho_0(|\overline{X}_s^u-x|),$ for $s\in [t,t+\delta],\ (t,x)\in [0,T)\times \mathbb{R}^n,\ u\in \mathcal
   {U}_{t,t+\delta}.$
  Therefore, we have, from equations (4.30) and (4.31), estimates
  (4.25), (4.35) and (4.39),
  \be\begin{array}{llll}
 & & |\overline{Y}_t^{1,u}-\overline{Y}_t^{2,u}|=|E[(\overline{Y}_t^{1,u}-\overline{Y}_t^{2,u})\mid\mathcal
  {F}_t]| \\
  &
  =&|E[\int_t^{t+\delta}(L(s,\overline{X}_s^{u},\overline{Y}_s^{1,u},\overline{Z}_s^{u},u_s)-L(s,x,0,\hat{Z}_s^{u},u_s))ds\mid\mathcal
  {F}_t]| \\
  & \leq &CE[\int_t^{t+\delta}(\rho_0(|\overline{X}_s^u-x|)+C(1+|\overline{X}_s^{u}|+|\overline{Y}_s^{1,u}|+|\overline{Z}_s^{u}|)|\overline{Y}_s^{1,u}|\\
  &&+C(1+|x|+|\overline{Z}_s^{u}|+|\overline{Z}_s^{u}-\hat{Z}_s^{u}|)|\overline{Z}_s^{u}-\hat{Z}_s^{u}|)ds\mid\mathcal
  {F}_t]\\
  & \leq &C\delta^\frac{5}{4}+CE[\int_t^{t+\delta}|\overline{Z}_s^{u}-\hat{Z}_s^{u}|ds\mid\mathcal
  {F}_t]+CE[\int_t^{t+\delta}|\overline{Z}_s^{u}||\overline{Z}_s^{u}-\hat{Z}_s^{u}|ds\mid\mathcal
  {F}_t]\\
  &&+CE[\int_t^{t+\delta}|\overline{Z}_s^{u}-\hat{Z}_s^{u}|^2ds\mid\mathcal
  {F}_t].\end{array}\ee
Furthermore,
 \\$\begin{array}{llll}
&&\langle\hat{Z}_s^{u}-\overline{Z}_s^{u},\hat{Z}_s^{u}-\overline{Z}_s^{u}\rangle\\
&=&\langle
D_{x_1}\varphi(s,x)\sigma_1(s,x,\overline{Y}_s^{1,u}+\varphi(s,x),\hat{Z}_s^{u})
-D_{x_1}\varphi(s,x)\sigma_1(s,x,\overline{Y}_s^{1,u}+\varphi(s,x),\overline{Z}_s^{u}),\hat{Z}_s^{u}-\overline{Z}_s^{u}\rangle\\
&&+\langle
D_{x_1}\varphi(s,x)\sigma_1(s,x,\overline{Y}_s^{1,u}+\varphi(s,x),\overline{Z}_s^{u})
-D_{x_1}\varphi(s,x)\sigma_1(s,\overline{X}_s^{u},\overline{Y}_s^{1,u}+\varphi(s,\overline{X}_s^{u}),\overline{Z}_s^{u}),\hat{Z}_s^{u}-\overline{Z}_s^{u}\rangle\\
&&+\langle
D_{x_1}\varphi(s,x)\sigma_1(s,\overline{X}_s^{u},\overline{Y}_s^{1,u}+\varphi(s,\overline{X}_s^{u}),\overline{Z}_s^{u})
-D_{x_1}\varphi(s,\overline{X}_s^{u})\sigma_1(s,\overline{X}_s^{u},\overline{Y}_s^{1,u}+\varphi(s,\overline{X}_s^{u}),\overline{Z}_s^{u}),
\hat{Z}_s^{u}-\overline{Z}_s^{u}\rangle\\
&\leq&-\beta_2|\hat{Z}_s^{u}-\overline{Z}_s^{u}|^2+C|\overline{X}_s^{u}-x||\hat{Z}_s^{u}-\overline{Z}_s^{u}|
+C|\overline{X}_s^{u}-x|(1+|\overline{X}_s^{u}|+|\overline{Y}_s^{1,u}|+|\overline{Z}_s^{u}|)|\hat{Z}_s^{u}-\overline{Z}_s^{u}|\\
&\leq&C|\overline{X}_s^{u}-x|^2+\frac{1}{2}|\hat{Z}_s^{u}-\overline{Z}_s^{u}|^2
+C|\overline{X}_s^{u}-x|^2(1+|\overline{X}_s^{u}|+|\overline{Y}_s^{1,u}|+|\overline{Z}_s^{u}|)^2,\\
\end{array}$\\
from Proposition 4.1 and Remark 4.4. Therefore, we have
 $$|\overline{Z}_s^{u}-\hat{Z}_s^{u}|\leq
C(1+|\overline{X}_s^{u}|)|\overline{X}_s^u-x|+C|\overline{X}_s^u-x|(|\overline{Y}_s^{1,u}|+|\overline{Z}_s^{u}|).$$
Then, from Lemma 4.5 and (4.25), the proof is complete. \endpf

We now consider the following equation \be
  \left \{
  \begin{array}{llll}
  d\overline{Y}_s^{3,u} & = &
  -L(s,x,0,\psi(s,x),u_s)ds+\overline{Z}_s^{3,u}dB_s,\ s\in
  [t,t+\delta],\\
\psi(s,x)&=&D\varphi(s,x).\sigma(s,x,\varphi(s,x),\psi(s,x)),\ s\in
  [t,t+\delta],\\
   \overline{Y}^{3,u}_{t+\delta}& = & 0.
   \end{array}
   \right.\ee
\bl For every $u \in \mathcal {U}_{t,t+\delta}$, we have
$$|\overline{Y}_t^{2,u}-\overline{Y}_t^{3,u}|\leq
C\delta^{\frac{5}{4}}, \ \mbox{P-a.s.},\ \
0\leq\delta\leq\bar{\delta}_1,$$ where $C$ is
 independent of the control process $u$ and of $\delta>0$.  \el
  \noindent \textbf{Proof}. From (4.31) and (4.41), we get \\$\begin{array}{llll}
 & & |\overline{Y}_t^{2,u}-\overline{Y}_t^{3,u}| = |E[\int_t^{t+\delta}(L(s,x,0,\hat{Z}_s^{u},u_s)-L(s,x,0,\psi(s,x),u_s))ds\mid\mathcal
  {F}_t]| \\
  & \leq &CE[\int_t^{t+\delta}(1+|x|+|\hat{Z}_s^{u}|)|\hat{Z}_s^{u}-\psi(s,x)|ds\mid\mathcal
  {F}_t].\end{array}$\\
From Remark 4.6,  $\hat{Z}_s^u=h(s, x, \overline{Y}_s^{1,u},
\overline{Z}_s^{1,u})$, and from Proposition 4.1,
$\psi(s,x)=h(s,x,0,0)$, hence we obtain
 $|\hat{Z}_s^{u}-\psi(s,x)|\leq
C(|\overline{Y}_s^{1,u}|+|\overline{Z}_s^{1,u}|).$\ Notice
 $$\begin{array}{llll}&&E[\int_t^{t+\delta}|\hat{Z}_s^{u}||\hat{Z}_s^{u}-\psi(s,x)|ds\mid\mathcal
  {F}_t]\\
  &\leq& C(E[\int_t^{t+\delta}|\hat{Z}_s^{u}|^2ds\mid\mathcal
  {F}_t])^\frac{1}{2}(E[\int_t^{t+\delta}(|\overline{Y}_s^{1,u}|+|\overline{Z}_s^{1,u}|)^2ds\mid\mathcal
  {F}_t])^\frac{1}{2} \leq  C\delta^\frac{5}{4},\ \ \mbox{P-a.s.},
\end{array}$$
where the last inequality is due to (4.38) and (i) in Remark 4.6. Furthermore, from Lemma 4.5,
we have $|\overline{Y}_t^{2,u}-\overline{Y}_t^{3,u}|\leq
C\delta^{\frac{5}{4}}, \ \mbox{P-a.s.}$ \endpf

\bl Let $\overline{Y}_0(\cdot)$ be the solution of the following
ordinary differential equation combined with an algebraic equation:
\be
  \left \{
  \begin{array}{llll}
  d\overline{Y}_0(s) & = &
  -L_0(s,x,0,\psi(s,x))ds,\ s\in
  [t,t+\delta],\\
\psi(s,x)&=&D\varphi(s,x).\sigma(s,x,\varphi(s,x),\psi(s,x)),\ s\in
  [t,t+\delta],\\
   \overline{Y}_0(t+\delta)& = & 0,
   \end{array}
   \right.\ee
   where the function $L_0$ is defined by $$L_0(s,x,0,z)=\sup\limits_{u\in U}L(s,x,0,z,u),\ (s,x,z)\in [t,t+\delta]\times\mathbb{R}^n\times\mathbb{R}^d.$$
Then, \mbox{P-a.s.},
$$\overline{Y}_0(t)=\esssup_{u\in \mathcal {U}_{t,t+\delta}}\overline{Y}_t^{3,u}.$$
   \el

 \noindent \textbf{Proof}. Since $L_0(s,x,0,z)=\sup\limits_{u\in
{U}}L(s,x,0,z,u),$ we have
$$L_0(s,x,0,\psi(s,x)) \geq L(s,x,0,\psi(s,x),u_s),\  s\in [t,t+\delta],\ \mbox{for all} \ u\in \mathcal{U}_{t,t+\delta},$$
we have $\overline{Y}_0(t)\geq \overline{Y}_t^{3,u},\
\mbox{P-a.s.},$ for any $u\in \mathcal{U}_{t,t+\delta}$. Indeed,
(4.42) can be regarded as a BSDE with the solution
$(Y_s,Z_s)=(\overline{Y}_0(s),0).$ This allows to apply the
comparison theorem of BSDEs.

 On the other hand, since
$L_0(s,x,0,z)=\sup\limits_{u\in U}L(s,x,0,z,u)$, there exists a
measurable function $\tilde{u}:[t,t+\delta]\times \mathbb{R}^n\times
\mathbb{R}^d\rightarrow U$, such that
$L_0(s,x,0,\psi(s,x))=L(s,x,0,\psi(s,x),\tilde{u}(s,x,\psi(s,x))),$
for all $s\in [t,t+\delta]$. We put
$\tilde{u}_s=\tilde{u}(s,x,\psi(s,x)), \ s\in[t,t+\delta]$,
obviously $\tilde{u}\in \mathcal {U}_{t,t+\delta},$ and
 $L_0(s,x,0,\psi(s,x))=L(s,x,0,\psi(s,x),\tilde{u}_s),\ \ s\in[t,t+\delta].$
Consequently, from the uniqueness of the solution of the BSDE, we
have
$(\overline{Y}_t^{3,\tilde{u}},\overline{Z}_t^{3,\tilde{u}})=(\overline{Y}_0(t),0)$,
 and particularly, $\overline{Y}_0(t)=\overline{Y}_t^{3,\tilde{u}},\
\mbox{P-a.s.}$ Therefore,
 $\overline{Y}_0(t)=\esssup_{u\in\mathcal {U}_{t,t+\delta}}\overline{Y}_t^{3,u}.$ \endpf
 We are now able to finish the proof of Theorem 4.2.
\\
Indeed, from (4.29) we know that
$$\esssup_{u\in \mathcal{U}_{t,t+\delta}}\overline{Y}_t^{1,u}\geq0, \ \mbox{P-a.s.}$$ Thus, from the Lemmas 4.6 and 4.7 we
get $\esssup_{u\in \mathcal{U}_{t,t+\delta}}\overline{Y}_t^{3,u} \geq
-C\delta^\frac{5}{4}, \ \mbox{P-a.s.}$\ Thus, by Lemma 4.8,
$\overline{Y}_0(t)\geq -C\delta^\frac{5}{4}, \ \mbox{P-a.s.},$ where
$\overline{Y}_0$ is the unique solution of (4.42). Consequently,
$$-C\delta^\frac{1}{4}\leq \frac{1}{\delta}\overline{Y}_0(t)=\frac{1}{\delta}\int_t^{t+\delta}L_0(s,x,0,\psi(s,x))ds,\ \delta>0,$$ from where, thanks to the continuity of $s\mapsto L_0(s,x,0,\psi(s,x))$,
it follows that
$$\sup\limits_{u\in U}L(t,x,0,\psi(t,x),u)=L_0(t,x,0,\psi(t,x))\geq0,$$where $ \psi(t,x)=D\varphi(t,x).\sigma(t,x,\varphi(t,x),\psi(t,x))$ and from
the definition of $L$ we see that $W$ is a viscosity subsolution of
(4.23). Similarly, we can prove that $W$ is a viscosity
supersolution of (4.23). Therefore, $W$ is a viscosity solution of
(4.23).\\
\endpf

\section{Examples}

 Now we give two examples associated with the two
cases studied above. For simplification, we set $m=n=d=1$,
and $G=1$. In the first example, $\sigma$ does not depend on $z$,
but depends on $u$.

 \bex \mbox{}We consider the following fully coupled FBSDE: \be\label{1000}
  \left \{
  \begin{array}{ll}
  dX_s^{t,x;u}  = (3X_s^{t,x;u}+5Z_s^{t,x;u})ds
         +(4X_s^{t,x;u}-5Y_s^{t,x;u}+u_s)dB_s, \\
  dY_s^{t,x;u} =  -(2X_s^{t,x;u}+3Y_s^{t,x;u}+4Z_s^{t,x;u}+u_s)ds + Z_s^{t,x;u}dB_s, \  \ \ s \in [t,T],\\
   X_t^{t,x;u} =  x,\
   Y_T^{t,x;u}  =  X_T^{t,x;u},
   \end{array}
   \right.
  \ee
  where $u\in \mathcal {U}$ is an admissible control.

For a given admissible control $u,$ the coefficients of equation
(\ref{1000}) satisfy the assumptions (B1), (B2) and (B4), then there exists a unique
solution $(X^{t,x;u}, Y^{t,x;u}, Z^{t,x;u}).$\ We define \be
W(t,x)=\esssup_{u\in\mathcal {U}_{t,T}}Y_t^{t,x;u}, \ee
it follows from Theorem 4.1 that $W(t,x)$ is the viscosity solution
of the following PDE: $$
  \left \{
  \begin{array}{llll}
  \frac{\partial}{\partial t}W(t,x)+\sup\limits_{u\in U}\{\frac{1}{2} \frac{\partial^2}{\partial x^2}W(t,x)(4x-5W(t,x)+u)^2+\frac{\partial}{\partial x}W(t,x)(3x+5\frac{\partial}{\partial x}W(t,x)(4x-5W(t,x)+u))\\
  \qquad\qquad\qquad\quad\ +2x+3W(t,x)+4\frac{\partial}{\partial x}W(t,x)(4x-5W(t,x)+u)\}=0,  \ \ (t,x)\in[0,T)\times\mathbb{R},\\
   W(T,x) =  x.
   \end{array}
   \right.
 $$
 \eex

In the following example, $\sigma$ depends on $z$, but does not depend on $u$.

\bex We consider the following fully coupled FBSDE: \be\label{1001}
  \left \{
  \begin{array}{ll}
  dX_s^{t,x;u} =(-{(X_s^{t,x;u})}^+-4Y_s^{t,x;u}+u_s)ds +
          (-X_s^{t,x;u}-L_\sigma Z_s^{t,x;u})dB_s, \\
  dY_s^{t,x;u}  =  -(2X_s^{t,x;u}-{(Y_s^{t,x;u})}^+-Z_s^{t,x;u}+u_s)ds +  Z_s^{t,x;u}dB_s,\ \ \ s \in [t,T], \\
   X_t^{t,x;u} =  x,\
   Y_T^{t,x;u} =  X_T^{t,x;u},
   \end{array}
   \right.
  \ee
  where the constant $L_\sigma>0$\ is sufficiently small, $u\in \mathcal {U}$ is an admissible control.

It is easy to check that the coefficients of equation (\ref{1001}) satisfy the assumptions (B1), (B2), (B4), (B5), (B6) and (B7), hence there
exists a unique solution $(X^{t,x;u},Y^{t,x;u},Z^{t,x;u}).$ We
define \be\label{1002} W(t,x)=\esssup_{u\in\mathcal {U}_{t,T}}Y_t^{t,x;u}, \ee
and we associate (\ref{1001}) with the following partial differential
equation, \be\label{1003}
  \left \{
  \begin{array}{llll}
  \frac{\partial}{\partial t}W(t,x)+\sup\limits_{u\in U}\{\frac{1}{2} \frac{\partial^2}{\partial x^2}W(t,x)(x+L_\sigma V(t,x))^2+\frac{\partial}{\partial x}W(t,x)(-4W(t,x)-x^++u)+2x-W^+(t,x)\\
  \qquad\qquad\qquad\quad\ -V(t,x)+u\}=0,\\
 V(t,x)  = \frac{\partial}{\partial x}W(t,x)(-x-L_\sigma V(t,x)), \ \ \ (t,x)\in[0,T)\times\mathbb{R},\\
   W(T,x) =  x.\
   \end{array}
   \right.
  \ee
Therefore, from Theorem 4.2, $W(t,x)$ defined by (\ref{1002}) is the
viscosity solution of (\ref{1003}). \eex

\section{\large{Appendix}}\label{App}
 \hskip1cm
In this subsection we prove some  basic important estimates for
fully coupled FBSDEs under monotonic assumptions, and present new
estimates and new generalized comparison theorem for FBSDEs on small time
interval. Let us now give four mappings:\\
  $b: \Omega \times [0,T] \times \mathbb{R}^n \times  \mathbb{R}
\times  \mathbb{R}^{d}
  \rightarrow \mathbb{R}^n, \ \ \ \ \ \ \ \ \ \  \sigma: \Omega \times [0,T] \times \mathbb{R}^n \times  \mathbb{R} \times  \mathbb{R}^{
  d}
  \rightarrow \mathbb{R}^{n\times d}, $\\
$f: \Omega \times [0,T] \times \mathbb{R}^n \times \mathbb{R}
\times  \mathbb{R}^{d}
  \rightarrow \mathbb{R},  \ \ \ \ \ \  \ \ \ \ \  \Phi: \Omega \times \mathbb{R}^n
  \rightarrow \mathbb{R}, $\\
$(b(t,x,y,z))_{t\in[0,T]},\ (\sigma(t,x,y,z))_{t\in[0,T]},\ (f(t,x,y,z))_{t\in[0,T]}$ are $\mathbb {F}$-progressively measurable for each $(x,y,z)\in \mathbb{R}^n \times  \mathbb{R} \times  \mathbb{R}^{d}$, and  $\Phi(x)$\ is ${\cal F}_T$-measurable for each $x\in \mathbb{R}^n$, which satisfy (B1) and (B2), and also \\
(C1) there exists a constant $K\geq 0$\ such that, for any $t\in [0,T],$\ $(x,y,z),\
(x',y',z')\in\mathbb{R}^n\times\mathbb{R}\times\mathbb{R}^d,\
\mbox{P-a.s.},$\\
$\begin{array}{llll}
&&|l(t,x,y,z)-l(t,x',y',z')|\leq K(|x-x'|+|y-y'|+|z-z'|),\ \mbox{where}\ l=b,\ \sigma,\ f, \mbox{respectively,\ and}\\
&&|\Phi(x)-\Phi(x')|\leq K(|x-x'|);
 \end{array}$\\
(C2) there exists a constant $L\geq 0$\ such that, for any $t\in
[0,T],$\ $(x,y,z)\in\mathbb{R}^n\times\mathbb{R}\times\mathbb{R}^d,\
\mbox{P-a.s.},$\\
$\begin{array}{llll}\mbox{    }\ \ \ \ \ |b(t,x,y,z)|+|\sigma(t,x,y,z)|+ |f(t,x,y,z)|+|\Phi(x)|\leq
L(1+|x|+|y|+|z|). \end{array}$\\

We consider the following FBSDE parameterized
by the initial condition $(t,\zeta) \in [0,T] \times
L^2(\Omega,\mathcal {F}_t,P;\mathbb{R}^n):$

\be\label{100}
  \left \{
  \begin{array}{ll}
  dX_s^{t,\zeta}  =  b(s,X_s^{t,\zeta},Y_s^{t,\zeta},Z_s^{t,\zeta})ds +
          \sigma(s,X_s^{t,\zeta},Y_s^{t,\zeta},Z_s^{t,\zeta}) dB_s, \\
  dY_s^{t,\zeta}  = -f(s,X_s^{t,\zeta},Y_s^{t,\zeta},Z_s^{t,\zeta})ds + Z_s^{t,\zeta}dB_s, \ \ \ \ \ s\in [t,T],\\
   X_t^{t,\zeta} =  \zeta,\
   Y_T^{t,\zeta}  =  \Phi(X_T^{t,\zeta}).
   \end{array}
   \right.
  \ee
\bp\label{App-Pro1} Under the assumptions (B1), (B2), (C1) and (C2),
for any $0 \leq t \leq T$ and the associated initial states $\zeta,\
\zeta' \in L^2(\Omega,\mathcal {F}_t,P;\mathbb{R}^n),$ we have the
following estimates, P-a.s.: $$
  \begin{array}{llll}
 {\rm(i)}&E[\sup\limits_{t\leq s\leq T}|{X}_s^{t,\zeta}-{X}_s^{t,\zeta'}|^2+\sup\limits_{t\leq s\leq T}|{Y}_s^{t,\zeta}-{Y}_s^{t,\zeta'}|^2 + \int_t^T|{Z}_s^{t,\zeta}-{Z}_s^{t,\zeta'}|^2ds\mid\mathcal {F}_t]  \leq  C|\zeta - \zeta'|^2; \\
   {\rm(ii)}&E[\sup\limits_{t\leq s\leq T}|X_s^{t,\zeta}|^2+\sup\limits_{t\leq s\leq T}|Y_s^{t,\zeta}|^2+\int_t^T|Z_s^{t,\zeta}|^2ds\mid\mathcal {F}_t]   \leq  C(1 +|\zeta|^2). \\
   \end{array}
 $$
If $\sigma$ also satisfies:

 {\rm(C3)} \ for any $t\in [0,T]$, for any
$(x,y,z)\in\mathbb{R}^n\times\mathbb{R}\times\mathbb{R}^d,\
\mbox{P-a.s.},$ $|\sigma(t,x,y,z)|\leq L(1+|x|+|y|),$ \\then we can
get

 {\rm(iii)}\ \ $E[\sup\limits_{t\leq s\leq
t+\delta}|X_s^{t,\zeta}-\zeta|^2\mid\mathcal {F}_t]\leq C\delta(1
+|\zeta|^2),\ \mbox{P-a.s.},\  0\leq\delta\leq T-t.$ \ep \noindent
\textbf{Proof}.  From Lemma 2.1 we know, for the initial states $\zeta,\
\zeta' \in L^2(\Omega,\mathcal {F}_t,P;\mathbb{R}^n),$ FBSDE
(\ref{100}) has a unique solution
$(X_s^{t,\zeta},Y_s^{t,\zeta},Z_s^{t,\zeta})_{s\in[t,T]}\in\mathcal
{S}^2(t,T;\mathbb{R}^n)\times\mathcal
{S}^2(t,T;\mathbb{R})\times\mathcal {H}^2(t,T;\mathbb{R}^d),$\ and
$(X_s^{t,\zeta'},Y_s^{t,\zeta'},Z_s^{t,\zeta'})_{s\in[t,T]}\in\mathcal
{S}^2(t,T;\mathbb{R}^n)\times\mathcal
{S}^2(t,T;\mathbb{R})\times\mathcal {H}^2(t,T;\mathbb{R}^d),$
respectively. We define $\hat{X}_s=X_s^{t,\zeta}-X_s^{t, {\zeta'}},
\ \hat{Y}_s=Y_s^{t,\zeta}-Y_s^{t,{\zeta'}},\
\hat{Z}_s=Z_s^{t,\zeta}-Z_s^{t,{\zeta'}},$ $\Delta
h(s)=h(s,X_s^{t,\zeta},Y_s^{t,\zeta},Z_s^{t,\zeta})-h(s,X_s^{t,{\zeta'}},Y_s^{t,{\zeta'}},Z_s^{t,{\zeta'}}),$
where $h=b,\ \sigma,\ f,\ A,$ respectively.\\
Applying It\^o's formula to $|\hat{X}_s|^2$ we have
$$\begin{array} [c]{llll}E[|\hat{X}_s|^2\mid\mathcal {F}_t]
&=& |\zeta-\zeta'|^2+E[\int_t^s(2\hat{X}_r\Delta
b(r)+|\Delta\sigma(r)|^2)dr\mid\mathcal {F}_t]  \\
&\leq
&|\zeta-\zeta'|^2+CE[\int_t^s(|\hat{X}_r|^2+|\hat{Y}_r|^2+|\hat{Z}_r|^2)dr\mid\mathcal
{F}_t],\  t\leq s\leq T.\\
\end{array} $$
 Then, from the Gronwall inequality, we obtain
\be\label{99} E[|\hat{X}_s|^2\mid\mathcal{F}_t] \leq
C(|\zeta-\zeta'|^2+E[\int_t^s(|\hat{Y}_r|^2+|\hat{Z}_r|^2)dr\mid\mathcal
{F}_t]),\ \mbox{P-a.s.},\ t\leq s\leq T. \ee
Apply It\^o's formula to $e^{\beta s}|\hat{Y}_s|^2$, taking $\beta$ large enough, using standard methods for BSDEs we can get\\
\be\label{98-1}\begin{array} [c]{llll}&& E[|\hat{Y}_s|^2\mid\mathcal
{F}_t]+E[\int_s^T|\hat{Y}_r|^2dr\mid\mathcal
{F}_t]+E[\int_s^T|\hat{Z}_r|^2dr\mid\mathcal
{F}_t] \\
&\leq &C(E[|\hat{X}_T|^2\mid\mathcal
{F}_t]+E[\int_t^T|\hat{X}_r|^2dr\mid\mathcal {F}_t]),\ t\leq s\leq
T.\end{array} \ee Then from (\ref{99}) it follows that
\be\label{98}\begin{array} [c]{llll}&& E[|\hat{Y}_s|^2\mid\mathcal
{F}_t]+E[\int_s^T|\hat{Y}_r|^2dr\mid\mathcal
{F}_t]+E[\int_s^T|\hat{Z}_r|^2dr\mid\mathcal
{F}_t] \\
&\leq&
C|\zeta-\zeta'|^2+CE[\int_t^T(|\hat{X}_r|^2+|\hat{Y}_r|^2+|\hat{Z}_r|^2)dr\mid\mathcal
{F}_t],\  t\leq s\leq T.\end{array}\ee
 On the other hand, applying It\^o's formula to $\langle
G\hat{X}_r,\hat{Y}_r\rangle$, from the assumption (B2) we get \\
\be\label{97}\begin{array} [c]{llll}&& \langle
G\hat{X}_s,\hat{Y}_s\rangle = E[\langle G\hat{X}_T,\hat{Y}_T\rangle\mid\mathcal {F}_s]-
E[\int_s^T\langle \Delta A(r),(\hat{X}_r,\hat{Y}_r,\hat{Z}_r)\rangle
dr\mid\mathcal {F}_s]\\
&\geq&E[\mu_1|G\hat{X}_T|^2\mid\mathcal
{F}_s]+E[\beta_1\int_s^T|G\hat{X}_r|^2dr\mid\mathcal {F}_s]
+E[\int_s^T\beta_2(|G^T\hat{Y}_r|^2+|G^T\hat{Z}_r|^2)dr\mid\mathcal
{F}_s].
\end{array} \ee
 Therefore, $ \langle
G\hat{X}_s,\hat{Y}_s\rangle\geq 0,\ t\leq s\leq T, \ \mbox{P-a.s.}$ \\
If $\beta_2>0$, then we get \be\label{96}\begin{array} [c]{llll}
\langle G\hat{X}_t,\hat{Y}_t\rangle &=&E[\langle
G\hat{X}_s,\hat{Y}_s\rangle\mid\mathcal {F}_t]- E[\int_t^s\langle
\Delta A(r),(\hat{X}_r,\hat{Y}_r,\hat{Z}_r)\rangle
dr\mid\mathcal {F}_t]\\
&\geq&
\beta_2E[\int_t^s(|G^T\hat{Y}_r|^2+|G^T\hat{Z}_r|^2)dr\mid\mathcal
{F}_t], \ t\leq s\leq T, \ \mbox{P-a.s.}
\end{array} \ee
Therefore, \be\label{95}
E[\int_t^s(|\hat{Y}_r|^2+|\hat{Z}_r|^2)dr\mid\mathcal {F}_t]\leq
C\langle G\hat{X}_t,\hat{Y}_t\rangle,\ t\leq s\leq T, \
\mbox{P-a.s.}\ee Then, from (\ref{99}) we can get \be\label{94}
E[|\hat{X}_s|^2\mid\mathcal{F}_t] \leq C|\zeta-\zeta'|^2+ C\langle
G\hat{X}_t,\hat{Y}_t\rangle, \ t\leq s\leq T, \ \mbox{P-a.s.}\ee
From (\ref{98}) we have \be\label{93} E[|\hat{Y}_s|^2\mid\mathcal
{F}_t]+E[\int_s^T(|\hat{Y}_r|^2+|\hat{Z}_r|^2)dr\mid\mathcal
{F}_t]\leq  C|\zeta-\zeta'|^2+ C\langle G\hat{X}_t,\hat{Y}_t\rangle,
\ t\leq s\leq T, \ \mbox{P-a.s.}\ee Therefore, recalling that
$\hat{X}_t=X_t^{t,\zeta}-X_t^{t, {\zeta'}}=\zeta-\zeta'$,
$$\begin{array} [c]{llll} |\hat{Y}_t|^2\leq C|\zeta-\zeta'|^2+
C\langle G\hat{X}_t,\hat{Y}_t\rangle\leq C|\zeta-\zeta'|^2+
C|\hat{X}_t||\hat{Y}_t|\leq C|\zeta-\zeta'|^2+ C|
\hat{X}_t|^2+\frac{1}{2}|\hat{Y}_t|^2,\ \mbox{P-a.s.}
\end{array}$$
which means $|\hat{Y}_t|\leq C|\zeta-\zeta'|,\ \mbox{P-a.s.}$\
Furthermore, from (\ref{94}), (\ref{93}), we can get $$
E[|\hat{X}_s|^2\mid\mathcal {F}_t]+E[|\hat{Y}_s|^2\mid\mathcal
{F}_t]+E[\int_s^T(|\hat{Y}_r|^2+|\hat{Z}_r|^2)dr\mid\mathcal
{F}_t]\leq C|\zeta-\zeta'|^2, \  t\leq s\leq T,\ \mbox{P-a.s.}$$ If
$\beta_2=0$, then according to assumption (B2), we have $\beta_1>0,\
\mu_1>0,\ m=n=1,$ i.e., $G\in \mathbb{R}$.
 From (\ref{97}), $E[|\hat{X}_T|^2\mid\mathcal
{F}_t]+E[\int_t^T|\hat{X}_r|^2dr\mid\mathcal {F}_t]\leq
CG\hat{X}_t\cdot\hat{Y}_t,\ C>0.$\ From (\ref{98-1}),
 $|\hat{Y}_t|^2+E[\int_t^T(|\hat{Y}_r|^2+|\hat{Z}_r|^2)dr\mid\mathcal
{F}_t]\leq CG\hat{X}_t\cdot\hat{Y}_t\leq
C|\zeta-\zeta'|^2+\frac{1}{2}|\hat{Y}_t|^2,$\ therefore,
 $|\hat{Y}_t|^2+E[\int_t^T(|\hat{Y}_r|^2+|\hat{Z}_r|^2)dr\mid\mathcal
{F}_t]\leq C|\zeta-\zeta'|^2.$ Furthermore, from (\ref{99}), $
E[|\hat{X}_s|^2\mid\mathcal {F}_t]\leq C|\zeta-\zeta'|^2,\ t\leq
s\leq T, \ \mbox{P-a.s.} $\ Then from (\ref{98-1}) we get
 $
E[|\hat{Y}_s|^2\mid\mathcal
{F}_t]+E[\int_s^T(|\hat{Y}_r|^2+|\hat{Z}_r|^2)dr\mid\mathcal
{F}_t]\leq C|\zeta-\zeta'|^2, \  t\leq s\leq T,\ \mbox{P-a.s.}$\
From above we always have \be\label{92}
E[|\hat{X}_s|^2\mid\mathcal {F}_t]+E[|\hat{Y}_s|^2\mid\mathcal
{F}_t]+E[\int_s^T(|\hat{Y}_r|^2+|\hat{Z}_r|^2)dr\mid\mathcal
{F}_t]\leq C|\zeta-\zeta'|^2, \  t\leq s\leq T,\ \mbox{P-a.s.}\ee
Finally, from equation (\ref{100}) and Buckholder-Davis-Gundy
inequality we have
\\$\begin{array} [c]{llll}E[\sup\limits_{t\leq s\leq
T}|\hat{X}_s|^2\mid\mathcal {F}_t]&\leq& 3|\zeta-\zeta'|^2+
CE[\int_t^T|\Delta b(r)|^2dr\mid\mathcal {F}_t]+ CE[\int_t^T|\Delta
\sigma(r)|^2dr\mid\mathcal
{F}_t] \\
&\leq&
3|\zeta-\zeta'|^2+CE[\int_t^T(|\hat{X}_r|^2+|\hat{Y}_r|^2+|\hat{Z}_r|^2)dr\mid\mathcal
{F}_t]\\
&\leq& C|\zeta-\zeta'|^2,\ \mbox{P-a.s.};
\end{array} $\\
and \\$\begin{array} [c]{llll}E[\sup\limits_{t\leq s\leq
T}|\hat{Y}_s|^2\mid\mathcal {F}_t]&\leq&
CE[|X_T^{t,\zeta}-X_T^{t,\zeta'}|^2\mid\mathcal
{F}_t]+CE[\int_t^T(|\hat{X}_r|^2+|\hat{Y}_r|^2+|\hat{Z}_r|^2)dr\mid\mathcal
{F}_t]\\
&\leq& C|\zeta-\zeta'|^2,\ \mbox{P-a.s.}
\end{array} $\\
Similarly we can prove (ii) by making full use of the monotonic assumption (B2).\\
Now it is not hard to prove (iii). Indeed, from (C3) we
obtain,\\$\begin{array} [c]{llll}&&E[\sup\limits_{t\leq s\leq
t+\delta}|X_s^{t,\zeta}-\zeta|^2\mid\mathcal {F}_t]\\
&\leq& 2E[|\int_t^{t+\delta}
b(r,X_r^{t,\zeta},Y_r^{t,\zeta},Z_r^{t,\zeta})dr|^2\mid\mathcal
{F}_t] + CE[\int_t^{t+\delta}
|\sigma(r,X_r^{t,\zeta},Y_r^{t,\zeta},Z_r^{t,\zeta})|^2dr\mid\mathcal
{F}_t]& \\
& \leq & C\delta
E[\int_t^{t+\delta}(1+|X_r^{t,\zeta}|^2+|Y_r^{t,\zeta}|^2+|Z_r^{t,\zeta}|^2)dr\mid\mathcal
{F}_t]+CE[\int_t^{t+\delta}(1+|X_r^{t,\zeta}|^2+|Y_r^{t,\zeta}|^2)dr\mid\mathcal
{F}_t]\\
 & \leq
&C\delta E[\sup\limits_{t\leq r\leq
t+\delta}(|X_r^{t,\zeta}|^2+|Y_r^{t,\zeta}|^2)+\int_t^{t+\delta}|Z_r^{t,\zeta}|^2dr\mid\mathcal
{F}_t]+C\delta+C\delta E[\sup\limits_{t\leq r\leq
t+\delta}(|X_r^{t,\zeta}|^2+|Y_r^{t,\zeta}|^2)\mid\mathcal {F}_t]\\
 & \leq
&C\delta(1+|\zeta|^2), \ \mbox{P-a.s.}
\end{array} $\\ \endpf

\br\label{App-Re1} From Proposition \ref{App-Pro1}, we have the
following inequalities: \be\label{91} |Y^{t,\zeta}_t| \leq
C(1+|\zeta|);\ \ \ \ \ \ \ |Y^{t,\zeta}_t-Y_t^{t,\zeta'}| \leq
C|\zeta-\zeta'|,\ \ \mbox{P-a.s.},\ee where the constant $C>0$
depends only on the Lipschitz constants and linear growth constants
of $b,\ \sigma,\ f \ and \ \Phi.$ \er

\br\label{App-Re2} Let $\Phi(x)=\xi\in L^2(\Omega,\mathcal {F}_T,P;\mathbb{R})$. From the proof of Proposition \ref{App-Pro1} we see that: \\
under the assumptions (B1), (B2)-(i), (C1) and (C2), the statements
of Proposition \ref{App-Pro1}-(i) and (ii) still hold true; if
furthermore, assumption (C3) holds, then we also have the same
result as Proposition \ref{App-Pro1}-(iii). \er

Now we consider the continuous dependence of the fully coupled FBSDE
on the terminal condition.
\bp\label{App-Pro2} Suppose the
assumptions (B1), (B2), (C1) and (C2) are satisfied. For any $0 \leq
t \leq T,$\ the associated initial state $\zeta \in
L^2(\Omega,\mathcal {F}_t,P;\mathbb{R}^n)$\ and $\xi\in
L^2(\Omega,\mathcal {F}_T,P;\mathbb{R})$, we let $(X^{t, \zeta}_s,
Y^{t, \zeta}_s, Z^{t, \zeta}_s)_{s\in [t, T]}$ be the solution of
FBSDE (\ref{100}) associated with $(b, \sigma, f, \zeta, \Phi)$, and
$(\overline{X}^{t, \zeta}_s, \overline{Y}^{t, \zeta}_s,
\overline{Z}^{t, \zeta}_s)_{s\in [t, T]}$ be the solution of FBSDE
(\ref{100}) associated with $(b, \sigma, f, \zeta, \Phi+\xi)$. Then
we have\ \
$$|Y^{t, \zeta}_t-\overline{Y}^{t, \zeta}_t|^2+E[\int_t^T| Y^{t, \zeta}_r-\overline{Y}^{t, \zeta}_r|^2dr\mid\mathcal
{F}_t]+E[\int_t^T|Z^{t, \zeta}_r-\overline{Z}^{t,
\zeta}_r|^2dr\mid\mathcal {F}_t]\leq CE[\xi^2\mid\mathcal {F}_t],\
\mbox{P-a.s.}$$ \ep \noindent The proof is similar to Proposition \ref{App-Pro1}, we omit it here.

 Let us now introduce the random field:
 $$u(t,x) = Y_s^{t,x}\mid_{s=t}, \ \ (t,x) \ \in \ [0,T]\times\mathbb{R}^n,$$
 where $Y^{t,x}$ is the solution of FBSDE (\ref{100}) with the initial state $x\in\mathbb{R}^n$.

As a consequence of Remark \ref{App-Re1} we have that, for all $t
\in [0,T],\ \mbox{P-a.s.},$ \be\label{79} \begin{array}[c]{llll}&&
{\rm(i)}\ |u(t,x)-u(t,y)|
\leq C|x-y|, \ \mbox{for} \ \mbox{all} \ x,y \in \mathbb{R}^n; \\
&&{\rm(ii)} \ |u(t,x)| \leq C(1+|x|), \ \mbox{for} \ \mbox{all} \ x\in
\mathbb{R}^n.\end{array}\ee

\br\label{App-Re3}  Moreover, it is well known that, under the
additional assumption that the functions \be\label{78} b,\sigma,f  \
\mbox{and} \ \Phi  \ are  \ deterministic, \ee $u$ is also a
deterministic function of $(t,x)$. \er The random field $u$ and
$Y^{t,\zeta}, (t,\zeta)\in[0,T]\times L^2(\Omega,\mathcal
{F}_t,P;\mathbb{R}^n),$ are related by the following theorem.

\bt\label{App-Th1} Under the assumptions (B1) and (B2), for any
$t\in[0,T]$ and $\zeta \in L^2(\Omega,\mathcal
{F}_t,P;\mathbb{R}^n),$ we have
$$u(t,\zeta)=Y_t^{t,\zeta}, \ \ \mbox{P-a.s.}$$ \et

The proof of Theorem \ref{App-Th1} is similar to that in
Peng~\cite{Pe4} for the decoupled FBSDE, or we can also refer to the
proof of Theorem A.2 in~\cite{BL}. \br\label{App-Re4} From Theorem
\ref{App-Th1}, we can obtain
$Y_s^{t,\zeta}=Y_s^{s,X_s^{t,\zeta}}=u(s,X_s^{t,\zeta}).$ \er

\bp\label{App-Pro3} Under the assumptions (B1), (B2), (C1), (C2) and
(C3), for any $p\geq 2$, $0 \leq t \leq T$ and the associated
initial state $\zeta \in L^p(\Omega,\mathcal {F}_t,P;\mathbb{R}^n),$
there exists $\tilde{\delta}_0>0$, which depends on $p$ and
Lipschitz constant $K$ and the linear growth constant $L$, such that
$$
  \begin{array}{llll}
{(\rm i)}&&E[\sup\limits_{t\leq s\leq t+\tilde{\delta}_0}|{X}_s^{t,\zeta}|^p+\sup\limits_{t\leq s\leq t+\tilde{\delta}_0}|{Y}_s^{t,\zeta}|^p + (\int_t^{t+\tilde{\delta}_0}|{Z}_s^{t,\zeta}|^2ds)^\frac{p}{2}\mid\mathcal {F}_t]  \leq  C_p(1+|\zeta |^p),\ \mbox{P-a.s.} ;\\
  {(\rm ii)}&& E[\sup\limits_{t\leq s\leq
t+\delta}|X_s^{t,\zeta}-\zeta|^p\mid\mathcal {F}_t]\leq
C_p\delta^\frac{p}{2}(1 +|\zeta|^p),\ \mbox{P-a.s.},\
0\leq\delta\leq \tilde{\delta}_0,
\end{array}
$$
where $(X^{t, \zeta}_s, Y^{t, \zeta}_s, Z^{t, \zeta}_s)_{s\in [t,
T]}$ is the solution of FBSDE (\ref{100}) associated with $(b,
\sigma, f, \zeta, \Phi)$. \ep \noindent \textbf{Proof}. From Remarks
\ref{App-Re1} and \ref{App-Re4}, we have \be\label{77}
|{Y}_s^{t,\zeta}|=|{Y}_s^{s,{X}_s^{t,\zeta}}|\leq
C(1+|{X}_s^{t,\zeta}|),\ \mbox{P-a.s.}\ee Since
${Y}_t^{t,\zeta}={Y}_s^{t,\zeta}+\int_t^s
f(r,{X}_r^{t,\zeta},{Y}_r^{t,\zeta},{Z}_r^{t,\zeta})dr-\int_t^s{Z}_r^{t,\zeta}dB_r,\
t\leq s\leq t+\delta,$\ from Buckholder-Davis-Gundy inequality, \\$
 \begin{array}{llll}
&&E[(\int_t^{t+\delta}|Z_s^{t,\zeta}|^2ds)^{\frac{p}{2}}\mid\mathcal
{F}_t]\leq C_pE[\sup\limits_{t\leq s\leq
t+\delta}|\int_t^s{Z}_r^{t,\zeta}dB_r|^p\mid\mathcal {F}_t]\\
&\leq&C_pE[\sup\limits_{t\leq s\leq
t+\delta}|{Y}_s^{t,\zeta}|^p+(\int_t^{t+\delta}
f(s,{X}_s^{t,\zeta},{Y}_s^{t,\zeta},{Z}_s^{t,\zeta})ds)^p\mid\mathcal
{F}_t]\\
&\leq&C_pE[\sup\limits_{t\leq s\leq
t+\delta}|{Y}_s^{t,\zeta}|^p\mid\mathcal
{F}_t]+C_pE[(\int_t^{t+\delta}(1+|{X}_s^{t,\zeta}|+|{Y}_s^{t,\zeta}|+|{Z}_s^{t,\zeta}|)ds)^p\mid\mathcal
{F}_t]\\
&\leq&C_pE[\sup\limits_{t\leq s\leq
t+\delta}|{Y}_s^{t,\zeta}|^p\mid\mathcal
{F}_t]+C_p\delta^p+C_p\delta^pE[\sup\limits_{t\leq s\leq
t+\delta}|{X}_s^{t,\zeta}|^p\mid\mathcal
{F}_t]+C_p\delta^pE[\sup\limits_{t\leq s\leq
t+\delta}|{Y}_s^{t,\zeta}|^p\mid\mathcal
{F}_t]\\
&&+C_p\delta^\frac{p}{2}E[(\int_t^{t+\delta}|Z^{t,\zeta}_s|^2ds)^{\frac{p}{2}}\mid\mathcal
{F}_t]\\
&= &C_p\delta^p+C_p\delta^pE[\sup\limits_{t\leq s\leq
t+\delta}|{X}_s^{t,\zeta}|^p\mid\mathcal
{F}_t]+(C_p+C_p\delta^p)E[\sup\limits_{t\leq s\leq
t+\delta}|{Y}_s^{t,\zeta}|^p\mid\mathcal
{F}_t]\\
&&+C_p\delta^\frac{p}{2}E[(\int_t^{t+\delta}|Z^{t,\zeta}_s|^2ds)^{\frac{p}{2}}\mid\mathcal
{F}_t],
\end{array}
  $\\
there exists $\delta_0>0,$ such that $1-C_p\delta_0^\frac{p}{2}>0,$
then we get, for any $0\leq\delta\leq\delta_0$, P-a.s.,\be\label{76}
E[(\int_t^{t+\delta}|Z_s^{t,\zeta}|^2ds)^{\frac{p}{2}}\mid\mathcal
{F}_t]\leq C_p\delta^{p}+C_p\delta^{p}E[\sup\limits_{t\leq s\leq
t+\delta}|{X}_s^{t,\zeta}|^p\mid\mathcal
{F}_t]+(C_p+C_p\delta^p)E[\sup\limits_{t\leq s\leq
t+\delta}|{Y}_s^{t,\zeta}|^p\mid\mathcal {F}_t].\ee On the other
hand, for $t\leq s\leq T$, from (\ref{100}) and (\ref{77}),
$$
 \begin{array}{llll}&&
E[\sup\limits_{t\leq r\leq
s}|X_r^{t,\zeta}-\zeta|^p\mid\mathcal {F}_t]\\
&\leq &C_pE[(\int_t^s
|b(r,X_r^{t,\zeta},Y_r^{t,\zeta},Z_r^{t,\zeta})|dr)^p\mid\mathcal
{F}_t]+C_pE[(\int_t^s|
\sigma(r,X_r^{t,\zeta},Y_r^{t,\zeta},Z_r^{t,\zeta})|^2dr)^\frac{p}{2}\mid\mathcal
{F}_t]\\
&\leq
&C_pE[(\int_t^s(1+|X_r^{t,\zeta}-\zeta|+|\zeta|+|Z_r^{t,\zeta}|)dr)^p\mid\mathcal
{F}_t]+C_pE[(\int_t^s(1+|X_r^{t,\zeta}|)^2dr)^\frac{p}{2}\mid\mathcal
{F}_t]\\
&\leq
&C_p(1+|\zeta|^p)(s-t)^\frac{p}{2}+C_p(s-t)^\frac{p}{2}E[(\int_t^{s}|Z_r^{t,\zeta}|^2dr)^{\frac{p}{2}}\mid\mathcal
{F}_t]+C_pE[\int_t^s|X_r^{t,\zeta}-\zeta|^pdr\mid\mathcal {F}_t],
\end{array}
  $$
from Gronwall inequality, \be\label{75} E[\sup\limits_{t\leq r\leq
s}|X_r^{t,\zeta}-\zeta|^p\mid\mathcal {F}_t]\leq
C_p(1+|\zeta|^p)(s-t)^\frac{p}{2}+C_p(s-t)^\frac{p}{2}E[(\int_t^{s}|Z_r^{t,\zeta}|^2dr)^{\frac{p}{2}}\mid\mathcal
{F}_t],\ \mbox{P-a.s.},\ t\leq s\leq T.\ee Then, from (\ref{76}),
(\ref{77}) and (\ref{75})  we have
 $$\begin{array}{llll}&&
E[(\int_t^{t+\delta}|Z_s^{t,\zeta}|^2ds)^{\frac{p}{2}}\mid\mathcal
{F}_t]\\
&\leq& C_p\delta^p(1+|\zeta|^p)+ C_p|\zeta|^p+C_p+(C_p+
C_p\delta^p)E[\sup\limits_{t\leq s\leq
t+\delta}|X_s^{t,\zeta}-\zeta|^p\mid\mathcal {F}_t]\\
&\leq& C_p+C_p|\zeta|^p+C_p\delta^\frac{p}{2}(1+|\zeta|^p)+(C_p+
C_p\delta^p)C_p\delta^\frac{p}{2}E[(\int_t^{s}|Z_r^{t,\zeta}|^2dr)^{\frac{p}{2}}\mid\mathcal
{F}_t],
\end{array}$$
taking $0<\widetilde{\delta}_0\leq \delta_0,$ such that $1-(C_p+
C_p\widetilde{\delta}_0^p)C_p\widetilde{\delta}_0^\frac{p}{2}>0$,
then $0\leq \delta\leq \widetilde{\delta}_0,$ $$
E[(\int_t^{t+\delta}|Z_s^{t,\zeta}|^2ds)^{\frac{p}{2}}\mid\mathcal
{F}_t]\leq  C_p(1+|\zeta|^p),\ \mbox{P-a.s.} $$ From (\ref{75}), we
get \
$E[\sup\limits_{t\leq s\leq
t+\delta}|X_s^{t,\zeta}-\zeta|^p\mid\mathcal {F}_t]\leq
C_p\delta^\frac{p}{2}(1 +|\zeta|^p),\ \mbox{P-a.s.},\
0\leq\delta\leq \widetilde{\delta}_0.$\\ From (\ref{77})  we have
\ $E[\sup\limits_{t\leq s\leq t+\delta}|Y_s^{t,\zeta}|^p\mid\mathcal
{F}_t]\leq C_p(1 +|\zeta|^p),\ \mbox{P-a.s.},\  0\leq\delta\leq
\widetilde{\delta}_0.$ \endpf

\bp\label{App-Pro4} We suppose the assumptions (C1), (C2) and (C4) hold true, where the assumption (C4) is the following hypothesis:\\
(C4) the Lipschitz constant $L_\sigma\geq0$\ of $\sigma$ with
respect to $z$\ is sufficiently small, i.e., there exists small
enough $L_\sigma\geq0$\ such that, for all $t\in[0, T],\ u\in U,\
x_1, x_2\in\mathbb{R}^n,\ y_1, y_2\in\mathbb{R},\ z_1,
z_2\in\mathbb{R}^d,$\ P-a.s.,
$$|\sigma(t,x_1,y_1,z_1,u)-\sigma(t,x_2,y_2,z_2,u)|\leq K(|x_1- x_2|+|y_1-y_2|)+L_\sigma|z_1-z_2|.$$
Then, there exists a constant $0< \delta_0$, only depending on the
Lipschitz constant $K$,\  such that for every $0\leq\delta\leq
\delta_0$ and $\zeta\in L^2(\Omega,\mathcal {F}_t,P;\mathbb{R}^n),$
FBSDE (\ref{100}) has a
 unique solution $(X^{t, \zeta}_s, Y^{t, \zeta}_s, Z^{t, \zeta}_s)_{s\in [t, t+\delta]}$\ on the time interval $[t,\ t+\delta]$.
\ep

\noindent \textbf{Proof}. Let us give any $0<T\leq T_0$, and observe that for any pair $v=(y, z)\in {\cal{H}}^{2}(t,T;{\mathbb{R}}^{1+d})$\ there exists a unique solution $V=(Y, Z)\in {\cal{H}}^{2}(t,T;{\mathbb{R}}^{1+d})$\ to the following decoupled FBSDE:
\be\label{101}
  \left \{
  \begin{array}{llll}
  dX_s   =  b(s,X_s ,y_s ,z_s)ds +
          \sigma(s,X_s ,y_s ,z_s) dB_s, \\
  dY_s  =  -f(s,X_s ,Y_s ,Z_s)ds + Z_sdB_s, \ \ \ \ \ s\in [t,T],\\
   X_t   =  \zeta,\
   Y_T   =  \Phi(X_T).   \end{array}
   \right.
  \ee
We are going to prove that there exists a constant $0< \delta_0$, only depending on the
Lipschitz constant $K$,\  such that for
every $0\leq\delta\leq \delta_0$ the mapping defined by
$$I(v)=V: {\cal{H}}^{2}(t,t+\delta;{\mathbb{R}}^{1+d})\rightarrow {\cal{H}}^{2}(t,t+\delta;{\mathbb{R}}^{1+d})$$
is a contraction.

Let $v^i=(y^i, z^i)\in {\cal{H}}^{2}(t,t+\delta;{\mathbb{R}}^{1+d})$,\ and $V^i=I(v^i),\ i=1, 2.$\ We define $\widehat{v}=(y^1-y^2, z^1-z^2),$\ and $\widehat{V}=(Y^1-Y^2, Z^1-Z^2),\ \widehat{X}=X^1-X^2$. Then, we get
\be\label{102}\begin{array}{rcl}
   E[\sup\limits_{t\leq s\leq r}|\widehat{X}_s|^2\mid\mathcal {F}_t]&\leq& 2E[(\int_t^r|b(s, X_s^1, y_s^1, z_s^1)-b(s, X_s^2, y_s^2, z_s^2)|ds)^2\mid\mathcal {F}_t]\\
   & &+8E[\int_t^r|\sigma(s, X_s^1, y_s^1, z_s^1)-\sigma(s, X_s^2, y_s^2, z_s^2)|^2ds\mid\mathcal {F}_t]\\
 & \leq & 6(r-t)K^2E[\int_t^r(|\widehat{X}_s|^2+|\widehat{y}_s|^2+|\widehat{z}_s|^2)ds\mid\mathcal {F}_t]\\
 & &+24E[\int_t^r(K^2|\widehat{X}_s|^2+K^2|\widehat{y}_s|^2+L_\sigma^2|\widehat{z}_s|^2)ds\mid\mathcal {F}_t],
  \end{array} \ee
and the Gronwall inequality yields \be \label{103}\begin{array}{rcl}
& &E[\sup\limits_{t\leq s\leq T}|\widehat{X}_s|^2\mid\mathcal
{F}_t]\leq CE[\int_t^T|\widehat{y}_s|^2ds\mid\mathcal {F}_t]+
(C(T-t)+C L_\sigma^2)E[\int_t^T|\widehat{z}_s|^2ds\mid\mathcal {F}_t]\\
& & \leq C(T-t)E[\sup\limits_{t\leq s\leq
T}|\widehat{y}_s|^2\mid\mathcal {F}_t]+
(C(T-t)+C L_\sigma^2)E[\int_t^T|\widehat{z}_s|^2ds\mid\mathcal
{F}_t].
 \end{array}\ee
On the other hand, by using BSDE standard estimates combined with the help of (\ref{103}), we get
$$ \begin{array}{rcl}
& &E[\sup\limits_{t\leq s\leq T}|\widehat{Y}_s|^2+\int_t^T|\widehat{Z}_s|^2ds]\\
&\leq & CE[|\Phi(X_T^1)-\Phi(X_T^2)|^2]+
CE[\int_t^T|f(r, X_r^1, Y_r^1, Z_r^1)-f(r, X_r^2, Y_r^1, Z_r^1)|^2dr]\\
&\leq & CE[|\widehat{X}_T|^2]+CE[\int_t^T|\widehat{X}_r|^2dr]\\
&\leq & C(T-t)E[\sup\limits_{t\leq s\leq T}|\widehat{y}_s|^2]+
(C(T-t)+CL_\sigma^2)E[\int_t^T|\widehat{z}_s|^2ds]\\
&\leq & (C(T-t)+CL_\sigma^2)(E[\sup\limits_{t\leq s\leq
T}|\widehat{y}_s|^2]+E[\int_t^T|\widehat{z}_s|^2ds]).
 \end{array}$$
Notice $L_\sigma$\ is sufficiently small such that $CL_\sigma^2\leq
\frac{1}{3}$. Then there exists $\delta_0>0$\ such that
$C\delta_0+CL_\sigma^2< \frac{1}{2}$, and therefore, for any $0\leq
\delta\leq \delta_0$, we have \be \label{104}\begin{array}{rcl}
E[\sup\limits_{t\leq s\leq
t+\delta}|\widehat{Y}_s|^2+\int_t^{t+\delta}|\widehat{Z}_s|^2ds]\leq
\frac{1}{2}(E[\sup\limits_{t\leq s\leq
t+\delta}|\widehat{y}_s|^2]+E[\int_t^{t+\delta}|\widehat{z}_s|^2ds]).
 \end{array}\ee
It follows immediately that for any $0\leq \delta\leq \delta_0$\
this mapping $I$ has a unique fixed point $I(V)=V$,  i.e., FBSDE
(\ref{100}) has a unique solution $(X^{t, \zeta}_s, Y^{t, \zeta}_s,
Z^{t, \zeta}_s)_{s\in [t, t+\delta]}$\ on the time interval $[t,\
t+\delta]$.\endpf

\br\label{App-Re5} In fact, from the proof we see that $L_\sigma\geq
0$\ with $CL_\sigma^2<1$\ is sufficient for Proposition
\ref{App-Pro4}. \er

\bt\label{App-Th2} (Generalized Comparison Theorem) We suppose the
assumptions (C1), (C2) and (C4) are satisfied. Let $\delta_0>0$\ be
a constant, only depending on the Lipschitz constant $K$,\  such
that for every $0\leq\delta\leq \delta_0$ and $\zeta\in
L^2(\Omega,\mathcal {F}_t,P;\mathbb{R}^n),$\ FBSDE (\ref{100}) has a
 unique solution $(X^i_s, Y^i_s, Z^i_s)_{s\in [t, t+\delta]}$\ associated with $(b, \sigma, f, \zeta, \Phi^i)$\ on the time interval $[t,\ t+\delta]$, respectively.
Then, if for any $0\leq\delta\leq \delta_0$\ we have $\Phi^1(X^2_{t+\delta})\geq \Phi^2(X^2_{t+\delta})$, P-a.s., (resp., $\Phi^1(X^1_{t+\delta})\geq \Phi^2(X^1_{t+\delta})$, P-a.s.), we also have $Y_t^1\geq Y_t^2$, P-a.s.
\et

\noindent \textbf{Proof}. The proof is similar to that of Theorem 3.1 in Wu~\cite{W}. For notational simplification, we assume $d=n=1$. We define $\widehat{X}=X^1-X^2,\ \widehat{Y}=Y^1-Y^2,\ \widehat{Z}=Z^1-Z^2$, then $(\widehat{X}, \widehat{Y}, \widehat{Z})$\ satisfies the following FBSDE:
\be \label{105}\begin{array}{rcl}
& &d\widehat{X}_s=(b_s^1\widehat{X}_s+b_s^2\widehat{Y}_s+b_s^3\widehat{Z}_s)ds
+(\sigma_s^1\widehat{X}_s+\sigma_s^2\widehat{Y}_s+\sigma_s^3\widehat{Z}_s)dB_s,\\
& & d\widehat{Y}_s=-(f_s^1\widehat{X}_s+f_s^2\widehat{Y}_s+f_s^3\widehat{Z}_s)ds+\widehat{Z}_sdB_s,\\
& &\widehat{X}_t=0,\ \widehat{Y}_{t+\delta}=\overline{\Phi}\widehat{X}_{t+\delta}+\Phi^1(X^2_{t+\delta})- \Phi^2(X^2_{t+\delta}),
 \end{array}\ee
where
$$
  \begin{array}{llll}
    & l_s^1=\frac{l(s,X^1_s ,Y^1_s ,Z^1_s) -
          l(s,X^2_s ,Y^1_s ,Z^1_s)}{X^1_s-X_s^2},\ \ \ \ \ \ \hfill \widehat{X}_s\neq 0; \\
   &l_s^1=0,\ \ \ \ \ \ \hfill \mbox{otherwise};
   \end{array}
$$
$$
  \begin{array}{llll}
    & l_s^2=\frac{l(s,X^2_s ,Y^1_s ,Z^1_s) -
          l(s,X^2_s ,Y^2_s ,Z^1_s)}{Y^1_s-Y_s^2},\ \ \ \ \ \ \hfill \widehat{Y}_s\neq 0; \\
   &l_s^2=0,\ \ \ \ \ \ \hfill \mbox{otherwise};
   \end{array}
$$
$$
  \begin{array}{llll}
    & l_s^3=\frac{l(s,X^2_s ,Y^2_s ,Z^1_s) -
          l(s,X^2_s ,Y^2_s ,Z^2_s)}{Z^1_s-Z_s^2},\ \ \ \ \ \ \hfill \widehat{Z}_s\neq 0; \\
   &l_s^3=0,\ \ \ \ \ \ \hfill \mbox{otherwise},
   \end{array}
$$
$l=b,\ \sigma,\ f$\ respectively, and
$$
  \begin{array}{llll}
    & \overline{\Phi}=\frac{\Phi^1(X^1_{t+\delta}) -
          \Phi^1(X^2_{t+\delta})}{X^1_{t+\delta}-X_{t+\delta}^2},\ \ \ \ \ \ \hfill \widehat{X}_{t+\delta}\neq 0; \\
   &\overline{\Phi}=0,\ \ \ \ \ \ \hfill \mbox{otherwise}.
   \end{array}
$$
It's easy to check that (\ref{105}) satisfies (C1), (C2) and (C4).
Therefore, from Proposition \ref{App-Pro4} there exists a constant
$0< \delta_1\leq \delta_0$, such that for every $0\leq\delta\leq
\delta_1$,\ (\ref{105}) has a unique solution on $[t, t+\delta]$,
i.e., $(\widehat{X}, \widehat{Y}, \widehat{Z})$\ is the unique
solution of (\ref{105}) on $[t, t+\delta]$, for every
$0\leq\delta\leq \delta_1$. Now we introduce the dual FBSDE \be
\label{106}\begin{array}{rcl} &
&dP_s=(f_s^2P_s-b_s^2Q_s-\sigma_s^2K_s)ds
+(f_s^3P_s-b_s^3Q_s-\sigma_s^3K_s)dB_s,\\
& & dQ_s=(f_s^1P_s-b_s^1Q_s-\sigma_s^1K_s)ds+K_sdB_s,\\
& &P_t=1,\ Q_{t+\delta}=-\overline{\Phi}P_{t+\delta}.
 \end{array}\ee
Similarly, (\ref{106}) satisfies (C1), (C2) and (C4). Consequently,
due to Proposition \ref{App-Pro4} there exists a constant $0<
\delta_2\leq \delta_1$, such that for every $0\leq\delta\leq
\delta_2$,\ (\ref{106}) has a unique solution $(P, Q, K)$\ on $[t,
t+\delta]$. Using It\^o's formula to
$\widehat{X}_sQ_s+\widehat{Y}_sP_s$, we deduce from the equations
(\ref{105}) and (\ref{106}) that,
$$ E[\widehat{X}_{t+\delta}(-\overline{\Phi}P_{t+\delta})|{\cal F}_t]+E[(\overline{\Phi}\widehat{X}_{t+\delta}+\Phi^1(X^2_{t+\delta})- \Phi^2(X^2_{t+\delta}))P_{t+\delta}|{\cal F}_t] =\widehat{Y}_t,$$
i.e.,
\be\label{107} \widehat{Y}_t=E[(\Phi^1(X^2_{t+\delta})- \Phi^2(X^2_{t+\delta}))P_{t+\delta}|{\cal F}_t].\ee
Since $\Phi^1(X^2_{t+\delta})\geq \Phi^2(X^2_{t+\delta})$, P-a.s., if we can prove $P_{t+\delta}\geq 0$, P-a.s., then we can get $\widehat{Y}_t\geq 0$, P-a.s.

For this we define the following stopping time: $\tau=\inf\{s>t: P_s=0\}\wedge (t+\delta),$\
and consider the following FBSDE (\ref{108}) on $[\tau, t+\delta]$ (notice that $\tau>t$, since $P$\ is continuous and $P_t=1$):
\be \label{108}\begin{array}{rcl}
& &d\widetilde{P}_s=(f_s^2\widetilde{P}_s-b_s^2\widetilde{Q}_s-\sigma_s^2\widetilde{K}_s)ds
+(f_s^3\widetilde{P}_s-b_s^3\widetilde{Q}_s-\sigma_s^3\widetilde{K}_s)dB_s,\\
& & d\widetilde{Q}_s=(f_s^1\widetilde{P}_s-b_s^1\widetilde{Q}_s-\sigma_s^1\widetilde{K}_s)ds+\widetilde{K}_sdB_s,\\
& &\widetilde{P}_\tau=0,\ \widetilde{Q}_{t+\delta}=-\overline{\Phi}\widetilde{P}_{t+\delta}.
 \end{array}\ee
Similarly to the equation (\ref{106}) we see that, (\ref{108})
satisfies (C1), (C2) and (C4), and therefore, from Proposition
\ref{App-Pro4} there exists $ 0<\delta_3\leq \delta_2$\ such that
for every $0\leq\delta\leq \delta_3$,\ (\ref{108}) has a unique
solution $(\widetilde{P}, \widetilde{Q}, \widetilde{K})$\ on $[\tau,
t+\delta]$. Clearly, $(\widetilde{P}_s, \widetilde{Q}_s,
\widetilde{K}_s)\equiv (0, 0, 0)$\ is the unique solution of
(\ref{108}). Let
$$\begin{array}{rcl}
& &\overline{P}_s=I_{[t, \tau]}(s)P_s+I_{(\tau, t+\delta]}(s)\widetilde{P}_s,\\
& & \overline{Q}_s=I_{[t, \tau]}(s)Q_s+I_{(\tau, t+\delta]}(s)\widetilde{Q}_s,\\
& &\overline{K}_s=I_{[t, \tau]}(s)K_s+I_{(\tau, t+\delta]}(s)\widetilde{K}_s,\ s\in [t, t+\delta].
\end{array}$$
It's easy to show that $(\overline{P}, \overline{Q}, \overline{K})$\ is a solution of FBSDE (\ref{106}). Therefore, from the uniqueness of solution of FBSDE (\ref{106}) on $[t, t+\delta]$, where $0\leq\delta\leq \delta_3$, we have
$\overline{P}_t=P_t=1>0$. Furthermore, from the definition of $\tau$\ we have $\overline{P}_{t+\delta}\geq 0$, P-a.s., that is, ${P}_{t+\delta}\geq 0$, P-a.s.
 Therefore, we have $Y_t^1\geq Y_t^2$, P-a.s. \endpf

\bp\label{App-Pro5} Let $\Phi$ be deterministic. We suppose the
assumptions (C1), (C2), and (C4) hold true. Then, for every $p\geq
2$, there exists sufficiently small constant $\widetilde{\delta}>0$,
only depending on the Lipschitz constant $K$, and some constant
$\tilde{C}_{p.K}$, only depending on $p$, the Lipschitz constant $K$
and the linear growth constant $L$, such that for every $0\leq
\delta\leq \widetilde{\delta}$ and $\zeta\in L^p(\Omega,\mathcal
{F}_t,P;\mathbb{R}^n),$
$$
  \begin{array}{llll}
{(\rm i)}&&E[\sup\limits_{t\leq s\leq t+\delta}|{X}_s^{t,\zeta}|^p+\sup\limits_{t\leq s\leq t+\delta}|{Y}_s^{t,\zeta}|^p + (\int_t^{t+\delta}|{Z}_s^{t,\zeta}|^2ds)^\frac{p}{2}\mid\mathcal {F}_t]  \leq  \tilde{C}_{p,K}(1+|\zeta |^p),\ \mbox{P-a.s.}; \\
{(\rm ii)}&& E[\sup\limits_{t\leq s\leq
t+\delta}|X_s^{t,\zeta}-\zeta|^p\mid\mathcal {F}_t]\leq
\tilde{C}_{p,K}\delta^\frac{p}{2}(1 +|\zeta|^p),\ \mbox{P-a.s.};\\
{(\rm
iii)}&&E[(\int_t^{t+\delta}|{Z}_s^{t,\zeta}|^2ds)^\frac{p}{2}\mid\mathcal
{F}_t]\leq \tilde{C}_{p,K}\delta^\frac{p}{2}(1 +|\zeta|^p),\
\mbox{P-a.s.},
\end{array}
$$
where $(X^{t, \zeta}_s, Y^{t, \zeta}_s, Z^{t, \zeta}_s)_{s\in [t, t+\delta]}$ is the solution of FBSDE (\ref{100})  associated with $(b, \sigma, f, \zeta, \Phi)$\ and with the time horizon $t+\delta$.\ep
\noindent \textbf{Proof}. Due to Proposition \ref{App-Pro4}
there exists a constant $ \delta_0>0$, such that for every
$0\leq\delta\leq \delta_0$,\ (\ref{100}) has a unique solution on
$[t, t+\delta]$, i.e.,
\be\label{74}{Y}_s^{t,\zeta}=\Phi({X}_{t+\delta}^{t,\zeta})+\int_s^{t+\delta}
f(r,{X}_r^{t,\zeta},{Y}_r^{t,\zeta},{Z}_r^{t,\zeta})dr-\int_s^{t+\delta}{Z}_r^{t,\zeta}dB_r,\
t\leq s\leq {t+\delta}.\ee We consider
$\widetilde{{Y}}_s^{t,\zeta}={Y}_s^{t,\zeta}-\Phi(\zeta)$, and for
any $\beta\geq 0$\ by applying It\^o's formula to $e^{\beta
s}|\widetilde{{Y}}_s^{t,\zeta}|^2$, we get \be\label{73}
\begin{array}{rcl}
&&E[e^{\beta s}|\widetilde{{Y}}_s^{t,\zeta}|^2|{\cal F}_s]+E[\int_s^{t+\delta}\{ \beta e^{\beta r}|\widetilde{{Y}}_r^{t,\zeta}|^2+e^{\beta r}| {{Z}}_r^{t,\zeta}|^2\}dr|{\cal F}_s]\\
&=&E[e^{\beta (t+\delta)}|\widetilde{{Y}}_{t+\delta}^{t,\zeta}|^2|{\cal F}_s]+E[\int_s^{t+\delta} e^{\beta r}2\widetilde{{Y}}_r^{t,\zeta}f(r,X^{t, \zeta}_r, Y^{t, \zeta}_r, Z^{t, \zeta}_r)dr|{\cal F}_s]\\
&=&E[e^{\beta ({t+\delta})}|\widetilde{{Y}}_{t+\delta}^{t,\zeta}|^2|{\cal F}_s]+E[\int_s^{t+\delta} e^{\beta r}2\widetilde{{Y}}_r^{t,\zeta}(f(r,X^{t, \zeta}_r, Y^{t, \zeta}_r, Z^{t, \zeta}_r)-f(r, \zeta, \Phi(\zeta), 0)+
f(r, \zeta, \Phi(\zeta), 0))dr|{\cal F}_s].\\
\end{array}
\ee
By taking $\beta$\ large enough and since $|f(r, \zeta, \Phi(\zeta), 0)|\leq C(1+|\zeta|)$, we get by using BSDE standard methods
\be\label{72}
\begin{array}{llll}
&& |\widetilde{{Y}}_s^{t,\zeta}|^2 +E[\int_s^{t+\delta}\{ |\widetilde{{Y}}_r^{t,\zeta}|^2+| {{Z}}_r^{t,\zeta}|^2\}dr|{\cal F}_s]\\
&&\leq CE[\sup\limits_{s\leq r\leq {t+\delta}}|X^{t, \zeta}_r-\zeta|^2|{\cal F}_s]+C({t+\delta}-s)(1+|\zeta|^2),\ \mbox{P-a.s.,}\\
\end{array}
\ee where $C$\ only depends on $K$\ and $L$. Therefore, from
(\ref{74}) and (\ref{72}) and Buckholder-Davis-Gundy inequality,
\be\label{71} E[\sup_{t\leq s\leq
{t+\delta}}|\widetilde{{Y}}_s^{t,\zeta}|^2|{\cal F}_t]\leq
CE[\sup_{t\leq r\leq {t+\delta}}|X^{t, \zeta}_r-\zeta|^2|{\cal
F}_t]+C\delta(1+|\zeta|^2), \ \mbox{P-a.s.} \ee On the other hand,
from (\ref{72}) \be\label{70} |\widetilde{{Y}}_s^{t,\zeta}|^2\leq
CE[\sup\limits_{t\leq r\leq {t+\delta}}|X^{t,
\zeta}_r-\zeta|^2|{\cal F}_s]+C\delta(1+|\zeta|^2),\ \mbox{P-a.s.},\
t\leq s\leq {t+\delta}. \ee When $p>2$, we define
$\eta=\sup\limits_{t\leq r\leq {t+\delta}}|X^{t,
\zeta}_r-\zeta|\in L^2(\Omega, {\cal F}_{t+\delta}, P)$. Then
$M_s:=E[\eta|{\cal F}_s],\ s\in[t, {t+\delta}]$, is a martingale,
and from Doob's martingale inequality, we have
\be\label{69}\begin{array}{llll}
&&E[\sup\limits_{t\leq s\leq {t+\delta}}|M_s|^{\frac{p}{2}}|{\cal F}_t]\leq C_p E[|M_{{t+\delta}}|^{\frac{p}{2}}|{\cal F}_t]\leq C_p E[\eta^{\frac{p}{2}}|{\cal F}_t]\\
&&=C_pE[\sup\limits_{t\leq r\leq {t+\delta}}|X^{t,
\zeta}_r-\zeta|^p|{\cal F}_t], \ \mbox{P-a.s.}
\end{array}
\ee
Therefore, from (\ref{70}) and (\ref{69})
\be\label{68} E[\sup_{t\leq s\leq {t+\delta}}|\widetilde{{Y}}_s^{t,\zeta}|^p|{\cal F}_t]\leq C_pE[\sup_{t\leq r\leq {t+\delta}}|X^{t, \zeta}_r-\zeta|^p|{\cal F}_t]+C_p\delta^{\frac{p}{2}}(1+|\zeta|^p), \ \mbox{P-a.s.}
\ee
Now we consider ${Y}_s^{t,\zeta}-\Phi(\zeta)=\Phi({X}_{t+\delta}^{t,\zeta})-\Phi(\zeta)+\int_s^{t+\delta}
f(r,{X}_r^{t,\zeta},{Y}_r^{t,\zeta},{Z}_r^{t,\zeta})dr-\int_s^{t+\delta}{Z}_r^{t,\zeta}dB_r,\
t\leq s\leq {t+\delta}.$
From Burkholder-Davis-Gundy inequality and (\ref{68}),
$$
 \begin{array}{llll}
&&E[(\int_t^{{t+\delta}}|Z_s^{t,\zeta}|^2ds)^{\frac{p}{2}}\mid\mathcal
{F}_t]\leq C_pE[\sup\limits_{t\leq s\leq
{t+\delta}}|\int_t^s{Z}_r^{t,\zeta}dB_r|^p\mid\mathcal {F}_t]\\
&\leq&C_pE[\sup\limits_{t\leq s\leq
{t+\delta}}|\widetilde{{Y}}_s^{t,\zeta}|^p+(\int_t^{{t+\delta}}
|f(s,{X}_s^{t,\zeta},{Y}_s^{t,\zeta},{Z}_s^{t,\zeta})|ds)^p\mid\mathcal
{F}_t]\\
&=&C_pE[\sup\limits_{t\leq s\leq
t+\delta}|\widetilde{{Y}}_s^{t,\zeta}|^p\mid\mathcal
{F}_t]\\
&&+C_pE[(\int_t^{t+\delta}
|f(s,{X}_s^{t,\zeta},{Y}_s^{t,\zeta},{Z}_s^{t,\zeta})-f(s, \zeta,
\Phi(\zeta),0)+f(s, \zeta, \Phi(\zeta),0)|ds)^p\mid\mathcal
{F}_t]\\
&\leq&(C_p+C_p\delta^p)E[\sup\limits_{t\leq s\leq
{t+\delta}}|{X}_s^{t,\zeta}-\zeta|^p\mid\mathcal
{F}_t]+C_p\delta^\frac{p}{2}(1+|\zeta|^p)+C_p\delta^\frac{p}{2}E[(\int_t^{t+\delta}|Z_s^{t,\zeta}|^2ds)^{\frac{p}{2}}\mid\mathcal
{F}_t].\\
\end{array}
  $$
By choosing then $0< {\delta}_1\leq \delta_0$\ such that
$1-C_p{\delta}_1^\frac{p}{2}>0$, we get, for any $0\leq \delta\leq
{\delta}_1$, P-a.s.,\be\label{67}
E[(\int_t^{{t+\delta}}|Z_s^{t,\zeta}|^2ds)^{\frac{p}{2}}\mid\mathcal
{F}_t]\leq(C_p+C_p\delta^p)E[\sup\limits_{t\leq s\leq
T}|{X}_s^{t,\zeta}-\zeta|^p\mid\mathcal
{F}_t]+C_p\delta^\frac{p}{2}(1+|\zeta|^p).\ee Similarly, equation
(\ref{100}) and the estimates (\ref{68}) and (\ref{67}) yield $$
 \begin{array}{llll}&&
E[\sup\limits_{t\leq r\leq
{t+\delta}}|X_r^{t,\zeta}-\zeta|^p\mid\mathcal {F}_t]\\
&\leq &C_pE[(\int_t^{t+\delta}
b(r,X_r^{t,\zeta},Y_r^{t,\zeta},Z_r^{t,\zeta})dr)^p\mid\mathcal
{F}_t]+C_pE[(\int_t^{t+\delta}|
\sigma(r,X_r^{t,\zeta},Y_r^{t,\zeta},Z_r^{t,\zeta})|^2dr)^\frac{p}{2}\mid\mathcal
{F}_t]\\
&\leq
&C_p(1+|\zeta|^p)\delta^\frac{p}{2}+(C_p\delta^\frac{p}{2}+C_pL_\sigma^p)E[(\int_t^{t+\delta}|Z_r^{t,\zeta}|^2dr)^{\frac{p}{2}}\mid\mathcal
{F}_t]\\
& &\ \ +C_p\delta^\frac{p}{2}E[\sup\limits_{t\leq r\leq
{t+\delta}}|X_r^{t,\zeta}-\zeta|^pdr\mid\mathcal
{F}_t]\\
&\leq
&C_p(1+|\zeta|^p)\delta^\frac{p}{2}+(C_p\delta^\frac{p}{2}+C_pL_\sigma^p)
E[\sup\limits_{t\leq r\leq
{t+\delta}}|X_r^{t,\zeta}-\zeta|^p\mid\mathcal {F}_t].
\end{array}
  $$
Let $L_\sigma>0$ be sufficiently small such that $C_pL_\sigma^p<1$.
Then there exists constant $0< \delta_2\leq \delta_1$\ such that
$1-(C_p\delta_2^\frac{p}{2}+C_pL_\sigma^p)>0,$ and we obtain,
\be\label{66} E[\sup\limits_{t\leq r\leq
{t+\delta}}|X_r^{t,\zeta}-\zeta|^p\mid\mathcal {F}_t]\leq
C_p(1+|\zeta|^p)\delta^\frac{p}{2},\ \mbox{P-a.s.},\ t\leq s\leq
{t+\delta}.\ee Finally, (\ref{68}) and (\ref{67}) allow to complete
the proof.\endpf
Similarly, we can prove the following proposition.
\bp\label{App-Pro6} Suppose that $(b_i,\sigma_i,f_i,\Phi_i),\ i=1,2$
all satisfy the assumptions (C1), (C2) and (C4). There exists a
constant $0< \delta_0$, only depending on the Lipschitz constant
$K$,\  such that for every $0\leq\delta\leq \delta_0$, the same
initial state $\zeta \in L^2(\Omega,\mathcal {F}_t,P;\mathbb{R}^n)$,
$(X^i_s, Y^i_s, Z^i_s)_{s\in [t, t+\delta]}$ is the solution of
FBSDE (\ref{100}) associated with $(b_i,\sigma_i,f_i,\Phi_i)$ on the
time interval $[t, t+\delta]$, respectively. Then we have: there
exists a constant $\delta_1>0,$ such that for every
$0\leq\delta\leq\delta_1,$
$$\begin{array}{llll}|Y^1_t-Y^2_t|^2&\leq& CE[|\Phi_1(t+\delta,X_{t+\delta}^1)-\Phi_2(t+\delta,X_{t+\delta}^1)|^2\mid\mathcal
{F}_t]\\
&&+C\delta
E[\int_t^{t+\delta}|(b_1-b_2)(s,X_s^1,Y_s^1,Z_s^1)|^2ds\mid\mathcal
{F}_t]\\
&&+C
E[\int_t^{t+\delta}|(\sigma_1-\sigma_2)(s,X_s^1,Y_s^1,Z_s^1)|^2ds\mid\mathcal
{F}_t]\\
&&+C\delta
E[\int_t^{t+\delta}|(f_1-f_2)(s,X_s^1,Y_s^1,Z_s^1)|^2ds\mid\mathcal
{F}_t],\ \mbox{P-a.s.}\end{array}$$ \ep

\br\label{App-Re6} When $(b_1,\sigma_1,f_1)=(b_2,\sigma_2,f_2)$ in
Proposition \ref{App-Pro6}, we have
$$|Y^1_t-Y^2_t|\leq C(E[|\Phi_1(t+\delta,X^1_{t+\delta})-\Phi_2(t+\delta,X^1_{t+\delta})|^2|{\cal
F}_t])^{\frac{1}{2}},\ \mbox{P-a.s.}$$ \er \bc\label{App-Cro1} Under
the assumptions (C1), (C2) and (C4), there exists a constant $0<
\delta_0$, only depending on the Lipschitz constant $K$,\  such that
for every $0\leq\delta\leq \delta_0$, the associated initial state
$\zeta \in L^2(\Omega,\mathcal {F}_t,P;\mathbb{R}^n)$\ and any
$\varepsilon>0$, $(X^{t, \zeta}_s, Y^{t, \zeta}_s, Z^{t,
\zeta}_s)_{s\in [t, t+\delta]}$ is the solution of FBSDE (\ref{100})
associated with $(b, \sigma, f, \zeta, \Phi)$,  and
$(\overline{X}^{t, \zeta}_s, \overline{Y}^{t, \zeta}_s,
\overline{Z}^{t, \zeta}_s)_{s\in [t, t+\delta]}$ is the solution of
FBSDE (\ref{100}) associated with $(b, \sigma, f, \zeta,
\Phi+\varepsilon)$\ on the time interval $[t, t+\delta]$,
respectively. Then we have\ \
$$|Y^{t, \zeta}_t-\overline{Y}^{t, \zeta}_t|\leq C\varepsilon,\ \mbox{P-a.s.}$$
\ec

\section*{Acknowledgments}  The authors thank the referees for their very helpful comments. Juan Li acknowledges financial support by the NSF of P.R.China (No. 10701050, 11071144, 11222110), Shandong
Province (Nos. Q2007A04, BS2011SF010, JQ201202), SRF for ROCS (SEM), 111 Project (No. B12023).

\end{document}